\documentclass[11pt]{article}%
\usepackage{multicol}
\usepackage{amsmath}
\usepackage{amsfonts}
\usepackage{float}
\usepackage{amssymb}
\usepackage{graphicx}
\usepackage{placeins}
\usepackage{xcolor}
\usepackage[letterpaper,margin=1.25in]{geometry}
\usepackage{url}
\usepackage{hyperref}
\usepackage{setspace}%
\providecommand{\U}[1]{\protect\rule{.1in}{.1in}}
\hypersetup{colorlinks=true,linkcolor=black,citecolor=black,urlcolor=black,filecolor=black,pdfborder={0 0 0}}
\newtheorem{theorem}{Theorem}
\newtheorem{assumption}{Assumption}

\newtheorem{corollary}{Corollary}

\newtheorem{lemma}{Lemma}

\newtheorem{proposition}{Proposition}
\newtheorem{remark}[theorem]{Remark}

\newtheorem{thm}{Theorem}

\newenvironment{proof}[1][Proof]{\noindent \textbf{#1.} }{\  \rule{0.5em}{0.5em}}
\begin{document}

\title{A High-Dimensional LM Test for Cross-Sectional Independence with Serially
Correlated Data}
\date{}
\author{Ulrich Hounyo\\Department of Economics\\University at Albany, SUNY\\Albany, NY 12222, USA
\and Chihwa Kao\\Department of Economics\\University of Connecticut\\Storrs, CT 06269, USA
\and Min Seong Kim\\Department of Economics\\University of Connecticut\\Storrs, CT 06269, USA}
\maketitle

\begin{abstract}
We develop a high-dimensional Lagrange multiplier (LM) test for
cross-sectional independence in Gaussian panels with a common stationary,
short-memory temporal covariance matrix. The statistic is based on
$\operatorname{Tr}(R^{2})$, where $R$ is the sample correlation matrix. When
the cross-sectional and time dimensions grow proportionally, we derive
explicit expansions for the null mean and variance in terms of temporal
spectral moments. Random diagonal normalization produces a non-negligible
centering correction, while squared sample correlations sharing a
cross-sectional unit contribute an additional variance term. Both corrections
vanish under temporal independence. We establish a central limit theorem under
the null, consistency against fixed alternatives, and a Gaussian-shift limit
with a detection boundary under local alternatives. Monte Carlo experiments
examine the size and power of the oracle test and feasible versions that
estimate the temporal spectral moments. The asymptotic theory applies to the
oracle statistic; feasible procedures are evaluated numerically.

\end{abstract}

\thispagestyle{plain}



\newpage

\noindent\textbf{Running head:} HIGH-DIMENSIONAL LM TEST

\medskip\noindent\textbf{Corresponding author and proofs:} Chihwa Kao,
Department of Economics, University of Connecticut, Storrs, CT 06269, USA;
e-mail: chih-hwa.kao@uconn.edu.

\bigskip

\newpage

\setcounter{equation}{0}

\section{Introduction}

\label{sec:intro}

The LM test of Breusch and Pagan (1980) aggregates squared sample correlations
and has a chi-squared null limit when the cross-sectional dimension is fixed
and the time dimension grows. In high dimensions, Schott (2005) derives a
Gaussian limit for the sum of squared correlations when both dimensions
increase proportionally under temporally independent Gaussian sampling.
Baltagi, Feng, and Kao (2012) extend high-dimensional LM testing to
fixed-effects panel regressions, and Huang, Li, and Yao (2025) develop unified
and robust LM-type tests for large panels; both assume serially independent errors.

The existing literature does not establish the centering and scaling of the
quadratic LM statistic under serial correlation when the cross-sectional and
time dimensions grow proportionally. Baltagi, Kao, and Peng (2016) document
size distortions of existing cross-sectional correlation tests under serial
dependence and develop a modified CD test; Hounyo, Kao, and Kim (2025) also
study CD-type statistics. Those procedures aggregate correlations linearly,
whereas the LM statistic aggregates their squares. Chen and Liu (2018)
consider a separable Gaussian model but test independence in the temporal
dimension. Thus, these results do not directly provide a high-dimensional null
distribution for the quadratic LM statistic with serially correlated observations.

We study an $n\times T$ Gaussian panel whose cross-sectional units are
independent under the null and share a stationary, short-memory Toeplitz
temporal covariance matrix $\Gamma$, normalized by $Tr(\Gamma)/T=1$. Write
$R=(r_{ij})$ for the sample correlation matrix formed using the maintained
zero population means. Let $\lambda=T/n$ denote the dimension ratio, with
$\lambda\to\lambda_{0}\in(0,\infty)$. We consider the statistic
\[
Tr(R^{2})=n+2\sum_{i<j}r_{ij}^{2}.
\]
Although this is a linear spectral statistic (LSS), existing CLTs for
high-dimensional sample correlation matrices (Gao, Han, Pan, and Yang, 2017;
Yin, Zheng, and Zou, 2023; Li, Pan, Xie, and Zhou, 2024) address temporally
independent sampling and cannot be applied directly to our setting.

The contribution of this paper is to show how to center and scale $Tr(R^{2})$
under serial correlation. Let $\tau_{k,T}$ denote the normalized temporal
spectral moments of $\Gamma$, and let $K_{T}$ be the function of $\tau
_{2,T},\tau_{3,T},\tau_{4,T}$ defined in
Section~\ref{subsec:model-specification}. We show that
\begin{align*}
\mathbb{E}\{Tr(R^{2})\}  &  =n+\frac{n-1}{\lambda}\tau_{2,T} +\frac
{4(n-1)}{n\lambda^{2}}(\tau_{2,T}^{2}-\tau_{3,T})+O(n^{-1}),\\
Var\{Tr(R^{2})\}  &  =\frac{4\tau_{2,T}^{2}}{\lambda^{2}} +\frac{8K_{T}%
}{\lambda^{3}}+O(n^{-1}),
\end{align*}
and that $Tr(R^{2})$, centered and scaled accordingly, is asymptotically
standard normal. The centering and scaling depend on the temporal dependence
only through $\tau_{2,T}$, $\tau_{3,T}$, and $\tau_{4,T}$, and reduce to the
known-mean version of Schott (2005) when $\Gamma=I_{T}$. Whereas the LSS
literature relies on random matrix theory, we obtain these moments directly by
Wick's theorem. This is the key step: once the centering and scaling are
known, asymptotic normality is established by a martingale central limit
theorem. Details are given in Section~\ref{sec:main}. The null CLT concerns
the oracle statistic, which uses the population spectral moments $\tau_{k,T}$;
feasible versions based on estimated moments are evaluated in the Monte Carlo
experiments, which also illustrate the size distortions that arise when serial
correlation is ignored.

We also characterize the behavior of the test under alternatives. We establish
consistency against fixed alternatives and a Gaussian shift under local
alternatives, with aggregate signal $\Vert R_{\ast}-I_{n}\Vert_{F}^{2}$, where
$R_{\ast}$ is the population cross-sectional correlation matrix. This yields
an explicit local-power expression and a detection boundary for the statistic
under the alternatives considered here; it is not asserted to be a minimax
boundary over all tests.

Section~\ref{sec:model} presents the model and statistic,
Section~\ref{sec:main} establishes the null expansions and CLT,
Section~\ref{sec:power} studies power, and Section~\ref{sec:mc} reports
simulations. Section~\ref{sec:conclusion} concludes. Proofs of the null
results (Lemmas~2 and~3 and Theorem~1) are in the Appendix; the remaining
proofs, including those of the power results, are in the Online Supplement.

\section{Model and Test Statistic}

\label{sec:model}

\subsection{Model Specification}

\label{subsec:model-specification}

Let $X=(x_{it})\in\mathbb{R}^{n\times T}$ denote observations on $n$
cross-sectional units over $T$ periods. Assume the separable Gaussian model
\[
X=U^{1/2}Z\Gamma^{1/2},
\]
where $U\in\mathbb{R}^{n\times n}$ is a positive-definite cross-sectional
covariance matrix, $\Gamma\in\mathbb{R}^{T\times T}$ is a temporal covariance
matrix, and $Z$ has independent $N(0,1)$ entries. We use symmetric matrix
square roots. Then $vec(X)\sim N(0,\Gamma\otimes U)$ under column-wise
vectorization, and $Cov(x_{it},x_{js})=U_{ij}\Gamma_{ts}$. We write
$\lambda=T/n$ and impose proportional growth in Assumption~3.

The temporal covariance matrix $\Gamma=\Gamma_{T}=(\gamma_{|t-s|})_{t,s=1}%
^{T}$, with $\gamma_{0}=1$, is stationary and short-memory; the precise
conditions are given in Assumption~2. We suppress the subscript $T$ except
where the dependence on $T$ matters. Define the normalized temporal spectral
moments by
\[
\tau_{k,T}=\frac{1}{T}Tr(\Gamma^{k})=\frac{1}{T}\sum_{t=1}^{T}\nu_{t,T}%
^{k},\qquad k=1,2,\ldots,
\]
where $\nu_{1,T},\ldots,\nu_{T,T}$ are the eigenvalues of $\Gamma$. Since
$\gamma_{0}=1$, $\tau_{1,T}=1$. The normalization $\gamma_{0}=1$ is without
loss of generality for the sample correlation matrix, which is invariant to
positive scalar rescaling of $\Gamma$. Define
\[
K_{T}=\tau_{4,T}-2\tau_{2,T}\tau_{3,T}+\tau_{2,T}^{3},\qquad K=\lim
_{T\rightarrow\infty}K_{T}=\tau_{4}-2\tau_{2}\tau_{3}+\tau_{2}^{3},
\]
where $\tau_{k}=\lim_{T\rightarrow\infty}\tau_{k,T}$.

\subsection{Hypotheses and Test Statistic}

Let $x_{i}=(x_{i1},\ldots,x_{iT})^{\prime}\in\mathbb{R}^{T}$ be the $i$th row
of $X$ written as a column vector. We test cross-sectional independence. Under
Gaussianity, this is equivalent to $U$ being diagonal, or to the population
correlation matrix $R_{\ast}=\operatorname{diag}(U)^{-1/2}U\operatorname{diag}%
(U)^{-1/2}$ being $I_{n}$. Because the sample correlation matrix is invariant
to positive row-wise rescaling, we use the normalization $U=I_{n}$ when
deriving the null distribution. The null hypothesis is therefore
$H_{0}:R_{\ast}=I_{n}$, with $U=I_{n}$ as its normalized representation.

With the population mean specified as zero, define the uncentered sample
covariance and correlation matrices by
\[
S=\frac{1}{T}XX^{\prime},\qquad R=\operatorname{diag}(S)^{-1/2}%
S\operatorname{diag}(S)^{-1/2}.
\]
In particular,
\[
r_{ij}=\frac{\sum_{t=1}^{T}x_{it}x_{jt}}{\left(  \sum_{t=1}^{T}x_{it}%
^{2}\right)  ^{1/2}\left(  \sum_{t=1}^{T}x_{jt}^{2}\right)  ^{1/2}},\qquad
Tr(R^{2})=n+2\sum_{i<j}r_{ij}^{2}.
\]
This identity expresses the quadratic spectral statistic as a sum of squared
off-diagonal sample correlations, the basis of the LM-type test.

Under $H_{0}$ and $\Gamma=I_{T}$, the known-mean version of the standardized
statistic of Schott (2005) is
\[
LM_{S}=\frac{\sum_{i<j}r_{ij}^{2}-\frac{n(n-1)}{2T}}{\sqrt{\frac
{n(n-1)(T-1)}{T^{2}(T+2)}}}=\frac{Tr(R^{2})-n-\frac{n(n-1)}{T}}{2\sqrt
{\frac{n(n-1)(T-1)}{T^{2}(T+2)}}}.
\]

\subsection{Assumptions}

The null results below use Assumptions~1--3. The analysis of alternatives
additionally invokes Assumption~4.

\begin{assumption}
The entries $Z_{it}$ are independent and identically distributed $N(0,1)$
random variables.
\end{assumption}

\begin{assumption}
(Stationary short-memory temporal covariance) For each $T$, $\Gamma
_{T}=(\gamma_{|t-s|})_{t,s=1}^{T}$ is a deterministic symmetric positive
semidefinite Toeplitz matrix, with $\gamma_{0}=1$ and $\sum_{h=0}^{\infty
}|\gamma_{h}|<\infty$. In particular, $T^{-1}Tr(\Gamma_{T})=1$ and $\sup
_{T}\|\Gamma_{T}\|_{op}<\infty$.
\end{assumption}

\begin{assumption}
(Proportional growth) As $n,T\to\infty$, $T/n\to\lambda_{0}\in(0,\infty)$.
\end{assumption}

Gaussianity permits exact higher-order moment calculations using Wick's
theorem (Isserlis' formula). The short-memory assumption controls $\tau_{k,T}$
for each fixed $k$, while proportional growth allows $n$ and $T$ to increase together.

\begin{remark}
[Limits of the temporal spectral moments]\label{rem:szego}Under Assumption~2,
$\Vert\Gamma_{T}\Vert_{op}\leq\sum_{h=-\infty}^{\infty}|\gamma_{|h|}|$, so
$\tau_{k,T}=O(1)$ for each fixed $k$, and the spectral symbol $f(\omega
)=\sum_{h=-\infty}^{\infty}\gamma_{|h|}e^{i\,h\omega}$ is continuous and
bounded. Szeg\H{o}'s theorem (Grenander \& Szeg\H{o}, 1958, Chap.~5) then
gives $\tau_{k}=\lim_{T\to\infty}\tau_{k,T}=(2\pi)^{-1}\int_{-\pi}^{\pi
}f(\omega)^{k}d\omega$ for each fixed $k$. Assumption~2 imposes no rate on
this convergence. Accordingly, terms with an $O(n)$ coefficient retain
$\tau_{k,T}$, while terms with an $O(1)$ coefficient may use $\tau_{k}$ at an
$o(1)$ cost. If in addition $\sum_{h\geq1}h|\gamma_{h}|<\infty$, as for
stationary ARMA processes, then $\tau_{k,T}-\tau_{k}=O(T^{-1})$ for each fixed
$k$ (Lemma~S.3 in the Online Supplement), and this cost is $O(n^{-1})$.
\end{remark}

\begin{remark}
[Bounds on the spectral moments and the AR(1) example]\label{rem:toeplitz}%
Since the eigenvalues $\nu_{t,T}$ are nonnegative with mean $\tau_{1,T}=1$,
Jensen's inequality gives $\tau_{k,T}\geq1$ for $k\geq2$, with equality if and
only if $\Gamma=I_{T}$. Moreover,
\[
K_{T}=\frac{1}{T}\sum_{t=1}^{T}\nu_{t,T}^{2}(\nu_{t,T}-\tau_{2,T})^{2}\geq0.
\]
For a stationary AR(1) correlation function $\gamma_{h}=\rho^{|h|}$ with
$|\rho|<1$, the symbol is $f(\omega)=(1-\rho^{2})/(1-2\rho\cos\omega+\rho
^{2})$, and
\begin{equation}
\tau_{2}=\frac{1+\rho^{2}}{1-\rho^{2}},\quad\tau_{3}=\frac{1+4\rho^{2}%
+\rho^{4}}{(1-\rho^{2})^{2}},\quad\tau_{4}=\frac{1+9\rho^{2}+9\rho^{4}%
+\rho^{6}}{(1-\rho^{2})^{3}}. \label{eq:AR1moments}%
\end{equation}
Since $\sum_{h\geq1}h|\rho|^{h}<\infty$, the finite-$T$ moments $\tau_{k,T}$
used in $\mu_{n}$ and $\sigma_{n}^{2}$ differ from these limits by $O(T^{-1})$ (Lemma~S.3).
\end{remark}

\section{Main Results}

\label{sec:main}

\subsection{Null mean and variance}

\label{sec:moments}

We expand each $r_{ij}^{2}=S_{ij}^{2}/(S_{ii}S_{jj})$ around $S_{ii}=S_{jj}%
=1$, which requires the sample variances to be uniformly close to one.

\begin{lemma}
[Diagonal Concentration]\label{lem:diag} Under Assumptions~1--3 and $H_{0}$,
\[
\Vert\operatorname{diag}(S)-I_{n}\Vert_{op}=O_{a.s.}\left(  \sqrt{\frac{\log
n}{T}}\right)  .
\]

\end{lemma}

The bound follows from Gaussian quadratic-form concentration and a union
bound; see Section~S.A.3 of the Online Supplement. Although the diagonal
elements of $S$ converge uniformly to one, this alone is insufficient to
establish the limiting distribution of $Tr(R^{2})$. The randomness of the
diagonal normalization must be accounted for in the asymptotic mean and
variance expansions.

\begin{lemma}
[Null Mean]\label{lem:mean} Under Assumptions~1--3 and $H_{0}$,
\[
\mathbb{E}\{Tr(R^{2})\}=n+\frac{n-1}{\lambda}\tau_{2,T} +\frac{4(n-1)}%
{n\lambda^{2}}(\tau_{2,T}^{2}-\tau_{3,T}) +O(n^{-1}).
\]

\end{lemma}

\begin{remark}
[Interpretation of Lemma~\ref{lem:mean}]\label{rem:mean-terms} Under $H_{0}$,
all ordered off-diagonal pairs are identically distributed, so
\[
\mathbb{E}\{Tr(R^{2})\}=n+n(n-1)\mathbb{E}(r_{12}^{2}).
\]
The pairwise expansion in Appendix~\ref{app:proof-mean} is
\[
\mathbb{E}(r_{12}^{2})=\frac{\tau_{2,T}}{T} +\frac{4(\tau_{2,T}^{2}-\tau
_{3,T})}{T^{2}}+O(T^{-3}).
\]
The leading term $\tau_{2,T}/T=\mathbb{E}(S_{12}^{2})$ captures temporal
correlation. Since $\tau_{2,T}\geq1$, with equality only when $\Gamma=I_{T}$,
centering as if observations were temporally independent ($\tau_{2,T}=1$)
incurs the order-$n$ error $(n-1)(\tau_{2,T}-1)/\lambda$. The second term is
the correction from the diagonal normalization: the linear denominator terms
contribute $-4\tau_{3,T}/T^{2}$, while the quadratic terms contribute
$4\tau_{2,T}^{2}/T^{2}+O(T^{-3})$. By Cauchy--Schwarz and $\tau_{1,T}=1$,
$\tau_{2,T}^{2}\leq\tau_{3,T}$; hence this second-order correction is
nonpositive. The $O(T^{-3})$ pairwise remainder sums to $O(n^{-1})$ when
$T\asymp n$.

For $\Gamma=I_{T}$, all $\tau_{k,T}=1$ and $r_{12}^{2}\sim Beta(1/2,(T-1)/2)$.
Thus $\mathbb{E}(r_{12}^{2})=1/T$ exactly, recovering $\mathbb{E}%
\{Tr(R^{2})\}=n+n(n-1)/T$.
\end{remark}

The variance has two components: the variance of individual squared
correlations and the covariance between squared correlations that share a
unit. Disjoint pairs are independent under the null. The shared-unit
covariance vanishes in the temporally independent benchmark; under temporal
dependence, its aggregate contribution can remain of order $O(1)$ as $n$ and
$T$ grow proportionally.

\begin{lemma}
[Null Variance]\label{lem:var} Under Assumptions~1--3 and $H_{0}$,
\[
Var\{Tr(R^{2})\}=\frac{4\tau_{2,T}^{2}}{\lambda^{2}} +\frac{8K_{T}}%
{\lambda^{3}}+O(n^{-1}).
\]

\end{lemma}

\begin{remark}
[Interpretation of Lemma~\ref{lem:var}]\label{rem:var-terms} Let
$Z_{i}^{\prime}$ denote the $i$th row of $Z$, so that $Z_{1},\ldots,Z_{n}$ are
i.i.d.\ $N(0,I_{T})$. Under $H_{0}$, disjoint pairs are independent, while
pairs sharing a unit are conditionally independent given their common vector
but not unconditionally independent. Consequently,
\[
Var\{Tr(R^{2})\} =2n(n-1)a_{T}+4n(n-1)(n-2)b_{T},
\]
where $a_{T}=Var(r_{12}^{2})$ and $b_{T}=Cov(r_{12}^{2},r_{13}^{2})$. The
expansions established in Appendix~\ref{app:proof-variance} are
\[
a_{T}=\frac{2\tau_{2,T}^{2}}{T^{2}}+O(T^{-3}),\qquad b_{T}=\frac{2K_{T}}%
{T^{3}}+O(T^{-4}).
\]
Thus the $O(n^{2})$ individual pairs generate $4\tau_{2,T}^{2}/\lambda^{2}$,
while the $O(n^{3})$ shared-unit configurations generate $8K_{T}/\lambda^{3}$.

By conditional independence, $b_{T}=Var\{\mathbb{E}(r_{12}^{2}\mid Z_{1})\}$.
With $q(z)=z^{\prime}\Gamma z$ and $B_{T}=\mathbb{E}\{ZZ^{\prime}/q(Z)\}$ for
$Z\sim N(0,I_{T})$,
\[
\mathbb{E}(r_{12}^{2}\mid Z_{1}=z)=\frac{z^{\prime}\Gamma B_{T}\Gamma z}%
{q(z)}.
\]
Its centered fluctuation is governed at leading order by
\[
W_{T}(Z_{1})=\frac{Z_{1}^{\prime}\Gamma^{2}Z_{1}}{T Z_{1}^{\prime}\Gamma
Z_{1}}, \qquad D_{1}=Z_{1}^{\prime}(\Gamma^{2}-\tau_{2,T}\Gamma)Z_{1},
\]
for which $Var(D_{1})=2TK_{T}$. The leading term of $W_{T}-\tau_{2,T}/T$ is
$D_{1}/T^{2}$, and $Var(W_{T})=Var(D_{1})/T^{4}+O(T^{-4})=2K_{T}%
/T^{3}+O(T^{-4})$. The exact conditional expectation differs from $W_{T}$ by a
constant plus a centered term of order $O_{P}(T^{-5/2})$, so it has the same
variance to this order (Appendix~\ref{app:proof-variance}, Steps~3--4). In
terms of the eigenvalues of $\Gamma$,
\[
K_{T}=\tau_{4,T}-2\tau_{2,T}\tau_{3,T}+\tau_{2,T}^{3} =\frac1T\sum_{t=1}%
^{T}\nu_{t,T}^{2}(\nu_{t,T}-\tau_{2,T})^{2}\geq0,
\]
a dispersion measure of the temporal eigenvalues weighted by their squares.
For $\Gamma=I_{T}$, $K_{T}=0$, the shared-unit covariance vanishes exactly,
and
\[
Var\{Tr(R^{2})\} =\frac{4n(n-1)(T-1)}{T^{2}(T+2)},
\]
which is the variance in the denominator of $LM_{S}$. Since $\tau_{k,T}\to
\tau_{k}$ and $K_{T}\to K$, replacing $\tau_{2,T}$ and $K_{T}$ by their limits
changes the leading variance by $o(1)$.

Carrying the expansions of Appendix~\ref{app:proof-variance} one order further
gives
\begin{align*}
a_{T}  &  =\frac{2\tau_{2,T}^{2}}{T^{2}}+\frac{28\tau_{2,T}^{3}-40\tau
_{2,T}\tau_{3,T}+6\tau_{4,T}}{T^{3}}+O(T^{-4}),\\
b_{T}  &  =\frac{2K_{T}}{T^{3}}+\frac{L_{T}}{T^{4}}+O(T^{-5}),
\end{align*}
where
\[
L_{T}=40\tau_{2,T}^{4}-104\tau_{2,T}^{2}\tau_{3,T}+60\tau_{2,T}\tau
_{4,T}+28\tau_{3,T}^{2}-24\tau_{5,T}.
\]
Consequently, with $T=n\lambda$,
\begin{align*}
Var\{Tr(R^{2})\}  &  =\frac{4\tau_{2,T}^{2}}{\lambda^{2}}+\frac{8K_{T}%
}{\lambda^{3}}+\frac{c_{T}}{n}+O(n^{-2}),\\
c_{T}  &  =-\frac{4\tau_{2,T}^{2}}{\lambda^{2}}+\frac{4(8\tau_{2,T}^{3}%
-8\tau_{2,T}\tau_{3,T}-3\tau_{4,T})}{\lambda^{3}}+\frac{4L_{T}}{\lambda^{4}}.
\end{align*}
For $\Gamma=I_{T}$, $L_{T}=0$ and $c_{T}=-4/\lambda^{2}-12/\lambda^{3}$, which
is the $n^{-1}$ term of the exact variance above. In the AR(1) designs of
Section~\ref{sec:mc} with $\rho>0$, $c_{T}$ is negative and large: for
$\rho=0.4$ and $n=T=100$, $c_{T}\approx-128$, and the leading variance exceeds
the exact variance by about $10\%$; for $\rho=0.7$, by about $26\%$. The
$O(1)$ terms are unaffected, so Theorem~\ref{thm:clt} is unchanged; the
$n^{-1}$ term accounts for the finite-sample conservativeness of the oracle
test reported in Section~\ref{sec:mc}.
\end{remark}

Relative to the temporally independent benchmark, temporal dependence changes
the leading mean and variance, adds an $O(1)$ mean correction from the
diagonal normalization, and introduces the $O(1)$ shared-unit variance
correction. The latter two corrections are of the same order as the null
fluctuations and cannot generally be omitted.

\subsection{Central limit theorem}

\label{sec:clt}

Let $\mu_{n}$ and $\sigma_{n}^{2}$ denote the centering and variance
expressions evaluated at the finite-$T$ spectral moments,
\begin{equation}
\mu_{n}=n+\frac{n-1}{\lambda}\tau_{2,T}+\frac{4(n-1)}{n\lambda^{2}}(\tau
_{2,T}^{2}-\tau_{3,T}),\qquad\sigma_{n}^{2}=\frac{4\tau_{2,T}^{2}}{\lambda
^{2}}+\frac{8K_{T}}{\lambda^{3}}. \label{eq:mu-sigma}%
\end{equation}
By Lemmas~\ref{lem:mean} and~\ref{lem:var}, $\mathbb{E}\{Tr(R^{2})\}=\mu
_{n}+O(n^{-1})$ and $Var\{Tr(R^{2})\}=\sigma_{n}^{2}+O(n^{-1})$. Replacing
$\tau_{k,T}$ and $K_{T}$ by their limits $\tau_{k}$ and $K$ in the $O(1)$
terms changes $\mu_{n}$ and $\sigma_{n}^{2}$ by $o(1)$
(Remarks~\ref{rem:szego} and~\ref{rem:var-terms}). Define the LM statistic
\[
LM=\frac{Tr(R^{2})-\mu_{n}}{\sigma_{n}}.
\]
The test rejects $H_{0}$ at level $\alpha$ when $|LM|>z_{\alpha/2}$, where
$z_{\alpha/2}$ is the $(1-\alpha/2)$ standard normal quantile.

\begin{thm}
[Null CLT]\label{thm:clt}Under Assumptions 1--3 and $H_{0}$,
$LM\overset{d}{\longrightarrow}N(0,1)$.
\end{thm}

Since $\tau_{2,T}\geq1$ and $K_{T}\geq0$ (Remark~\ref{rem:toeplitz}),
$\sigma_{n}^{2}\geq4/\lambda^{2}$ is bounded away from zero.
Theorem~\ref{thm:clt} continues to hold with the $O(1)$ terms of $\mu_{n}$ and
$\sigma_{n}^{2}$ evaluated at the limits $\tau_{k}$ and $K$; only the leading
mean term must retain $\tau_{2,T}$. When $\Gamma=I_{T}$, $\mu_{n}$ is the
centering of $LM_{S}$ and $\sigma_{n}^{2}$ equals its variance up to
$O(n^{-1})$, so Theorem~\ref{thm:clt} recovers the known-mean version of
Schott (2005). The test is an oracle test because it uses the population
spectral moments of $\Gamma$; feasible versions are investigated numerically
in Section~\ref{sec:mc}.

\begin{remark}
[Relation to Schott (2005)]\label{rem:schott} Schott (2005) proves a
high-dimensional CLT for sums of squared sample correlations under Gaussian
cross-sectional independence. Our result recovers the corresponding
\emph{known-mean} benchmark when $\Gamma=I_{T}$. Schott uses the usual sample
correlations after subtracting sample means, whereas our $R$ is formed without
temporal demeaning.

If $\Gamma=I_{T}$, the vectors $u_{i}=Z_{i}/\Vert Z_{i}\Vert$ are independent
and uniform on the unit sphere in $\mathbb{R}^{T}$. Hence $r_{ij}%
=u_{i}^{\prime}u_{j}$ and
\[
r_{ij}^{2}\sim Beta\left(  1/2,(T-1)/2\right)  ,\qquad\mathbb{E}(r_{ij}%
^{2}\mid u_{i})=1/T.
\]
The first-order Hoeffding projection therefore vanishes, and $\sum
_{i<k}(r_{ik}^{2}-1/T)$ is a martingale difference with respect to the
filtration generated by the first $k$ rows. In this benchmark $\tau_{2,T}=1$,
$K_{T}=0$, and the leading variance is $4/\lambda^{2}$.

For nonscalar $\Gamma$, normalized Gaussian rows are generally not uniform on
the sphere, and $\mathbb{E}(r_{ij}^{2}\mid Z_{i})$ need not be constant.
Writing $\theta_{T}=\mathbb{E}(r_{ij}^{2})$ and $\xi_{i}=\mathbb{E}(r_{ij}%
^{2}\mid Z_{i})-\theta_{T}$, the increments $\sum_{i<k}(r_{ik}^{2}-\theta
_{T})$ satisfy
\[
\mathbb{E}\left\{  \sum_{i<k}(r_{ik}^{2}-\theta_{T})|Z_{1},\ldots
,Z_{k-1}\right\}  =\sum_{i<k}\xi_{i},
\]
so they are not generally martingale differences. The Hoeffding--Doob
increments in Appendix~\ref{app:thm1} instead account for both projections:
\[
Y_{k}=2(n-1)\xi_{k}+2\sum_{i<k}h_{ik}^{\circ},
\]
where $h_{ik}^{\circ}=r_{ik}^{2}-\theta_{T}-\xi_{i}-\xi_{k}$ is the degenerate
component. The nondegenerate first-order projection yields the additional
$8K_{T}/\lambda^{3}$ variance component. Thus, replacing the
dependence-adjusted centering and variance by their temporally independent
benchmarks can lead to substantial size distortions.
\end{remark}

\section{Power of the LM Test}

\label{sec:power}

This section examines the LM statistic under cross-sectional dependence. We
derive its mean and variance under alternatives, establish consistency when
the aggregate signal $\sum_{i\neq j}(\rho_{ij}^{\ast})^{2}$ diverges, and
characterize local asymptotic power and the detection boundary. The proofs of
Lemmas~5 and~6, Theorem~\ref{thm:power}, and Corollary~\ref{cor:boundary} are
in Sections~S.B--S.E of the Online Supplement.

%-----------------------------------------------------------

\subsection{Setup Under the Alternative}

%-----------------------------------------------------------

\noindent Under the alternative, the data model remains $X=U^{1/2}%
Z\Gamma^{1/2}$ with $Z$ having i.i.d.\ $N(0,1)$ entries, but now $U\neq I_{n}%
$. We parameterize the cross-sectional dependence through the population
correlation matrix
\[
R_{\ast}=\left(  \operatorname{diag}(U)\right)  ^{-1/2}U\left(
\operatorname{diag}(U)\right)  ^{-1/2}%
\]
with $(i,j)$ entry $\rho_{ij}^{\ast}=U_{ij}/\sqrt{U_{ii}U_{jj}}$. The
hypotheses are $H_{0}:R_{\ast}=I_{n}$ and $H_{A}:R_{\ast}\neq I_{n}$. We write
$\mathbb{E}_{A}$, $Var_{A}$, and $P_{A}$ for moments and probabilities under
the alternative, and $\mathbb{E}_{0}$, $Var_{0}$ for those under $H_{0}$ with
the same $\Gamma$; thus $\mathbb{E}_{0}$ and $Var_{0}$ are the moments in
Lemmas~\ref{lem:mean} and~\ref{lem:var}. We define the cross-sectional
spectral moments of $R_{\ast}$ by
\[
\phi_{k}=\frac{1}{n}Tr(R_{\ast}^{k}),\qquad k=1,2,3,4.
\]
By construction, $\phi_{1}=1$. Under $H_{0}$, $\phi_{k}=1$ for all $k$. Under
$H_{A}$, $\phi_{2}>1$ and the aggregate signal is
\[
Tr(R_{\ast}^{2})-n=n(\phi_{2}-1)=\sum_{i\neq j}(\rho_{ij}^{\ast})^{2}>0.
\]
We call an alternative sequence local when $n(\phi_{2}-1)=O(1)$. By
Assumptions~2--3 and Remark~\ref{rem:toeplitz}, $\sigma_{n}$ is bounded above
and away from zero, so under local alternatives the aggregate signal is at
most of the order of the null standard deviation. We add one assumption
governing the alternative.

\begin{assumption}
\label{ass:U}(Cross-Sectional Covariance) The matrix $U$ is symmetric positive
definite with $c_{1}\leq\min_{i}U_{ii}\leq\Vert U\Vert_{op}\leq c_{2}$ for
constants $0<c_{1}\leq c_{2}<\infty$ independent of $n$.
\end{assumption}

\noindent Assumption~\ref{ass:U} is condition (1.1) of Heiny (2022): the lower
bound keeps each cross-sectional variance away from zero, and the
operator-norm bound controls the largest eigenvalue of $S$. Under this
condition Heiny (2022, Thm.~1.2) shows that $\operatorname{diag}(S)$
concentrates around the population diagonal for temporally independent data.
Hounyo, Kao, and Kim (2026a) extend this result to temporal dependence, under
finite fourth moments (Thm.~3) and, with the explicit rate $O(\sqrt{\log
n/T})$, under Gaussianity (Thm.~5); Lemma~\ref{lem:diag_alt} below is the
special case used here. Assumption~\ref{ass:U} implies (i) $\Vert R_{\ast
}\Vert_{op}\leq c_{2}/c_{1}<\infty$; (ii) $\phi_{k}=O(1)$ for each $k$; and
(iii) $Tr(R_{\ast}^{2})-n=O(n)$.

\begin{remark}
In testing $H_{0}:R_{\ast}=I_{n}$, the temporal covariance matrix $\Gamma$ is
a nuisance parameter: $H_{0}$ leaves it unrestricted, yet its moments
$\tau_{k,T}$ determine the null centering $\mu_{n}$ and variance $\sigma
_{n}^{2}$. The departure from $H_{0}$ is measured by $\phi_{2}$. Under local
alternatives, Lemmas~\ref{lem:mean_alt} and~\ref{lem:var_alt} give
$\mathbb{E}_{A}[Tr(R^{2})]=\mu_{n}+n(\phi_{2}-1)+O(n^{-1})$ and $Var_{A}%
(Tr(R^{2}))=\sigma_{n}^{2}+O(n^{-1})$: the alternative shifts the mean
additively and leaves the leading variance unchanged, so power is governed by
$\Delta_{n}=n(\phi_{2}-1)/\sigma_{n}$.
\end{remark}

%-----------------------------------------------------------

\subsection{Moments Under the Alternative}

%-----------------------------------------------------------

\noindent The first step extends the null diagonal-concentration result to
general $U$. Under $H_{A}$, $\mathbb{E}[S_{ii}]=U_{ii}$, and the diagonal
concentrates around $\operatorname{diag}(U)$.

\begin{lemma}
[Diagonal Concentration Under the Alternative]\label{lem:diag_alt} Under
Assumptions 1--4,
\[
\left\Vert \operatorname{diag}(S)-\operatorname{diag}(U)\right\Vert
_{op}=O_{a.s.}\left(  \sqrt{\frac{\log n}{T}}\right)  .
\]

\end{lemma}

\noindent The proof is in Section~S.A.3 of the Online Supplement.

Under the separable Gaussian model, for each $i$,
\[
\frac{S_{ii}}{U_{ii}}\overset{d}{=}\frac{1}{T}z_{i}^{\prime}\Gamma z_{i},
\qquad z_{i}\sim N(0,I_{T}),
\]
where $\overset{d}{=}$ denotes equality in distribution. Thus the standardized
diagonal elements have the same marginal distribution as under the null. The
uniform bounds on $U_{ii}$ in Assumption~4 yield the same concentration rate
as in Lemma~1.

\begin{lemma}
[Alternative Mean]\label{lem:mean_alt}Under Assumptions 1--4, the expectation
of $Tr(R^{2})$ under $H_{A}$ satisfies
\begin{align}
\mathbb{E}_{A}\left[  Tr(R^{2})\right]   &  =Tr(R_{\ast}^{2})+\frac
{n-1}{\lambda}\tau_{2,T}+\frac{4(n-1)}{\lambda^{2}n}\left(  \tau_{2,T}%
^{2}-\tau_{3,T}\right) \nonumber\\
&  -\frac{3\tau_{2,T}(\phi_{2}-1)}{\lambda}+\frac{2\tau_{2,T}}{T}\sum_{i\neq
j}(\rho_{ij}^{\ast})^{4}+O(n^{-1}). \label{eq:mean_alt}%
\end{align}

\end{lemma}

\begin{remark}
In terms of $\mu_{n}$ in \eqref{eq:mu-sigma}, Lemma~\ref{lem:mean_alt} reads
\[
\mathbb{E}_{A}[Tr(R^{2})]=\underbrace{\mu_{n}}_{\text{null mean}%
}+\underbrace{n(\phi_{2}-1)}_{\text{leading signal}}\underbrace{-\frac
{3\tau_{2,T}(\phi_{2}-1)}{\lambda}+\frac{2\tau_{2,T}}{T}\sum_{i\neq j}%
(\rho_{ij}^{\ast})^{4}}_{\text{interaction terms, }O(\phi_{2}-1)}+O(n^{-1}).
\]
Since $\sum_{i\neq j}(\rho_{ij}^{\ast})^{4}\leq n(\phi_{2}-1)$, the
interaction terms vanish under $H_{0}$ and are $O(n^{-1})$ under local
alternatives with $n(\phi_{2}-1)=O(1)$. Pair by pair, $\mathbb{E}_{A}%
(r_{ij}^{2})=(\rho_{ij}^{\ast})^{2}+\tau_{2,T}(1-(\rho_{ij}^{\ast}%
)^{2})(1-2(\rho_{ij}^{\ast})^{2})/T+O(T^{-2})$, which for $\Gamma=I_{T}$ is
the standard expansion for a bivariate Gaussian sample with known means. Since
$(1-\rho^{2})(1-2\rho^{2})=1-3\rho^{2}+2\rho^{4}$, summing the $T^{-1}$ term
over the $n(n-1)$ ordered pairs gives $(n-1)\tau_{2,T}/\lambda$ and the two
interaction terms. The $O(T^{-2})$ term equals $4(\tau_{2,T}^{2}-\tau
_{3,T})/T^{2}+O((\rho_{ij}^{\ast})^{2}T^{-2})+O(T^{-3})$, and summation over
pairs gives the remainder $O(n^{-1})$.
\end{remark}

\begin{remark}
\label{remark:power} By Lemma~\ref{lem:mean_alt}, the standardized mean shift
of $LM$ is
\[
\frac{\mathbb{E}_{A}[Tr(R^{2})]-\mu_{n}}{\sigma_{n}}=\Delta_{n}+O(\phi
_{2}-1)+O(n^{-1}),\qquad\Delta_{n}=\frac{n(\phi_{2}-1)}{\sigma_{n}}%
=\frac{\lambda\,n(\phi_{2}-1)}{2\sqrt{\tau_{2,T}^{2}+2K_{T}/\lambda}},
\]
and Theorem~\ref{thm:power} shows that $\Delta_{n}$ is the noncentrality
parameter of $LM$ under local alternatives. Temporal dependence reduces power
by inflating $\sigma_{n}$ through both $\tau_{2,T}$ and $K_{T}$. Using only
the first variance term, i.e.\ $2\tau_{2,T}/\lambda$ as the standard
deviation, would overstate $\Delta_{n}$ by the factor $\sqrt{1+2K_{T}%
/(\lambda\tau_{2,T}^{2})}$. For an AR(1) process with persistence $\rho$,
$\tau_{2,T}\rightarrow\tau_{2}=(1+\rho^{2})/(1-\rho^{2})$ (see
\eqref{eq:AR1moments}) and $K_{T}\rightarrow K=2\rho^{2}(1+\rho^{2}%
)/(1-\rho^{2})^{3}$, so that
\[
\sigma_{n}^{2}\rightarrow\frac{4(1+\rho^{2})^{2}}{(1-\rho^{2})^{2}\lambda^{2}%
}+\frac{16\rho^{2}(1+\rho^{2})}{(1-\rho^{2})^{3}\lambda^{3}}.
\]
At $\lambda=1$, the limiting factor $\sqrt{1+2K/(\lambda\tau_{2}^{2})}$ is
about $1.29$ for $\rho=0.4$ and $1.89$ for $\rho=0.7$ (about $1.72$ and $2.96$
at $\lambda=1/3$). Within the local-alternative regime, the leading power
depends on the alternative through the aggregate signal $n(\phi_{2}%
-1)=\sum_{i\neq j}(\rho_{ij}^{\ast})^{2}$. Many weak correlations are
therefore detected with power tending to one once the aggregate signal
diverges, whereas a fixed number of pairs with fixed correlations gives
$n(\phi_{2}-1)=O(1)$ and power bounded away from one.
\end{remark}

\begin{lemma}
[Alternative Variance]\label{lem:var_alt}Under Assumptions 1--4,
$Var_{A}(Tr(R^{2}))=O(1)$. If in addition $n(\phi_{2}-1)=O(1)$, then
\begin{equation}
Var_{A}(Tr(R^{2}))=\frac{4\tau_{2,T}^{2}}{\lambda^{2}}+\frac{8K_{T}}%
{\lambda^{3}}+O(n^{-1})=\sigma_{n}^{2}+O(n^{-1}), \label{eq:var_alt}%
\end{equation}
where $\sigma_{n}^{2}$ is the null variance in \eqref{eq:mu-sigma}.
\end{lemma}

\begin{remark}
The first statement holds uniformly over alternatives satisfying the stated
assumptions and is sufficient, together with a diverging mean shift, for
consistency; it follows from the Gaussian Poincar\'{e} inequality (Lemma~S.7
in the Online Supplement). The sharper local variance statement follows from
the pairwise covariance expansion in Lemma~S.9 in the Online Supplement,
combined with Lemma~\ref{lem:var}.
\end{remark}

\begin{remark}
[Persistence of the null variance]\label{rem:var-same} Under local
alternatives, the aggregate signal $n(\phi_{2}-1)=\sum_{i\neq j}(\rho
_{ij}^{\ast})^{2}$ is $O(1)$. It shifts the mean of $Tr(R^{2})$ by an $O(1)$
amount, whereas Lemma~S.9 shows that its variance changes by only $O(n^{-1}%
)$:
\[
Var_{A}\{Tr(R^{2})\}-Var_{0}\{Tr(R^{2})\}=O(n^{-1}).
\]
Thus both components of the null variance, $4\tau_{2,T}^{2}/\lambda^{2}$ and
$8K_{T}/\lambda^{3}$, remain the leading components under local alternatives.

The distinction reflects how cross-sectional correlations enter the two
moments. The leading mean shift sums $(\rho_{ij}^{\ast})^{2}$ directly. By
contrast, variance corrections involve covariances of squared sample
correlations, including pairs that share a unit and pairs that are disjoint.
Sign invariance restricts the correlation monomials in these covariances: each
index must have even total degree, although this condition alone does not
determine their asymptotic order. The uniform coefficient estimates and graph
summation in Lemma~S.9 are what establish the aggregate $O(n^{-1})$ bound.

A single pair shows the mechanism. Each $(\rho_{ij}^{\ast})^{2}$ enters
$\mathbb{E}_{A}(r_{ij}^{2})$ with coefficient one, but the variance of
$r_{ij}^{2}$ is already of order $T^{-1}$. A delta-method expansion gives
$Var_{A}(r_{ij})=(1-(\rho_{ij}^{\ast})^{2})^{2}\tau_{2,T}/T+O(T^{-2})$, so the
leading term of $Var_{A}(r_{ij}^{2})-Var_{0}(r_{ij}^{2})$ is $4(\rho
_{ij}^{\ast})^{2}(1-(\rho_{ij}^{\ast})^{2})^{2}\tau_{2,T}/T$. Summed over
pairs, these terms contribute at most of order $\tau_{2,T}n(\phi_{2}-1)/T$ to
$Var_{A}\{Tr(R^{2})\}$, which is $O(n^{-1})$ under local alternatives, whereas
the same signal shifts the mean by $n(\phi_{2}-1)$. Equivalently, the coupling
in Lemma~S.8 writes the alternative statistic as the null statistic plus the
constant $n(\phi_{2}-1)$ and an $L^{2}$-negligible remainder.

For fixed alternatives with $\phi_{2}-1$ bounded away from zero, $\sum_{i\neq
j}(\rho_{ij}^{\ast})^{2}$ is of order $n$, so the local variance comparison
need not apply. The local-power approximation below is therefore stated under
$n(\phi_{2}-1)=O(1)$.
\end{remark}

%-----------------------------------------------------------

\subsection{Consistency and Local Power}

\begin{proposition}
[Consistency Under Fixed Alternatives]\label{prop:fixed_consistency} Under
Assumptions 1--4, if $\phi_{2}-1$ is bounded away from zero, then
$P_{A}(|LM|>z_{\alpha/2})\to1$.
\end{proposition}

\noindent Lemma~\ref{lem:mean_alt} gives a mean shift of order $n$ under the
stated fixed-signal condition, while Lemma~\ref{lem:var_alt} bounds the
alternative variance by $O(1)$. Chebyshev's inequality therefore establishes
the result. The following theorem treats local alternatives, for which the
mean shift remains bounded.

\begin{thm}
[CLT Under Local Alternatives]\label{thm:power}Under Assumptions 1--4, suppose
that $n(\phi_{2}-1)=O(1)$. Then
\[
LM-\Delta_{n}\overset{d}{\longrightarrow}N(0,1),\qquad\Delta_{n}=\frac
{n(\phi_{2}-1)}{\sigma_{n}}.
\]

\end{thm}

\begin{remark}
[Role of the Coupling Argument]The null CLT alone does not imply the
Gaussian-shift approximation in Theorem~\ref{thm:power}. Under local
alternatives, the sample correlations change through both their numerators and
their diagonal normalization; moreover, covariances between squared
correlations sharing a unit can accumulate as $n$ and $T$ grow together. Even
equality of the leading null and alternative variances would not, by itself,
establish it.

The required distributional comparison follows instead from Lemma~S.8: using
common Gaussian innovations, it shows that the difference between the
alternative statistic and its null counterpart, after subtracting $n(\phi
_{2}-1)$, is $o_{P_{A}}(1)$. Theorem~\ref{thm:power} then follows from the
null CLT and Slutsky's theorem (Section~S.D of the Online Supplement).
Separately, Lemma~S.9 establishes the sharper $O(n^{-1})$ comparison of null
and local-alternative variances. The coupling and variance comparison
therefore serve distinct purposes.
\end{remark}

The next result, proved in Section~S.E of the Online Supplement, gives the
signal scale $n(\phi_{2}-1)$ at which the LM test has nontrivial asymptotic power.

\begin{corollary}
[Detection Boundary]\label{cor:boundary}Under Assumptions 1--4:

\noindent(a) the LM test has power tending to one if and only if $n(\phi
_{2}-1)\rightarrow\infty$, equivalently,
\[
\frac{Tr(R_{\ast}^{2})-n}{\sigma_{n}}\to\infty;
\]
\noindent(b) the test has asymptotic power equal to $\alpha$ if and only if
$n(\phi_{2}-1)\rightarrow0$;

\noindent(c) if $n(\phi_{2}-1)/\sigma_{n}\rightarrow\Delta\in(0,\infty)$, the
power converges to $\Phi(\Delta-z_{\alpha/2})+\Phi(-\Delta-z_{\alpha/2}%
)\in(\alpha,1)$.
\end{corollary}

\begin{remark}
The detection boundary for $n(\phi_{2}-1)$ is of order
\[
\sigma_{n}=\frac{2}{\lambda}\sqrt{\tau_{2,T}^{2}+2K_{T}/\lambda}.
\]
For fixed $\tau_{2,T}$ and $K_{T}$, it decreases as $\lambda$ increases. Since
$\tau_{2,T}\geq1$ and $K_{T}\geq0$, it is at least $2/\lambda=2n/T$, with
equality if and only if $\Gamma=I_{T}$ (Remark~\ref{rem:toeplitz}). Temporal
dependence therefore raises the boundary even for the appropriately
standardized $LM$ statistic. For $\Gamma=I_{T}$, $\sigma_{n}^{2}=4/\lambda
^{2}$ agrees with the independence-test case of Zou, Yang, Lin, Tian, and
Zheng (forthcoming, Cor.~1). In this case the boundary scale $n(\phi
_{2}-1)\asymp1$, equivalently $\Vert R_{\ast}-I_{n}\Vert_{F}\asymp1$ when
$\lambda$ stays bounded above and away from zero (in rate notation, $\Vert
R_{\ast}-I_{n}\Vert_{F}\asymp\sqrt{n/T}$), coincides with the Frobenius-norm
scale of the minimax lower bound of Zheng, Cheng, Guo, and Zhu (2019,
Thm.~2.3) for independent observations, so in that case the test detects dense
alternatives at the optimal rate.
\end{remark}

Section~\ref{sec:mc} reports the finite-sample power of the test under factor
and spatial alternatives (Tables~\ref{tab:power-independent}
and~\ref{tab:power-dependent}).

\section{Monte Carlo Simulation}

\label{sec:mc}

This section examines the finite-sample size and power of the proposed test
and compares it with Schott's (2005) test $LM_{S}$. All experiments use 5,000
Monte Carlo replications and the 5\% nominal level.

\subsection{Design and implementation}

Under the null hypothesis, the data are generated as
\[
x_{it} = \rho x_{it-1} + \sqrt{1-\rho^{2}} v_{it},
\]
where $x_{i0}$ and $v_{it}$ are i.i.d. draws from one of the following
distributions: (i) $N(0,1)$, (ii) $\chi^{2}(1)-1$, and (iii) $t(5)$. Under the
null hypothesis, the rows $x_{1},\ldots,x_{n}$ are cross-sectionally
independent AR(1) processes with autoregressive parameter $\rho\in
\{0,0.4,0.7\}$, and the temporal covariance matrix is the Toeplitz matrix
$\Gamma_{ts}=\rho^{|t-s|}$, whose limiting spectral moments are given in
\eqref{eq:AR1moments}. Since $x_{i0}$ and $v_{it}$ have a common variance,
each row has covariance matrix proportional to $\Gamma$; because $R$ is
invariant to rescaling the rows, the innovations need not be standardized. We
consider $n,T\in\{100,200,300\}$.

We compute three versions of the proposed test. The oracle test $LM_{\inf}$
uses $\mu_{n}$ and $\sigma_{n}^{2}$ in \eqref{eq:mu-sigma}, evaluated at the
finite-$T$ spectral moments $\tau_{k,T}$, $k=2,3,4$, of the data-generating
$\Gamma$. The feasible tests $LM_{UT}$ and $LM_{T}$ replace $\tau_{k,T}$ by
plug-in estimates based on the pooled untapered and tapered autocovariance
estimators of Hounyo, Kao, and Kim (2026b), normalized as in
Section~\ref{sec:model}.

The tapered estimator is defined as
\[
\hat{\Gamma}_{\mathcal{K}}=\frac{1}{n}\sum_{i=1}^{n}\hat{\Gamma}%
_{i,\mathcal{K}}%
\]
where $\hat{\Gamma}_{i,\mathcal{K}}$ is the $T\times T$ Toeplitz matrix whose
$h$-th diagonal elements, for $h=0,\ldots,T-1$, are given by
\[
\mathcal{K}\left(  \frac{h}{m_{T}}\right)  \frac{1}{T}\sum_{t=1}^{T-h}%
x_{it}x_{it+h}.
\]
Here $\mathcal{K}(\cdot)$ is a kernel function and $m_{T}$ is the bandwidth
parameter. In our simulations, we employ the rectangular kernel,
$\mathcal{K}(z)=1\{|z|\leq1\}$, which exhibited substantially better size
control than the Bartlett and Parzen kernels in additional simulations not
reported here. The untapered estimator corresponds to the full-bandwidth
rectangular-kernel case with $m_{T} = T.$

For bandwidth selection in the tapered estimator, we adopt the following
data-driven rule from Hounyo et al.\ (2026b):
\begin{align*}
\hat{m}_{T} = \max\{1,\min\{T-1,\lfloor\tilde{m}_{T} \rfloor\} \}
\end{align*}
where
\[
\tilde{m}_{T}=0.8\left\{  \frac{\hat{\alpha}^{2}nT}{\log(nT)} \right\}  ^{1/3}
\quad\hat{\alpha}=\frac{2\sum_{h=1}^{l} h \left|  \hat{\gamma}(h)\right|
}{\left|  \hat{\gamma}(0)\right|  + 2\sum_{h=1}^{l} \left|  \hat{\gamma
}(h)\right|  }, \quad\hat\gamma(h) = \frac{1}{nT} \sum_{i=1}^{n}\sum
_{t=1}^{T-h} x_{it}x_{i,t+h}.
\]
with $l = \min\{T-1,\lfloor(nT)^{1/5} \rfloor\}$. This rule was developed for
the Bartlett kernel; we use it for the rectangular kernel as a practical choice.

\subsection{Size}

Table~\ref{tab:size} reports the empirical sizes. We summarize four findings.

\medskip\noindent\textbf{Schott's test.} $LM_{S}$ is accurately sized when
$\rho=0$, but its empirical size is 1.000 for $\rho\in\{0.4,0.7\}$ in every
design and error distribution.

\medskip\noindent\textbf{Oracle test.} The $LM_{\inf}$ test provides accurate
empirical sizes with or without temporal dependence, as Theorem~1 predicts.
Its slight under-rejection under strong persistence reflects the $O(n^{-1})$
term in Remark~\ref{rem:var-terms}: both variance components fall short of
their leading approximations by a relative amount of order $T^{-1}$. For
$\rho=0.4$, the second-order terms of $a_{T}$ and $b_{T}$ in
Remark~\ref{rem:var-terms} are $-7.7\%$ for $a_{T}$ and $-9.6\%$ for $b_{T}$
at $T=100$, and $-3.9\%$ and $-4.8\%$ at $T=200$; for $\rho=0.7$ they are
about $-23\%$ at $T=100$. The leading formula therefore overstates
$Var(Tr(R^{2}))$ in small samples. Under normal innovations, the largest
deviations from the nominal level are at $\rho=0.7$ and $(n,T)=(200,100)$ and
$(100,100)$, where the empirical sizes are 0.030 and 0.031. When $\rho=0$,
$LM_{\inf}$ and $LM_{S}$ have nearly identical sizes.

\medskip\noindent\textbf{Non-Gaussian innovations.} Theorem~1 relies on
Gaussianity. Under the $t(5)$ distribution, the empirical sizes of $LM_{\inf}$
remain close to those under normality. Under the $\chi^{2}(1)-1$ distribution,
however, both $LM_{\inf}$ and $LM_{S}$ over-reject even without temporal
dependence (empirical sizes between 0.059 and 0.092 at $\rho=0$), consistent
with the fourth-cumulant term that enters the variance of $Tr(R^{2})$ for
non-Gaussian data (Gao et al., 2017). These distortions become larger when
temporal dependence is strong and $T$ is relatively small. For example, when
$\rho=0.7$ and $(n,T)=(300,100)$, the empirical size of $LM_{\inf}$ reaches
0.135 under the $\chi^{2}(1)-1$ distribution, compared with 0.034 under
normality. These distortions diminish as $T$ increases.

\medskip\noindent\textbf{Feasible tests.} When $\rho=0$, $LM_{UT}$ and
$LM_{T}$ have empirical sizes close to the nominal level. In the presence of
temporal dependence and under normal or $t(5)$ innovations, both tests are
markedly conservative: across the designs with $\rho\in\{0.4,0.7\}$ their
empirical sizes range from 0.003 to 0.035, except at $(n,T)=(300,100)$ with
$\rho=0.7$, where they lie between 0.048 and 0.103. This conservativeness does
not diminish when $n$ and $T$ grow together. For example, when $\rho=0.4$ and
the innovations are normally distributed, the empirical sizes of $LM_{UT}$ and
$LM_{T}$ are 0.008 and 0.014 at $(n,T)=(100,100)$ and 0.012 and 0.018 at
$(n,T)=(300,300)$. It diminishes, however, as $T$ grows with $n$ fixed. The
tapered test is generally less conservative than the untapered test. Under the
$\chi^{2}(1)-1$ distribution with strong temporal dependence and $T$ small
relative to $n$, the feasible tests instead over-reject substantially: when
$\rho=0.7$ and $(n,T)=(300,100)$, the empirical sizes of $LM_{UT}$ and
$LM_{T}$ are 0.540 and 0.583, respectively; these over-rejections diminish as
$T$ increases.

Table~\ref{tab:hybrid-size} separates the effects of estimating the centering
and the scaling terms, for $\rho=0.4$ and normal innovations. The hybrid tests
$LM_{T,\mu}$ and $LM_{UT,\mu}$, which use the true centering and estimated
scaling, have empirical sizes close to those of $LM_{\inf}$. The hybrid tests
$LM_{T,\sigma}$ and $LM_{UT,\sigma}$, which use the estimated centering and
true scaling, remain conservative, with rejection frequencies similar to those
of $LM_{T}$ and $LM_{UT}$. Estimation of the centering term, rather than the
scaling term, is therefore primarily responsible for the under-rejection.

A first-order expansion of the normalized spectral-moment estimator gives
\[
\widehat{\tau}_{2,T}-\tau_{2,T}\approx\frac{2}{nT}\sum_{i=1}^{n}D_{i}, \qquad
D_{i}=Z_{i}^{\prime}(\Gamma^{2}-\tau_{2,T}\Gamma)Z_{i}.
\]
Because the leading null centering contains $(n-1)\tau_{2,T}/\lambda$, its
plug-in estimation error is approximately
\[
\frac{n-1}{\lambda}(\widehat{\tau}_{2,T}-\tau_{2,T}) \approx\frac
{2(n-1)}{T^{2}}\sum_{i=1}^{n}D_{i}.
\]
This matches, to first order, the first-order Hoeffding projection of
$Tr(R^{2})$ that generates the additional oracle variance $8K_{T}/\lambda^{3}%
$. Thus, estimating $\tau_{2,T}$ from the same observations may absorb that
projection, leaving a plug-in-centered fluctuation with variance closer to the
degenerate component $4\tau_{2,T}^{2}/\lambda^{2}$ than to the full oracle
variance. Standardizing by the full oracle variance would then make the
feasible tests conservative. This argument is heuristic.

\subsection{Power}

We consider two alternative specifications of cross-sectional dependence: a
factor model and a spatial autoregressive (SAR(1)) model.

For the factor model, the data are generated as
\[
x_{it}=\alpha_{i} f_{t}+\varepsilon_{it}, \qquad\varepsilon_{it}
=\rho\varepsilon_{i,t-1} +\sqrt{1-\rho^{2}}v_{it},
\]
where $\varepsilon_{i0}$ and $v_{it}$ are i.i.d.\ $N(0,1)$, $\alpha
_{i}\overset{iid}{\sim} \kappa U(0,1)$, and $f_{t}\overset{iid}{\sim}N(0,1)$,
with all random variables mutually independent. Thus, $\kappa$ determines the
strength of cross-sectional dependence induced by the common factor.
Conditional on the factor loadings, the covariance matrix is
$Var(\operatorname{vec}(X)) = I_{T}\otimes\alpha\alpha^{\prime}+\Gamma\otimes
I_{n},$ where $\Gamma$ is the AR(1) temporal covariance matrix. This
covariance structure is generally not separable when $\rho\neq0$. Moreover,
the factor model generates a spiked cross-sectional covariance matrix, with
$\|U\|_{op}\asymp n$ for fixed $\kappa>0$, and therefore lies outside Assumption~4.

For the SAR(1) model, the data are generated as
\begin{align*}
x_{t}  &  = \psi W_{n} x_{t} + \varepsilon_{t}, \ \ x_{t}=(x_{1t}%
,\ldots,x_{nt})^{\prime}, \ \ \varepsilon_{t}=(\varepsilon_{1t},\ldots
,\varepsilon_{nt})^{\prime}\\
\varepsilon_{it}  &  = \rho\varepsilon_{it-1} + \sqrt{1-\rho^{2}} v_{it},
\quad\ \varepsilon_{i0},v_{it} \sim^{iid}N(0,1),\\
W_{n}  &  = \left[  w_{ij} \right]  _{ij=1}^{n} \text{ is an }(n \times
n)\text{ spatial contiguity matrix,}%
\end{align*}
where the spatial coefficient $\psi$ determines the degree of cross-sectional
dependence. In the SAR(1) DGP, the units are placed on a finite square integer
lattice containing $n$ sites. $W_{n}$ is a contiguity matrix in which
individuals i and j are considered to be neighbors if their $L_{1}$ distance
equals 1. The matrix is row-standardized and has zero diagonal entries. This
design satisfies Assumption~4 and the separable Gaussian representation.

\medskip\noindent\textbf{No temporal dependence.}
Table~\ref{tab:power-independent} reports power for $\rho=0$. Since all four
tests generally exhibit empirical sizes close to the nominal level under
normality when $\rho=0$, we compare their rejection probabilities without size
adjustment. Under the factor model, $LM_{\inf}$ exhibits power comparable to
$LM_{S}$, with rejection probabilities increasing toward one as the sample
sizes and/or the strength of cross-sectional dependence increase. The feasible
tapered test $LM_{T}$ performs almost as well as $LM_{\inf}$, whereas the
untapered test $LM_{UT}$ exhibits power losses, particularly when
cross-sectional dependence is relatively weak. Similar patterns emerge under
the SAR(1) model: $LM_{T}$ achieves power close to that of $LM_{\inf}$ and
$LM_{S}$, while $LM_{UT}$ tends to be less powerful.

\medskip\noindent\textbf{Temporal dependence.} Table~\ref{tab:power-dependent}
reports power for $\rho=0.4$. Because $LM_{S}$ is invalid under temporal
dependence, we report only the proposed tests. The feasible tests are
conservative in this design (Table~\ref{tab:size}), so their unadjusted power
is understated relative to tests with exact size. The power of $LM_{\inf}$ is
generally lower than that reported in Table~\ref{tab:power-independent},
consistent with Remark~\ref{remark:power}: temporal dependence inflates
$\sigma_{n}$ and lowers $\Delta_{n}$. The rejection probabilities approach one
as the strength of cross-sectional dependence increases and, in many
configurations, as the sample sizes increase. The tapered test $LM_{T}$ tends
to reject more often than the untapered test $LM_{UT}$ under both
alternatives. Moreover, when the cross-sectional dependence is weaker
($\kappa=0.2$ or $\psi=0.2$), $LM_{T}$ rejects more often than $LM_{\inf}$ in
every configuration, although it is conservative under the null.

\medskip\noindent\textbf{Accuracy of the local-alternative approximation.}
Theorem~2 rests on the coupling in Lemma~S.8: under a local alternative,
$Tr(R^{2})=Tr(R_{0}^{2})+n(\phi_{2}-1)+D_{n}$, where $R_{0}$ is the null
sample correlation matrix built from the same innovations and $\mathbb{E}%
(D_{n}^{2})\rightarrow0$. To gauge the size of $D_{n}$, we generate
$X=R_{\ast}^{1/2}Z\Gamma^{1/2}$ and $X_{0}=Z\Gamma^{1/2}$ from the same matrix
$Z$ of i.i.d.\ $N(0,1)$ entries, where $\Gamma$ is the AR(1) temporal
covariance matrix with $\rho=0.4$ and $R_{\ast}=I_{n}+E$ is equicorrelated,
with $E_{ij}=e_{n}=2/\sqrt{n(n-1)}$ for $i\neq j$, so that $n(\phi
_{2}-1)=\Vert E\Vert_{F}^{2}=4$ for every $n$. We take $n=T\in\{100,200,400\}$%
. The Monte Carlo mean of $D_{n}$ is at most $0.21$ in absolute value, and its
standard deviations are $2.19$, $1.63$, and $1.16$, respectively. The standard
deviation falls at the rate $n^{-1/2}$ ($\sqrt{n}$ times the standard
deviation is $21.9$, $23.1$, and $23.1$), as the bound $\mathbb{E}(D_{n}%
^{2})=O(T^{-1})$ in the proof of Lemma~S.8 predicts. Relative to the null
standard deviation $\sigma_{n}\approx3.5$ in this design, the approximation
error is not negligible at $n=T=100$ but shrinks steadily with the sample size.

\section{Conclusion and Discussion}

\label{sec:conclusion}

We develop a high-dimensional LM test of cross-sectional independence for
Gaussian panels with a common stationary, short-memory temporal covariance
matrix. Under proportional growth of the cross-sectional and time dimensions,
we derive explicit expansions for the null mean and variance of
$\operatorname{Tr}(R^{2})$ and establish its asymptotic normality. The results
identify two distinct consequences of serial correlation: random diagonal
normalization produces an $O(1)$ centering correction, and squared sample
correlations sharing a cross-sectional unit generate an additional $O(1)$
variance component. Both vanish under temporal independence.

The proofs combine uniform concentration of sample variances, Gaussian
quadratic-form expansions, a Hoeffding decomposition, and a martingale CLT.
Under alternatives, we establish consistency, a local Gaussian-shift limit,
and a detection boundary for the proposed statistic. The local-power analysis
separates the cross-sectional signal $n(\phi_{2}-1)$ from the temporal
spectral quantities determining the null variance.

The theoretical results are established under Gaussian innovations and a
common stationary, short-memory temporal covariance structure. Extending the
asymptotic theory to non-Gaussian innovations and heterogeneous temporal
covariance structures, and establishing the limiting distribution of feasible
statistics based on estimated temporal spectral moments, are natural
directions for future research. The extension to fixed-effects panel
regressions with serially correlated disturbances also requires accounting for
the effects of within transformation and slope estimation.

\section*{Acknowledgments}

During the preparation of this manuscript, the authors used Claude (Anthropic)
and ChatGPT (OpenAI) to assist with manuscript drafting and proofreading. No
AI tool was used to write the simulation code or to produce the numerical
results. The authors reviewed and edited all generated text and remain fully
responsible for the final version.

\section*{Competing interests}

The authors declare none.

\section*{Funding statement}

This research received no specific grant from any funding agency, commercial
or not-for-profit sectors.

\section*{References}

%References: uniform 2em hanging indentation; no list environment.
\begingroup
\setlength{\leftskip}{2em}\setlength{\parindent}{-2em}\setlength{\parskip}{0.4\baselineskip}\raggedright
Baltagi, B. H., Feng, Q., \& Kao, C. (2012). A Lagrange multiplier test for
cross-sectional dependence in a fixed effects panel data model.
\textit{Journal of Econometrics}, \textit{170}(1), 164--177.
\href{https://doi.org/10.1016/j.jeconom.2012.04.004}{\mbox{doi:10.1016/j.jeconom.2012.04.004}}%

Baltagi, B. H., Kao, C., \& Peng, B. (2016). Testing cross-sectional
correlation in large panel data models with serial correlation.
\textit{Econometrics}, \textit{4}(4), Article 44.
\href{https://doi.org/10.3390/econometrics4040044}{\mbox{doi:10.3390/econometrics4040044}}%

Breusch, T. S., \& Pagan, A. R. (1980). The Lagrange multiplier test and its
applications to model specification in econometrics. \textit{The Review of
Economic Studies}, \textit{47}(1), 239--253.
\href{https://doi.org/10.2307/2297111}{\mbox{doi:10.2307/2297111}}

Chen, X., \& Liu, W. (2018). Testing independence with high-dimensional
correlated samples. \textit{The Annals of Statistics}, \textit{46}(2),
866--894.
\href{https://doi.org/10.1214/17-AOS1571}{\mbox{doi:10.1214/17-AOS1571}}

Gao, J., Han, X., Pan, G., \& Yang, Y. (2017). High dimensional correlation
matrices: The central limit theorem and its applications. \textit{Journal of
the Royal Statistical Society: Series B (Statistical Methodology)},
\textit{79}(3), 677--693.
\href{https://doi.org/10.1111/rssb.12189}{\mbox{doi:10.1111/rssb.12189}}

Grenander, U., \& Szeg\H{o}, G. (1958). \textit{Toeplitz forms and their
applications}. University of California Press.

Hall, P., \& Heyde, C. C. (1980). \textit{Martingale limit theory and its
application}. Academic Press.

Heiny, J. (2022). Large sample correlation matrices: A comparison theorem and
its applications. \textit{Electronic Journal of Probability}, \textit{27},
1--20. \href{https://doi.org/10.1214/22-EJP817}{\mbox{doi:10.1214/22-EJP817}}

Hounyo, U., Kao, C., \& Kim, M. S. (2025). Serial dependence robust bootstrap
test for cross-sectional correlation. \textit{The Econometrics Journal},
\textit{28}(3), 442--461.
\href{https://doi.org/10.1093/ectj/utaf009}{\mbox{doi:10.1093/ectj/utaf009}}

Hounyo, U., Kao, C., \& Kim, M. S. (2026a). \textit{Limiting spectral
distribution for large sample correlation matrices under temporal dependence}
[Working paper]. University of Connecticut.

Hounyo, U., Kao, C., \& Kim, M. S. (2026b). \textit{Pooled tapered
autocovariance estimators in large panel model} [Working paper]. University of Connecticut.

Huang, Z., Li, Z., \& Yao, J. (2025). Unified and robust tests for cross
sectional independence in large panel data models. \textit{Electronic Journal
of Statistics}, \textit{19}(2), 4867--4913.
\href{https://doi.org/10.1214/25-EJS2456}{\mbox{doi:10.1214/25-EJS2456}}

Li, Y., Pan, G., Xie, J., \& Zhou, W. (2024). \textit{Necessary and sufficient
condition for CLT of linear spectral statistics of sample correlation
matrices} (arXiv:2409.12536). arXiv.
\href{https://doi.org/10.48550/arXiv.2409.12536}{\mbox{doi:10.48550/arXiv.2409.12536}}%

Schott, J. R. (2005). Testing for complete independence in high dimensions.
\textit{Biometrika}, \textit{92}(4), 951--956.
\href{https://doi.org/10.1093/biomet/92.4.951}{\mbox{doi:10.1093/biomet/92.4.951}}%

Yin, Y., Zheng, S., \& Zou, T. (2023). Central limit theorem of linear
spectral statistics of high-dimensional sample correlation matrices.
\textit{Bernoulli}, \textit{29}(2), 984--1006.
\href{https://doi.org/10.3150/22-BEJ1487}{\mbox{doi:10.3150/22-BEJ1487}}

Zheng, S., Cheng, G., Guo, J., \& Zhu, H. (2019). Test for high-dimensional
correlation matrices. \textit{The Annals of Statistics}, \textit{47}(5),
2887--2921.
\href{https://doi.org/10.1214/18-AOS1768}{\mbox{doi:10.1214/18-AOS1768}}

Zou, T., Yang, G., Lin, R., Tian, G.-L., \& Zheng, S. (forthcoming).
Structural testing of high-dimensional correlation matrices.
\textit{Statistica Sinica}.
\href{https://doi.org/10.5705/ss.202024.0078}{\mbox{doi:10.5705/ss.202024.0078}}%

\endgroup

%\pagebreak \setcounter{page}{1}
\newpage

\begin{singlespace}
\begin{center}
\refstepcounter{table}\label{tab:size}\textbf{Table \thetable.
Empirical size at the 5\% nominal level}
\end{center}
\begin{center}
\textit{Panel A: $N(0,1)$ and $\chi^{2}(1)-1$ innovations}
\medskip
\setlength{\tabcolsep}{4pt}\renewcommand{\arraystretch}{1.05}
\begin{tabular}[c]{c|cccc|cccc}\hline\hline
$(n,T)$ & \multicolumn{4}{c|}{$N(0,1)$} & \multicolumn{4}{c}{$\chi^{2}(1)-1$}\\
& $LM_{\inf}$ & $LM_{UT}$ & $LM_{T}$ & $LM_{S}$ & $LM_{\inf}$ & $LM_{UT}$ & $LM_{T}$ & $LM_{S}$\\\hline
& \multicolumn{8}{c}{\underline{$\rho=0.0$}}\\
(100,100) & 0.047 & 0.053 & 0.046 & 0.051 & 0.084 & 0.096 & 0.083 & 0.092\\
(200,100) & 0.045 & 0.050 & 0.046 & 0.049 & 0.083 & 0.092 & 0.083 & 0.088\\
(300,100) & 0.044 & 0.049 & 0.042 & 0.049 & 0.087 & 0.088 & 0.086 & 0.092\\
(100,200) & 0.041 & 0.051 & 0.041 & 0.044 & 0.070 & 0.081 & 0.069 & 0.074\\
(200,200) & 0.044 & 0.056 & 0.044 & 0.046 & 0.071 & 0.084 & 0.073 & 0.074\\
(300,200) & 0.049 & 0.063 & 0.050 & 0.052 & 0.067 & 0.076 & 0.066 & 0.069\\
(100,300) & 0.046 & 0.054 & 0.046 & 0.049 & 0.065 & 0.076 & 0.065 & 0.068\\
(200,300) & 0.045 & 0.057 & 0.045 & 0.047 & 0.068 & 0.080 & 0.067 & 0.070\\
(300,300) & 0.051 & 0.059 & 0.050 & 0.052 & 0.059 & 0.070 & 0.060 & 0.060\\\hline
& \multicolumn{8}{c}{\underline{$\rho=0.4$}}\\
(100,100) & 0.040 & 0.008 & 0.014 & 1.000 & 0.064 & 0.026 & 0.045 & 1.000\\
(200,100) & 0.038 & 0.006 & 0.010 & 1.000 & 0.066 & 0.046 & 0.067 & 1.000\\
(300,100) & 0.038 & 0.003 & 0.007 & 1.000 & 0.082 & 0.073 & 0.099 & 1.000\\
(100,200) & 0.040 & 0.016 & 0.021 & 1.000 & 0.052 & 0.030 & 0.039 & 1.000\\
(200,200) & 0.038 & 0.008 & 0.012 & 1.000 & 0.052 & 0.019 & 0.032 & 1.000\\
(300,200) & 0.043 & 0.006 & 0.011 & 1.000 & 0.051 & 0.015 & 0.029 & 1.000\\
(100,300) & 0.047 & 0.028 & 0.028 & 1.000 & 0.058 & 0.033 & 0.038 & 1.000\\
(200,300) & 0.039 & 0.013 & 0.015 & 1.000 & 0.056 & 0.022 & 0.031 & 1.000\\
(300,300) & 0.044 & 0.012 & 0.018 & 1.000 & 0.045 & 0.012 & 0.022 & 1.000\\\hline
& \multicolumn{8}{c}{\underline{$\rho=0.7$}}\\
(100,100) & 0.031 & 0.006 & 0.012 & 1.000 & 0.066 & 0.048 & 0.076 & 1.000\\
(200,100) & 0.030 & 0.017 & 0.023 & 1.000 & 0.090 & 0.230 & 0.277 & 1.000\\
(300,100) & 0.034 & 0.048 & 0.062 & 1.000 & 0.135 & 0.540 & 0.583 & 1.000\\
(100,200) & 0.036 & 0.004 & 0.009 & 1.000 & 0.041 & 0.009 & 0.023 & 1.000\\
(200,200) & 0.035 & 0.004 & 0.009 & 1.000 & 0.049 & 0.025 & 0.039 & 1.000\\
(300,200) & 0.040 & 0.006 & 0.011 & 1.000 & 0.056 & 0.044 & 0.068 & 1.000\\
(100,300) & 0.037 & 0.007 & 0.014 & 1.000 & 0.046 & 0.011 & 0.019 & 1.000\\
(200,300) & 0.038 & 0.005 & 0.010 & 1.000 & 0.042 & 0.009 & 0.019 & 1.000\\
(300,300) & 0.039 & 0.006 & 0.010 & 1.000 & 0.040 & 0.010 & 0.021 & 1.000\\\hline\hline
\end{tabular}
\end{center}
\newpage
\begin{center}
\textbf{Table \ref{tab:size} (continued).
Empirical size at the 5\% nominal level}
\end{center}
\begin{center}
\textit{Panel B: $t(5)$ innovations}\par
\medskip
\setlength{\tabcolsep}{4pt}\renewcommand{\arraystretch}{1.05}
\begin{tabular}[c]{c|cccc}\hline\hline
$(n,T)$ & \multicolumn{4}{c}{$t(5)$}\\
& $LM_{\inf}$ & $LM_{UT}$ & $LM_{T}$ & $LM_{S}$\\\hline
& \multicolumn{4}{c}{\underline{$\rho=0.0$}}\\
(100,100) & 0.055 & 0.061 & 0.055 & 0.059\\
(200,100) & 0.050 & 0.059 & 0.049 & 0.054\\
(300,100) & 0.049 & 0.056 & 0.049 & 0.053\\
(100,200) & 0.048 & 0.054 & 0.048 & 0.051\\
(200,200) & 0.051 & 0.057 & 0.050 & 0.052\\
(300,200) & 0.047 & 0.059 & 0.046 & 0.049\\
(100,300) & 0.051 & 0.056 & 0.051 & 0.055\\
(200,300) & 0.054 & 0.058 & 0.055 & 0.056\\
(300,300) & 0.050 & 0.062 & 0.051 & 0.052\\\hline
& \multicolumn{4}{c}{\underline{$\rho=0.4$}}\\
(100,100) & 0.041 & 0.013 & 0.020 & 1.000\\
(200,100) & 0.039 & 0.006 & 0.013 & 1.000\\
(300,100) & 0.041 & 0.006 & 0.013 & 1.000\\
(100,200) & 0.043 & 0.021 & 0.024 & 1.000\\
(200,200) & 0.044 & 0.010 & 0.016 & 1.000\\
(300,200) & 0.047 & 0.009 & 0.013 & 1.000\\
(100,300) & 0.043 & 0.029 & 0.029 & 1.000\\
(200,300) & 0.051 & 0.018 & 0.023 & 1.000\\
(300,300) & 0.043 & 0.011 & 0.014 & 1.000\\\hline
& \multicolumn{4}{c}{\underline{$\rho=0.7$}}\\
(100,100) & 0.030 & 0.011 & 0.018 & 1.000\\
(200,100) & 0.036 & 0.024 & 0.035 & 1.000\\
(300,100) & 0.033 & 0.080 & 0.103 & 1.000\\
(100,200) & 0.033 & 0.006 & 0.012 & 1.000\\
(200,200) & 0.040 & 0.007 & 0.013 & 1.000\\
(300,200) & 0.039 & 0.007 & 0.014 & 1.000\\
(100,300) & 0.038 & 0.007 & 0.014 & 1.000\\
(200,300) & 0.046 & 0.005 & 0.011 & 1.000\\
(300,300) & 0.041 & 0.005 & 0.011 & 1.000\\\hline\hline
\end{tabular}
\end{center}
\medskip
\noindent\textit{Notes:} $LM_{\inf}$ denotes the oracle LM test using the true
spectral moments. $LM_{UT}$ and $LM_T$ denote the feasible LM tests using
untapered and tapered autocovariance estimators, respectively. $LM_S$ denotes
Schott's (2005) test. All results are based on 5,000 Monte Carlo replications.
\newpage
\begin{center}
\refstepcounter{table}\label{tab:hybrid-size} \noindent\textbf{Table
\thetable. Empirical size of hybrid LM tests at the 5\% nominal level
\newline($\rho=0.4$, $N(0,1)$)}
\FloatBarrier
\centering
\medskip
{\setlength{\tabcolsep}{5pt} \renewcommand{\arraystretch}{1.1} }
\end{center}
{
\begin{tabular}
[c]{c|ccc|cc|cc}\hline\hline
& \multicolumn{3}{c|}{LM tests} & \multicolumn{2}{c|}{Tapered hybrid} &
\multicolumn{2}{c}{Untapered hybrid}\\
$(n,T)$ & $LM_{\inf}$ & $LM_{UT}$ & $LM_{T}$ & $LM_{T,\mu}$ & $LM_{T,\sigma}$
& $LM_{UT,\mu}$ & $LM_{UT,\sigma}$\\\hline
(100,100) & 0.040 & 0.008 & 0.014 & 0.039 & 0.013 & 0.036 & 0.009\\
(200,100) & 0.038 & 0.006 & 0.010 & 0.041 & 0.009 & 0.039 & 0.006\\
(300,100) & 0.038 & 0.003 & 0.007 & 0.039 & 0.007 & 0.039 & 0.004\\
(100,200) & 0.040 & 0.016 & 0.021 & 0.040 & 0.020 & 0.038 & 0.016\\
(200,200) & 0.038 & 0.008 & 0.012 & 0.040 & 0.012 & 0.038 & 0.009\\
(300,200) & 0.043 & 0.006 & 0.011 & 0.043 & 0.010 & 0.041 & 0.006\\
(100,300) & 0.047 & 0.028 & 0.028 & 0.046 & 0.028 & 0.045 & 0.029\\
(200,300) & 0.039 & 0.013 & 0.015 & 0.040 & 0.015 & 0.038 & 0.013\\
(300,300) & 0.044 & 0.012 & 0.018 & 0.045 & 0.018 & 0.044 &
0.013\\\hline\hline
\end{tabular}
}
\begin{center}
\medskip
\parbox{\textwidth}{
\noindent\textit{Notes:}
$LM_{\inf}$ denotes the oracle LM test using the true
spectral moments. $LM_{UT}$ and $LM_T$ denote the
feasible LM tests using untapered and tapered
autocovariance estimators, respectively.
The hybrid tests $LM_{T,\mu}$ and $LM_{UT,\mu}$
use the true centering term and the estimated scaling
term, whereas $LM_{T,\sigma}$ and $LM_{UT,\sigma}$
use the estimated centering term and the true scaling
term. The subscripts $T$ and $UT$ indicate tapered
and untapered estimation, respectively.
}
\end{center}
\newpage
\begin{center}
\refstepcounter{table}\label{tab:power-independent} \noindent\textbf{Table
\thetable. Size-unadjusted empirical power at the 5\% nominal level
($\rho=0.0$)}
\medskip
{\setlength{\tabcolsep}{5pt} \renewcommand{\arraystretch}{1.1} }
\end{center}
{
\begin{tabular}
[c]{c|cccc|cccc}\hline\hline
& \multicolumn{8}{c}{Panel A: Factor model}\\\hline
& \multicolumn{4}{c|}{$\kappa=0.2$} & \multicolumn{4}{c}{$\kappa=0.4$}\\
$(n,T)$ & $LM_{\inf}$ & $LM_{UT}$ & $LM_{T}$ & $LM_{S}$ & $LM_{\inf}$ &
$LM_{UT}$ & $LM_{T}$ & $LM_{S}$\\\hline
(100,100) & 0.154 & 0.098 & 0.147 & 0.161 & 0.999 & 0.998 & 0.999 & 0.999\\
(200,100) & 0.368 & 0.265 & 0.357 & 0.379 & 1.000 & 1.000 & 1.000 & 1.000\\
(300,100) & 0.611 & 0.503 & 0.597 & 0.622 & 1.000 & 1.000 & 1.000 & 1.000\\
(100,200) & 0.390 & 0.288 & 0.384 & 0.397 & 1.000 & 1.000 & 1.000 & 1.000\\
(200,200) & 0.812 & 0.742 & 0.806 & 0.816 & 1.000 & 1.000 & 1.000 & 1.000\\
(300,200) & 0.966 & 0.944 & 0.965 & 0.967 & 1.000 & 1.000 & 1.000 & 1.000\\
(100,300) & 0.651 & 0.547 & 0.648 & 0.658 & 1.000 & 1.000 & 1.000 & 1.000\\
(200,300) & 0.979 & 0.963 & 0.978 & 0.979 & 1.000 & 1.000 & 1.000 & 1.000\\
(300,300) & 0.999 & 0.999 & 0.999 & 0.999 & 1.000 & 1.000 & 1.000 &
1.000\\\hline
& \multicolumn{8}{c}{Panel B: SAR(1) model}\\\hline
& \multicolumn{4}{c|}{$\psi=0.2$} & \multicolumn{4}{c}{$\psi=0.4$}\\
$(n,T)$ & $LM_{\inf}$ & $LM_{UT}$ & $LM_{T}$ & $LM_{S}$ & $LM_{\inf}$ &
$LM_{UT}$ & $LM_{T}$ & $LM_{S}$\\\hline
(100,100) & 0.633 & 0.503 & 0.621 & 0.645 & 1.000 & 1.000 & 1.000 & 1.000\\
(196,100) & 0.594 & 0.458 & 0.577 & 0.606 & 1.000 & 1.000 & 1.000 & 1.000\\
(289,100) & 0.593 & 0.465 & 0.578 & 0.605 & 1.000 & 1.000 & 1.000 & 1.000\\
(100,200) & 0.995 & 0.987 & 0.995 & 0.995 & 1.000 & 1.000 & 1.000 & 1.000\\
(196,200) & 0.993 & 0.984 & 0.993 & 0.993 & 1.000 & 1.000 & 1.000 & 1.000\\
(289,200) & 0.992 & 0.980 & 0.992 & 0.992 & 1.000 & 1.000 & 1.000 & 1.000\\
(100,300) & 1.000 & 1.000 & 1.000 & 1.000 & 1.000 & 1.000 & 1.000 & 1.000\\
(196,300) & 1.000 & 1.000 & 1.000 & 1.000 & 1.000 & 1.000 & 1.000 & 1.000\\
(289,300) & 1.000 & 1.000 & 1.000 & 1.000 & 1.000 & 1.000 & 1.000 &
1.000\\\hline\hline
\end{tabular}
}
\begin{center}
\medskip
\parbox{\textwidth}{
\noindent\textit{Notes:}
$LM_{\inf}$ denotes the oracle LM test.
$LM_{UT}$ and $LM_T$ denote the feasible LM tests using
untapered and tapered autocovariance estimators,
respectively. $LM_S$ denotes Schott's (2005) test.
The parameters $\kappa$ and $\psi$ control the strength
of cross-sectional dependence in the factor and
SAR(1) models, respectively.
}
\end{center}
\FloatBarrier
\newpage
\begin{center}
\refstepcounter{table}\label{tab:power-dependent} \noindent\textbf{Table
\thetable. Size-unadjusted empirical power at the 5\% nominal level
($\rho=0.4$)}
\medskip
{\setlength{\tabcolsep}{6pt} \renewcommand{\arraystretch}{1.1} }
\end{center}
{
\begin{tabular}
[c]{c|ccc|ccc}\hline\hline
& \multicolumn{6}{c}{Panel A: Factor model}\\\hline
& \multicolumn{3}{c|}{$\kappa=0.2$} & \multicolumn{3}{c}{$\kappa=0.4$}\\
$(n,T)$ & $LM_{\inf}$ & $LM_{UT}$ & $LM_{T}$ & $LM_{\inf}$ & $LM_{UT}$ &
$LM_{T}$\\\hline
(100,100) & 0.054 & 0.038 & 0.064 & 0.962 & 0.990 & 0.995\\
(200,100) & 0.068 & 0.098 & 0.146 & 0.999 & 1.000 & 1.000\\
(300,100) & 0.074 & 0.225 & 0.286 & 1.000 & 1.000 & 1.000\\
(100,200) & 0.126 & 0.138 & 0.210 & 1.000 & 1.000 & 1.000\\
(200,200) & 0.284 & 0.449 & 0.546 & 1.000 & 1.000 & 1.000\\
(300,200) & 0.439 & 0.749 & 0.806 & 1.000 & 1.000 & 1.000\\
(100,300) & 0.287 & 0.315 & 0.421 & 1.000 & 1.000 & 1.000\\
(200,300) & 0.647 & 0.802 & 0.866 & 1.000 & 1.000 & 1.000\\
(300,300) & 0.854 & 0.971 & 0.979 & 1.000 & 1.000 & 1.000\\\hline
& \multicolumn{6}{c}{Panel B: SAR(1) model}\\\hline
& \multicolumn{3}{c|}{$\psi=0.2$} & \multicolumn{3}{c}{$\psi=0.4$}\\
$(n,T)$ & $LM_{\inf}$ & $LM_{UT}$ & $LM_{T}$ & $LM_{\inf}$ & $LM_{UT}$ &
$LM_{T}$\\\hline
(100,100) & 0.263 & 0.207 & 0.298 & 1.000 & 1.000 & 1.000\\
(196,100) & 0.191 & 0.162 & 0.232 & 0.998 & 1.000 & 1.000\\
(289,100) & 0.166 & 0.157 & 0.215 & 0.989 & 1.000 & 1.000\\
(100,200) & 0.823 & 0.813 & 0.886 & 1.000 & 1.000 & 1.000\\
(196,200) & 0.717 & 0.764 & 0.843 & 1.000 & 1.000 & 1.000\\
(289,200) & 0.650 & 0.737 & 0.817 & 1.000 & 1.000 & 1.000\\
(100,300) & 0.991 & 0.992 & 0.996 & 1.000 & 1.000 & 1.000\\
(196,300) & 0.984 & 0.993 & 0.996 & 1.000 & 1.000 & 1.000\\
(289,300) & 0.968 & 0.988 & 0.994 & 1.000 & 1.000 & 1.000\\\hline\hline
\end{tabular}
}
\begin{center}
\medskip
\parbox{\textwidth}{
\noindent\textit{Notes:}
$LM_{\inf}$ denotes the oracle LM test.
$LM_{UT}$ and $LM_T$ denote the feasible LM tests using
untapered and tapered autocovariance estimators,
respectively. The parameters $\kappa$ and $\psi$
control the strength of cross-sectional dependence
in the factor and SAR(1) models, respectively.
Schott's test is excluded because it is not robust
to temporal dependence.
}
\end{center}
\FloatBarrier
\newpage
%\begin{figure}[htbp]
%\centering
%\includegraphics[width=0.75\textwidth]{img/ERP_Bandwidth.png}
%\caption{Empirical rejection probabilities across bandwidths
%at the 5\% nominal level ($n=T=200$)}
%figure label omitted (figure disabled)
%\end{figure}
%\FloatBarrier
%\newpage
\end{singlespace}

\pagebreak\appendix{} \numberwithin{equation}{section} \numberwithin{lemma}{section}\numberwithin{theorem}{section}

\noindent This Appendix proves the null results: Lemma~2
(Appendix~\ref{app:lemma2}), Lemma~3 (Appendix~\ref{app:proof-variance}), and
Theorem~1 (Appendix~\ref{app:thm1}). The model, notation, Gaussian moment
identities, and preliminary lemmas used here (Lemmas~S.1--S.5) are in
Section~S.A of the Online Supplement, which also proves Lemmas~1 and~4. The
proofs of the power results (Lemmas~5 and~6, Theorem~2, and Corollary~1) and
the auxiliary inequalities are in the Online Supplement, whose sections,
lemmas, and equations carry the prefix ``S.''.

\section{Proof of Lemma 2}

\label{app:lemma2}\label{app:proof-mean} Under $H_{0}$ the rows of $Z$ are
i.i.d., so all off-diagonal pairs are identically distributed; since
$r_{ii}=1$,
\begin{equation}
\mathbb{E}\left[  Tr(R^{2})\right]  =n+n(n-1)\mathbb{E}(r_{12}^{2}%
).\label{eq:mean-decomp-new}%
\end{equation}
Throughout, $S_{12}=T^{-1}Z_{1}^{\prime}\Gamma Z_{2}$, $\delta_{i}=S_{ii}-1$,
$P_{1}=Z_{1}^{\prime}\Gamma^{2}Z_{1}$, $d_{i}=T\delta_{i}=Z_{i}^{\prime}\Gamma
Z_{i}-T$, and $t_{k}=Tr(\Gamma^{k})=T\tau_{k,T}$, so that $t_{k}=O(T)$ for
each fixed $k$. Let $\mathcal{E}_{12}=\{|\delta_{1}|\leq1/2,\ |\delta_{2}%
|\leq1/2\}$.

\emph{Step 1: Expansion.} On $\mathcal{E}_{12}$, Taylor's formula for
$(1+x)^{-1}$ gives
\[
\frac{1}{(1+\delta_{1})(1+\delta_{2})}=\Pi_{3}(\delta_{1},\delta
_{2})+\mathcal{Q}_{4},\qquad|\mathcal{Q}_{4}|\leq C(|\delta_{1}|+|\delta
_{2}|)^{4},
\]
where
\[
\Pi_{3}(\delta_{1},\delta_{2})=1-(\delta_{1}+\delta_{2})+(\delta_{1}%
^{2}+\delta_{1}\delta_{2}+\delta_{2}^{2})-(\delta_{1}^{3}+\delta_{1}^{2}%
\delta_{2}+\delta_{1}\delta_{2}^{2}+\delta_{2}^{3}).
\]
We retain the cubic terms because bounding them in absolute value gives only
$O(T^{-5/2})$ after multiplication by $S_{12}^{2}$, whereas their expectations
are $O(T^{-3})$ (Step~2). Since $r_{12}^{2}=S_{12}^{2}\Pi_{3}+S_{12}%
^{2}\mathcal{Q}_{4}$ on $\mathcal{E}_{12}$,
\[
\mathbb{E}(r_{12}^{2})=\mathbb{E}(S_{12}^{2}\Pi_{3})-\mathbb{E}(S_{12}^{2}%
\Pi_{3}1_{\mathcal{E}_{12}^{c}})+\mathbb{E}(S_{12}^{2}\mathcal{Q}%
_{4}1_{\mathcal{E}_{12}})+\mathbb{E}(r_{12}^{2}1_{\mathcal{E}_{12}^{c}}).
\]
By H\"older's inequality and Lemma~S.1, $\mathbb{E}\{S_{12}^{2}(|\delta
_{1}|+|\delta_{2}|)^{4}\}\leq C\Vert S_{12}\Vert_{4}^{2}(\Vert\delta_{1}%
\Vert_{8}+\Vert\delta_{2}\Vert_{8})^{4}=O(T^{-3})$. By the per-unit tail bound
following the proof of Lemma~4, $P(\mathcal{E}_{12}^{c})\leq Ce^{-cT}$. Since
$0\leq r_{12}^{2}\leq1$, $\mathbb{E}(r_{12}^{2}1_{\mathcal{E}_{12}^{c}%
})=O(e^{-cT})$; by Cauchy--Schwarz and $\Vert S_{12}^{2}\Pi_{3}\Vert
_{2}=O(T^{-1})$, $\mathbb{E}(S_{12}^{2}\Pi_{3}1_{\mathcal{E}_{12}^{c}%
})=O(e^{-cT})$ as well. Consequently,
\begin{equation}
\label{eq:mean-rem-bound}\mathbb{E}(r_{12}^{2})=\mathbb{E}(S_{12}^{2}\Pi
_{3})+O(T^{-3}).
\end{equation}

\emph{Step 2: Moments.} Since all matrices involved are powers of $\Gamma$,
the joint cumulants of Gaussian quadratic forms are $\kappa(Z^{\prime}%
A_{1}Z,\ldots,Z^{\prime}A_{r}Z)=2^{r-1}(r-1)!Tr(A_{1}\cdots A_{r})$; for
$r=2,3$ these follow from (S.A.3)--(S.A.4). Conditioning on $Z_{1}$, Lemma~S.4
gives $\mathbb{E}(S_{12}^{2}\mid Z_{1})=P_{1}/T^{2}$, so the monomials in
$\delta_{1}$ alone reduce to $\mathbb{E}(S_{12}^{2}\delta_{1}^{k}%
)=\mathbb{E}(P_{1}d_{1}^{k})/T^{k+2}$. Since $d_{1}$ is centered, the
moment--cumulant formula gives $\mathbb{E}(P_{1}d_{1})=Cov(P_{1},d_{1}%
)=2t_{3}$,
\[
\mathbb{E}(P_{1}d_{1}^{2})=\mathbb{E}(P_{1})Var(d_{1})+\kappa(P_{1}%
,d_{1},d_{1})=2t_{2}^{2}+8t_{4},
\]
and
\begin{align*}
\mathbb{E}(P_{1}d_{1}^{3})  &  =\mathbb{E}(P_{1})\kappa_{3}(d_{1}%
)+3Cov(P_{1},d_{1})Var(d_{1})+\kappa(P_{1},d_{1},d_{1},d_{1})\\
&  =8t_{2}t_{3}+12t_{2}t_{3}+48t_{5}=20t_{2}t_{3}+48t_{5}.
\end{align*}
For the mixed monomials, conditionally on $Z_{1}$, $Cov\{(Z_{2}^{\prime}\Gamma
Z_{1})^{2},Z_{2}^{\prime}\Gamma Z_{2}\mid Z_{1}\}=2Z_{1}^{\prime}\Gamma
^{3}Z_{1}$, so $\mathbb{E}(S_{12}^{2}d_{2}\mid Z_{1})=2Z_{1}^{\prime}%
\Gamma^{3}Z_{1}/T^{2}$; hence $\mathbb{E}(S_{12}^{2}\delta_{1}\delta
_{2})=2\mathbb{E}(Z_{1}^{\prime}\Gamma^{3}Z_{1}d_{1})/T^{4}$ and
$\mathbb{E}(S_{12}^{2}\delta_{1}^{2}\delta_{2})=2\mathbb{E}(Z_{1}^{\prime
}\Gamma^{3}Z_{1}d_{1}^{2})/T^{5}$. By symmetry in the two rows,
\begin{align*}
\mathbb{E}(S_{12}^{2})  &  =\frac{t_{2}}{T^{2}}, & \mathbb{E}(S_{12}^{2}%
\delta_{i})  &  =\frac{2t_{3}}{T^{3}},\\
\mathbb{E}(S_{12}^{2}\delta_{i}^{2})  &  =\frac{2t_{2}^{2}+8t_{4}}{T^{4}}, &
\mathbb{E}(S_{12}^{2}\delta_{1}\delta_{2})  &  =\frac{4t_{4}}{T^{4}},\\
\mathbb{E}(S_{12}^{2}\delta_{i}^{3})  &  =\frac{20t_{2}t_{3}+48t_{5}}{T^{5}},
& \mathbb{E}(S_{12}^{2}\delta_{1}^{2}\delta_{2})  &  =\mathbb{E}(S_{12}%
^{2}\delta_{1}\delta_{2}^{2})=\frac{4t_{2}t_{3}+16t_{5}}{T^{5}}.
\end{align*}
Since $t_{k}=O(T)$, every cubic monomial has expectation $O(T^{-3})$ uniformly.

\emph{Step 3: Assembly.} Collecting the monomials of $\Pi_{3}$ by symmetry,
\begin{align*}
\mathbb{E}(S_{12}^{2}\Pi_{3})  &  =\mathbb{E}(S_{12}^{2})-2\mathbb{E}%
(S_{12}^{2}\delta_{1})+\{2\mathbb{E}(S_{12}^{2}\delta_{1}^{2})+\mathbb{E}%
(S_{12}^{2}\delta_{1}\delta_{2})\}\\
&  \quad-\{2\mathbb{E}(S_{12}^{2}\delta_{1}^{3})+2\mathbb{E}(S_{12}^{2}%
\delta_{1}^{2}\delta_{2})\}\\
&  =\frac{t_{2}}{T^{2}}-\frac{4t_{3}}{T^{3}}+\frac{4t_{2}^{2}+20t_{4}}{T^{4}%
}-\frac{48t_{2}t_{3}+128t_{5}}{T^{5}}.
\end{align*}
Substituting into \eqref{eq:mean-rem-bound},
\[
\mathbb{E}(r_{12}^{2})=\frac{t_{2}}{T^{2}}-\frac{4t_{3}}{T^{3}}+\frac
{4t_{2}^{2}+20t_{4}}{T^{4}}+O(T^{-3})=\frac{\tau_{2,T}}{T}+\frac{4(\tau
_{2,T}^{2}-\tau_{3,T})}{T^{2}}+O(T^{-3}),
\]
because $20t_{4}/T^{4}$ and the cubic contribution are $O(T^{-3})$. Insert
this into \eqref{eq:mean-decomp-new} with $T=n\lambda$: $n(n-1)\tau
_{2,T}/T=(n-1)\tau_{2,T}/\lambda$, $4n(n-1)/T^{2}=4(n-1)/(n\lambda^{2})$, and
$n(n-1)O(T^{-3})=O(n^{-1})$. All constants are uniform under Assumptions~1--3.

\section{Proof of Lemma 3}

\label{app:proof-variance} With $t_{k}$, $d_{i}$, and $\mathcal{E}_{12}$ as in
Appendix~\ref{app:lemma2}, write $h_{ij}=r_{ij}^{2}$, $q(z)=z^{\prime}\Gamma
z$, $P(z)=z^{\prime}\Gamma^{2}z$ (so $P_{i}=P(Z_{i})$), and $g(z)=\mathbb{E}%
(h_{12}\mid Z_{1}=z)$. For a generic $Z\sim N(0,I_{T})$, put $d=q(Z)-T$.

\emph{Step 1: Exact decomposition.} Since $Tr(R^{2})=n+2\sum_{i<j}h_{ij}$,
$Var\{Tr(R^{2})\}=4Var(\sum_{i<j}h_{ij})$. Under $H_{0}$, disjoint pairs are
independent, and two pairs sharing one index are conditionally independent
given the shared row, so that $Cov(h_{12},h_{13})=Var\{g(Z_{1})\}$. There are
$n(n-1)/2$ unordered pairs and $n(n-1)(n-2)$ ordered pairs of distinct
unordered pairs with exactly one common index. Hence
\begin{equation}
\label{eq:var-exact-count}Var\{Tr(R^{2})\} =2n(n-1)a_{T}+4n(n-1)(n-2)b_{T},
\quad a_{T}=Var(h_{12}),\quad b_{T}=Var\{g(Z_{1})\}.
\end{equation}
We show $a_{T}=2\tau_{2,T}^{2}/T^{2}+O(T^{-3})$ and $b_{T}=2K_{T}%
/T^{3}+O(T^{-4})$.

\emph{Step 2: Self-pair term.} Put $A=Z_{1}^{\prime}\Gamma Z_{2}=TS_{12}$, so
that $h_{12}^{2}=T^{-4}A^{4}(1+\delta_{1})^{-2}(1+\delta_{2})^{-2}$. On
$\mathcal{E}_{12}$,
\[
(1+\delta_{1})^{-2}(1+\delta_{2})^{-2}=1-2\delta_{1}-2\delta_{2}%
+\mathcal{R}_{2},\qquad|\mathcal{R}_{2}|\leq C(\delta_{1}^{2}+\delta_{2}%
^{2}).
\]
Since $\mathbb{E}(A^{4}\mid Z_{1})=3P_{1}^{2}$, Lemma~S.2 and the
moment--cumulant formula give, for $i=1,2$ (by symmetry in the two rows),
\[
\mathbb{E}(A^{4})=3(t_{2}^{2}+2t_{4}),\qquad\mathbb{E}(A^{4}d_{i}%
)=3\mathbb{E}(P_{1}^{2}d_{1})=3(4t_{2}t_{3}+8t_{5})=O(T^{2}),
\]
so the linear terms contribute $T^{-5}\mathbb{E}(A^{4}d_{i})=O(T^{-3})$. By
H\"older's inequality and Lemma~S.1, $\mathbb{E}(A^{4}\delta_{i}^{2})\leq\Vert
A\Vert_{8}^{4}\Vert\delta_{i}\Vert_{4}^{2}=O(T)$, so the remainder contributes
$O(T^{-3})$. The exceptional event $\mathcal{E}_{12}^{c}$ contributes
$O(e^{-cT})$, as in Appendix~\ref{app:lemma2}, using $0\leq h_{12}^{2}\leq1$.
Hence $\mathbb{E}(h_{12}^{2})=(3t_{2}^{2}+6t_{4})/T^{4}+O(T^{-3})=3\tau
_{2,T}^{2}/T^{2}+O(T^{-3})$. Lemma~2 gives $\mathbb{E}(h_{12})=\tau
_{2,T}/T+O(T^{-2})$, so $a_{T}=2\tau_{2,T}^{2}/T^{2}+O(T^{-3})$.

\emph{Step 3: Representation of $g$.} Let $\nu_{1,T},\ldots,\nu_{T,T}$ be the
eigenvalues of $\Gamma$ and $m_{T}=\#\{t:\nu_{t,T}\geq1/2\}$. Since
$Tr\Gamma=T$ and $\Vert\Gamma\Vert_{op}\leq K_{\Gamma}$, $T\leq T/2+K_{\Gamma
}m_{T}$, so $m_{T}\geq T/(2K_{\Gamma})$. Hence $q(Z)$ stochastically dominates
$\chi_{m_{T}}^{2}/2$, and for each fixed $p\geq1$ and all $T$ with $m_{T}%
>2p$,
\begin{equation}
\label{eq:q-inverse}\mathbb{E}q(Z)^{-p}=O(T^{-p}).
\end{equation}
Let $B_{T}=\mathbb{E}\{ZZ^{\prime}/q(Z)\}$. Independence of $Z_{1}$ and
$Z_{2}$ gives the exact representation $g(z)=z^{\prime}\Gamma B_{T}\Gamma
z/q(z)$. The identity $T/q=1-d/T+d^{2}/(Tq)$ and $\mathbb{E}(ZZ^{\prime
}d)=2\Gamma$ yield
\[
TB_{T}-I_{T}=-\frac{2\Gamma}{T}+\frac{1}{T}\mathbb{E}\left\{  \frac
{ZZ^{\prime}d^{2}}{q(Z)}\right\}  .
\]
For a unit vector $u$, H\"older's inequality and \eqref{eq:q-inverse} give
\begin{align*}
\mathbb{E}\left\{  \frac{(u^{\prime}Z)^{2}d^{2}}{q(Z)}\right\}   &  \leq\Vert
u^{\prime}Z\Vert_{8}^{2}\Vert d\Vert_{8}^{2} \Vert q(Z)^{-1}\Vert_{2}\\
&  =O(1).
\end{align*}
Since $\mathbb{E}\{ZZ^{\prime}d^{2}/q(Z)\}$ is positive semidefinite, its
operator norm is the supremum of this quantity over unit $u$, and is therefore
$O(1)$. Hence
\begin{equation}
\label{eq:BT-bound}\Vert B_{T}-T^{-1}I_{T}\Vert_{op}=O(T^{-2}),\qquad\text{in
particular}\qquad\Vert B_{T}\Vert_{op}=O(T^{-1}).
\end{equation}
Set $N_{T}=\Gamma(B_{T}-T^{-1}I_{T})\Gamma$ and $\varpi_{T}=Tr(N_{T})/T$. Then
$\Vert N_{T}\Vert_{op}=O(T^{-2})$, $|\varpi_{T}|\leq\Vert N_{T}\Vert_{op}$,
$Tr(N_{T}-\varpi_{T}\Gamma)=0$, and $\Vert N_{T}-\varpi_{T}\Gamma\Vert_{F}\leq
T^{1/2}\Vert N_{T}-\varpi_{T}\Gamma\Vert_{op}=O(T^{-3/2})$. With
$W=P(Z)/\{Tq(Z)\}$,
\[
g(Z)=W+\varpi_{T}+\frac{Z^{\prime}(N_{T}-\varpi_{T}\Gamma)Z}{q(Z)}.
\]
For a trace-zero symmetric $H$, $\mathbb{E}(Z^{\prime}HZ)^{4}\leq C\Vert
H\Vert_{F}^{4}$; Cauchy--Schwarz and \eqref{eq:q-inverse} therefore bound the
$L^{2}$ norm of the last term by $O(T^{-3/2})\cdot O(T^{-1})=O(T^{-5/2})$.
Thus $g(Z)-\mathbb{E}g(Z)=W-\mathbb{E}W+e_{T}$ with $\Vert e_{T}\Vert
_{2}=O(T^{-5/2})$.

\emph{Step 4: Variance of $g$.} Put $D=P(Z)-\tau_{2,T}q(Z)=Z^{\prime}H_{D}Z$
with $H_{D}=\Gamma^{2}-\tau_{2,T}\Gamma$. Then $Tr(H_{D})=0$ and $Tr(H_{D}%
^{2})=TK_{T}$. The exact identity
\[
W-\frac{\tau_{2,T}}{T}=\frac{D}{T^{2}}-\frac{Dd}{T^{3}}+\frac{Dd^{2}}%
{T^{3}q(Z)}%
\]
holds, and by H\"older's inequality and \eqref{eq:q-inverse} the three terms
on the right have $L^{2}$ norms $O(T^{-3/2})$, $O(T^{-2})$, and $O(T^{-5/2})$.
Since $D$ and $d$ are centered, (S.A.3) and (S.A.4) give
\[
\mathbb{E}D^{2}=2TK_{T},\qquad\mathbb{E}(D^{2}d)=8Tr(H_{D}^{2}\Gamma)=8\left(
t_{5}-\frac{2t_{2}t_{4}}{T}+\frac{t_{2}^{2}t_{3}}{T^{2}}\right)  =O(T).
\]
The first term therefore has variance $2K_{T}/T^{3}$; its covariance with the
second term is $-\mathbb{E}(D^{2}d)/T^{5}=O(T^{-4})$; and all remaining
variance and covariance terms are $O(T^{-4})$ by Cauchy--Schwarz. Hence
$Var(W)=2K_{T}/T^{3}+O(T^{-4})$. Finally, $b_{T}=Var(W)+2Cov(W,e_{T}%
)+\mathbb{E}e_{T}^{2}$ and $\Vert W-\mathbb{E}W\Vert_{2}=O(T^{-3/2})$, so
\begin{equation}
\label{eq:bT-expansion}b_{T}=\frac{2K_{T}}{T^{3}}+O(T^{-4}).
\end{equation}
These bounds are absolute and do not require $K_{T}>0$.

\emph{Step 5: Assembly.} Insert both expansions
into~\eqref{eq:var-exact-count}. Because $T=n\lambda$, the aggregate errors
are $O(n^{2}T^{-3}+n^{3}T^{-4})=O(n^{-1})$. Replacing $n(n-1)$ and
$n(n-1)(n-2)$ by $n^{2}$ and $n^{3}$ costs $O(n^{-1})$, since $\tau_{2,T}$ and
$K_{T}$ are uniformly bounded. Hence
\[
Var\{Tr(R^{2})\} =\frac{4\tau_{2,T}^{2}}{\lambda^{2}}+ \frac{8K_{T}}%
{\lambda^{3}}+O(n^{-1}).
\]

\emph{Second-order terms.} The second-order expansions in
Remark~\ref{rem:var-terms} follow from the same argument with longer Taylor
polynomials. For $a_{T}$, expand $(1+\delta_{1})^{-2}(1+\delta_{2})^{-2}$
through total degree three in $(\delta_{1},\delta_{2})$; by H\"older's
inequality the remainder contributes $O(\Vert A\Vert_{8}^{4}\Vert\delta
\Vert_{8}^{4}T^{-4})=O(T^{-4})$. For $b_{T}$, write $b_{T}=\mathbb{E}%
(h_{12}h_{13})-(\mathbb{E}h_{12})^{2}$, expand $(1+\delta_{1})^{-2}%
(1+\delta_{2})^{-1}(1+\delta_{3})^{-1}$ through degree five and $(1+\delta
_{1})^{-1}(1+\delta_{2})^{-1}$ through degree seven; the remainders are
$O(T^{-5})$, and the event $\{\max_{a\leq3}|\delta_{a}|>1/2\}$ contributes
$O(e^{-cT})$ as before. The retained terms are products of centered Gaussian
quadratic forms in $(Z_{1},Z_{2},Z_{3})$. In an eigenbasis of $\Gamma$ these
forms are $T^{-1}\sum_{t}\nu_{t,T}\zeta_{t}$ with i.i.d.\ copies $\zeta_{t}$
of a vector of products of independent standard normals, so every $k$th joint
cumulant equals $T^{1-k}\tau_{k,T}$ times a rational constant. The
moment--cumulant formula therefore gives each expectation as an exact
polynomial in $T^{-1}$ with coefficients polynomial in $\tau_{2,T},\ldots
,\tau_{5,T}$. An expectation of a product of $d$ centered forms is
$O(T^{-\lceil d/2\rceil})$, and retaining those up to $T^{-3}$ for $a_{T}$ and
$T^{-4}$ for $b_{T}$ yields the second-order expansions of $a_{T}$ and $b_{T}%
$. The same computation also gives $\mathbb{E}(r_{12}^{2})=\tau_{2,T}%
/T+4(\tau_{2,T}^{2}-\tau_{3,T})/T^{2} +(28\tau_{2,T}^{3}-48\tau_{2,T}%
\tau_{3,T}+20\tau_{4,T})/T^{3}+O(T^{-4})$. In $b_{T}$ the $T^{-2}$ terms of
$\mathbb{E}(h_{12}h_{13})$ and $(\mathbb{E}h_{12})^{2}$ cancel exactly, in
agreement with Step~4.

\section{Proof of Theorem 1}

\label{app:thm1}

\noindent We first establish the exact Hoeffding decomposition of $Tr(R^{2})$
and moment bounds for its Doob increments
(Lemma~\ref{lem:hoeffding-projection}), and then apply the martingale CLT.
Put
\[
q_{i}=Z_{i}^{\prime}\Gamma Z_{i},\qquad u_{i}=\frac{\Gamma^{1/2}Z_{i}}%
{q_{i}^{1/2}},\qquad h_{ij}=h_{T}(Z_{i},Z_{j})=(u_{i}^{\prime}u_{j}%
)^{2}=r_{ij}^{2}.
\]
The vectors $u_{1},\ldots,u_{n}$ are i.i.d. and $\Vert u_{i}\Vert=1$. Define
\[
M_{T}=\mathbb{E}(u_{1}u_{1}^{\prime}),\qquad\theta_{T}=\mathbb{E}%
(h_{12})=Tr(M_{T}^{2}),
\]
and the first-order and second-order degenerate Hoeffding projections
\begin{align*}
\xi_{i}  &  =\mathbb{E}(h_{ij}\mid Z_{i})-\theta_{T}=u_{i}^{\prime}M_{T}%
u_{i}-\theta_{T},\\
h_{ij}^{\circ}  &  =h_{ij}-\theta_{T}-\xi_{i}-\xi_{j}.
\end{align*}
Then
\[
\mathbb{E}\xi_{i}=0,\qquad\mathbb{E}(h_{ij}^{\circ}\mid Z_{i})=\mathbb{E}%
(h_{ij}^{\circ}\mid Z_{j})=0.
\]

\begin{lemma}
[Hoeffding-projection bounds]\label{lem:hoeffding-projection} Under
Assumptions 1--3 and $H_{0}$,
\begin{gather}
\mathbb{E}\xi_{1}^{2} =\frac{2K_{T}}{T^{3}}+O(T^{-4}),\label{eq:xi-var}\\
\mathbb{E}[(h_{12}^{\circ})^{2}] =\frac{2\tau_{2,T}^{2}}{T^{2}}+O(T^{-3}%
),\label{eq:hdeg-var}\\
\mathbb{E}\xi_{1}^{4} =O(T^{-6}),\qquad\mathbb{E}[(h_{12}^{\circ}%
)^{4}]=O(T^{-4}). \label{eq:projection-fourth}%
\end{gather}
For $\mathcal{F}_{k}=\sigma(Z_{1},\ldots,Z_{k})$, the Doob increment of
$Tr(R^{2})$ is
\begin{equation}
Y_{k}=2(n-1)\xi_{k}+2\sum_{i<k}h_{ik}^{\circ}. \label{eq:full-doob}%
\end{equation}
Moreover,
\begin{gather}
\sum_{k=1}^{n}\mathbb{E}(Y_{k}^{4}) =o(1),\label{eq:fourth-sum}\\
\sum_{k=1}^{n}\mathbb{E}(Y_{k}^{2}\mid\mathcal{F}_{k-1})- \sum_{k=1}%
^{n}\mathbb{E}(Y_{k}^{2}) \overset{p}{\longrightarrow}0.
\label{eq:pqv-concentration}%
\end{gather}
Consequently,
\[
\max_{1\leq k\leq n}|Y_{k}|\overset{p}{\longrightarrow}0, \qquad\max_{1\leq
k\leq n}\frac{|Y_{k}|}{\sigma_{n}} \overset{p}{\longrightarrow}0,
\]
where $\sigma_{n}^{2}$ is given in \eqref{eq:mu-sigma}.
\end{lemma}

\begin{proof}
\emph{Second moments.} For $Z\sim N(0,I_{T})$ let $q=Z^{\prime}\Gamma Z$ and
$\Psi=\Gamma^{1/2}ZZ^{\prime}\Gamma^{1/2}/q$, so that $\Psi_{i}=u_{i}%
u_{i}^{\prime}$ and $M_{T}=\mathbb{E}\Psi=\Gamma^{1/2}B_{T}\Gamma^{1/2}$ with
$B_{T}$ as in Appendix~\ref{app:proof-variance}. Since $\xi_{1}=\langle
\Psi_{1},M_{T}\rangle_{F}-\theta_{T}$ is the centered conditional mean
$g(Z_{1})-\mathbb{E}g(Z_{1})$ in the proof of Lemma~3
(Appendix~\ref{app:proof-variance}), $\mathbb{E}\xi_{1}^{2}=b_{T}$, and
\eqref{eq:bT-expansion} gives \eqref{eq:xi-var}. Orthogonality of the
Hoeffding components gives $\mathbb{E}[(h_{12}^{\circ})^{2}]=Var(h_{12}%
)-2\mathbb{E}\xi_{1}^{2}=a_{T}-2b_{T}$, which is \eqref{eq:hdeg-var} by Lemma~3.

\emph{A uniform fourth-moment bound.} For a deterministic symmetric $H$, put
$Q_{H}=\Gamma^{1/2}H\Gamma^{1/2}$, $c_{H}=Tr(Q_{H})/T$, and $B_{H}=Q_{H}%
-c_{H}\Gamma$. Then
\[
\langle\Psi,H\rangle_{F}-c_{H}=\frac{Z^{\prime}B_{H}Z}{q},\qquad
Tr(B_{H})=0,\qquad\Vert B_{H}\Vert_{F}\leq C\Vert H\Vert_{F},
\]
because $|c_{H}|\leq\Vert H\Vert_{F}\Vert\Gamma\Vert_{F}/T$ and $\Vert
\Gamma\Vert_{F}\leq K_{\Gamma}T^{1/2}$. Gaussian quadratic-form moments give
$\mathbb{E}(Z^{\prime}B_{H}Z)^{8}\leq C\Vert B_{H}\Vert_{F}^{8}$, and
\eqref{eq:q-inverse} gives $\mathbb{E}q^{-8}=O(T^{-8})$. By Cauchy--Schwarz,
and since centering at the mean instead of $c_{H}$ costs at most a factor
$2^{4}$,
\begin{equation}
\mathbb{E}|\langle\Psi,H\rangle_{F}-\mathbb{E}\langle\Psi,H\rangle_{F}%
|^{4}\leq C\Vert B_{H}\Vert_{F}^{4}T^{-4}. \label{eq:psi-fourth}%
\end{equation}
In particular, by Lyapunov's inequality, $Var\langle\Psi,H\rangle_{F}\leq
CT^{-2}\Vert H\Vert_{F}^{2}$.

\emph{Fourth moments.} By \eqref{eq:BT-bound}, $\Vert M_{T}\Vert_{op}%
=O(T^{-1})$, so $\Vert Q_{M_{T}}\Vert_{F}=O(T^{-1/2})$, $|c_{M_{T}}%
|=O(T^{-1})$, and $\Vert B_{M_{T}}\Vert_{F}=O(T^{-1/2})$. Applying
\eqref{eq:psi-fourth} with $H=M_{T}$, now with $\Vert B_{M_{T}}\Vert_{F}%
^{4}=O(T^{-2})$, gives $\mathbb{E}\xi_{1}^{4}=O(T^{-6})$. Next, write
$h_{12}=\ell^{2}/q_{2}$ with $\ell=u_{1}^{\prime}\Gamma^{1/2}Z_{2}$.
Conditionally on $u_{1}$, $\ell$ is centered normal with variance
$u_{1}^{\prime}\Gamma u_{1}\leq K_{\Gamma}$, so $\mathbb{E}(\ell^{16}\mid
u_{1})\leq C$, and Cauchy--Schwarz with $\mathbb{E}q_{2}^{-8}=O(T^{-8})$
gives
\[
\mathbb{E}(h_{12}^{4}\mid u_{1})\leq\{\mathbb{E}(\ell^{16}\mid u_{1}%
)\}^{1/2}\{Eq_{2}^{-8}\}^{1/2}=O(T^{-4}).
\]
Since $h_{12}^{\circ}=h_{12}-\theta_{T}-\xi_{1}-\xi_{2}$ and $\theta
_{T}=O(T^{-1})$, $\mathbb{E}(h_{12}^{\circ})^{4}=O(T^{-4})$. This proves \eqref{eq:projection-fourth}.

\emph{Doob increments.} The exact Hoeffding identity is
\begin{equation}
2\sum_{i<j}(h_{ij}-\theta_{T}) =2(n-1)\sum_{i=1}^{n}\xi_{i}+2\sum_{i<j}%
h_{ij}^{\circ}. \label{eq:hoeffding-exact}%
\end{equation}
Its left side is $Tr(R^{2})-\mathbb{E} Tr(R^{2})$. Taking the difference of
its conditional expectations given $\mathcal{F}_{k}$ and $\mathcal{F}_{k-1}$,
the linear part contributes $2(n-1)\xi_{k}$, and, by degeneracy, the only
surviving degenerate terms are $2\sum_{i<k}h_{ik}^{\circ}$. This proves \eqref{eq:full-doob}.

Write $\Psi_{i}=u_{i}u_{i}^{\prime}$, $A_{i}=\Psi_{i}-M_{T}$, and $\bar
{M}=(n-1)M_{T}$. Since $\xi_{k}=\langle M_{T},A_{k}\rangle_{F}$ and
$h_{ik}^{\circ}=\langle A_{i},A_{k}\rangle_{F}$, the increment
\eqref{eq:full-doob} is $Y_{k}=2\langle\bar{M}+\sum_{i<k}A_{i},A_{k}%
\rangle_{F}$. Conditional on $A_{k}$, the random variables $X_{ik}=\langle
A_{i},A_{k}\rangle_{F}$ ($i<k$) are independent and centered. Uniformly in
$u_{k}$,
\[
\mathbb{E}(X_{ik}^{2}\mid A_{k})\leq CT^{-2},\qquad\mathbb{E}(X_{ik}^{4}\mid
A_{k})\leq CT^{-4}.
\]
These bounds hold because, conditional on $u_{k}=v$, $X_{ik}=\langle\Psi
_{i},H_{v}\rangle_{F}-\mathbb{E}\langle\Psi_{i},H_{v}\rangle_{F}$ with
$H_{v}=vv^{\prime}-M_{T}$, and $\|H_{v}\|_{F}\leq2$ uniformly in $v$, so
\eqref{eq:psi-fourth} applies. Moreover, $\langle\bar{M},A_{k}\rangle
_{F}=(n-1)\xi_{k}$ and $\mathbb{E}\xi_{k}^{4}=O(T^{-6})$ by
\eqref{eq:projection-fourth}. Since $\mathbb{E}\{(\sum_{i<k}X_{ik})^{4}\mid
A_{k}\}=\sum_{i<k}\mathbb{E}(X_{ik}^{4}\mid A_{k})+3\sum_{i\neq j<k}%
\mathbb{E}(X_{ik}^{2}\mid A_{k})\mathbb{E}(X_{jk}^{2}\mid A_{k})$ and
$(x+y)^{4}\leq8(x^{4}+y^{4})$,
\[
\mathbb{E} Y_{k}^{4}\leq C\left\{  n^{4}T^{-6}+kT^{-4}+k^{2}T^{-4}\right\}
\leq C\{n^{4}T^{-6}+n^{2}T^{-4}\}=O(n^{-2}).
\]
Summing over $k\leq n$ proves \eqref{eq:fourth-sum}.

We next verify predictable quadratic-variation convergence. Define the
covariance operator $\mathcal{C}_{T}(H)=\mathbb{E}\{\langle A_{1},H\rangle
_{F}A_{1}\}$ on the real Hilbert space of symmetric $T\times T$ matrices with
the Frobenius inner product; operator traces $Tr_{op}$ are taken on this
space. The operator $\mathcal{C}_{T}$ is positive and self-adjoint, and its
Rayleigh quotient at $H$ is $Var\langle\Psi,H\rangle_{F}$. By
\eqref{eq:psi-fourth} and $Tr_{op}\mathcal{C}_{T}=\mathbb{E}\Vert A_{1}%
\Vert_{F}^{2}=1-Tr(M_{T}^{2})\leq1$, the eigenvalues $d_{j}\geq0$ of
$\mathcal{C}_{T}$ satisfy $\max_{j}d_{j}\leq CT^{-2}$ and $\sum_{j}d_{j}\leq
1$. Hence
\begin{gather*}
\|\mathcal{C}_{T}\|_{op}\leq CT^{-2},\\
Tr_{op}(\mathcal{C}_{T}^{2})\leq(\max_{j}d_{j})\sum_{j}d_{j}\leq
CT^{-2},\qquad Tr_{op}(\mathcal{C}_{T}^{4})\leq(\max_{j}d_{j})^{3}\sum
_{j}d_{j}\leq CT^{-6}.
\end{gather*}

Since $Y_{k}=2\langle\bar{M}+\sum_{i<k}A_{i},A_{k}\rangle_{F}$ and $A_{k}$ is
independent of $\mathcal{F}_{k-1}$,
\[
\mathbb{E}(Y_{k}^{2}\mid\mathcal{F}_{k-1}) =4\left\langle \bar{M}+\sum
_{i<k}A_{i}, \mathcal{C}_{T}\left(  \bar{M}+\sum_{i<k}A_{i}\right)
\right\rangle _{F}.
\]
After summing over $k$ and subtracting expectations, the fluctuation
$\mathcal{V}_{n}=\sum_{k}\{\mathbb{E}(Y_{k}^{2}\mid\mathcal{F}_{k-1}%
)-\mathbb{E}Y_{k}^{2}\}$ is exactly
\begin{align*}
\mathcal{V}_{n}  &  =8\sum_{i=1}^{n-1}(n-i)\langle\bar{M},\mathcal{C}_{T}
A_{i}\rangle_{F}\\
&  \quad+4\sum_{i=1}^{n-1}(n-i)\{\langle A_{i},\mathcal{C}_{T}A_{i}\rangle_{F}
-\mathbb{E}\langle A_{i},\mathcal{C}_{T}A_{i}\rangle_{F}\}\\
&  \quad+8\sum_{i<j\leq n-1}(n-j)\langle A_{i},\mathcal{C}_{T}A_{j}\rangle
_{F}.
\end{align*}
All three sums are centered. The third is orthogonal in $L^{2}$ to the first
two, and its distinct unordered pairs are mutually orthogonal, including pairs
sharing one index, by conditioning on the unshared index. Each of the first
two sums has independent summands. Hence $\mathbb{E}\mathcal{V}_{n}^{2}$ is at
most a constant times the sum of the three variances. Since $\|M_{T}\|_{F}%
^{2}=Tr(M_{T}^{2}) =O(T^{-1})$ and $\mathbb{E}\langle\bar{M},\mathcal{C}%
_{T}A_{i}\rangle_{F}^{2}=\langle\mathcal{C}_{T}\bar{M},\mathcal{C}%
_{T}(\mathcal{C}_{T}\bar{M})\rangle_{F}=\langle\bar{M},\mathcal{C}_{T}^{3}%
\bar{M}\rangle_{F}$, the variance of the first sum is at most
\[
C n^{3}\langle\bar{M},\mathcal{C}_{T}^{3}\bar{M}\rangle_{F} \leq C
n^{3}\|\mathcal{C}_{T}\|_{op}^{3}\|\bar{M}\|_{F}^{2} =O(T^{-2}).
\]
For the second sum, $\|A_{i}\|_{F}\leq2$ and $|\langle A_{i},\mathcal{C}_{T}
A_{i}\rangle_{F}|\leq4\|\mathcal{C}_{T}\|_{op}$, giving variance at most
$Cn^{3}T^{-4}=O(T^{-1})$. For the third, distinct unordered pairs are
orthogonal and
\[
\mathbb{E}\langle A_{i},\mathcal{C}_{T}A_{j}\rangle_{F}^{2}=Tr_{op}%
(\mathcal{C}_{T}^{4})=O(T^{-6}),
\]
so its variance is bounded by $Cn^{4}T^{-6}=O(T^{-2})$. Consequently,
\[
\mathbb{E}\mathcal{V}_{n}^{2}=O(T^{-1})=o(1),
\]
which proves \eqref{eq:pqv-concentration}.

The fourth-moment bound also gives the maximal-increment assertion:
\[
P\left(  \max_{k\leq n}|Y_{k}|>\varepsilon\right)  \leq\varepsilon^{-4}%
\sum_{k=1}^{n}\mathbb{E}(Y_{k}^{4})=o(1).
\]
The normalized assertion follows because $\sigma_{n}^{2}\geq4/\lambda^{2}$ is
bounded away from zero as $\lambda\to\lambda_{0}\in(0,\infty)$.
\end{proof}

\medskip\noindent\textbf{Proof of Theorem 1.} By \eqref{eq:hoeffding-exact},
\[
Tr(R^{2})-\mathbb{E}[Tr(R^{2})]=2(n-1)\sum_{i=1}^{n}\xi_{i}+2\sum_{i<j}%
h_{ij}^{\circ}=\sum_{k=1}^{n}Y_{k}.
\]
Orthogonality of the Hoeffding components and
\eqref{eq:xi-var}--\eqref{eq:hdeg-var} give
\begin{align}
\sum_{k=1}^{n}\mathbb{E}(Y_{k}^{2})  &  =4n(n-1)^{2}\mathbb{E}\xi_{1}%
^{2}+4\binom{n}{2}\mathbb{E}[(h_{12}^{\circ})^{2}]\nonumber\\
&  =\frac{4n(n-1)\tau_{2,T}^{2}}{T^{2}}+\frac{8n(n-1)^{2}K_{T}}{T^{3}%
}+O(n^{-1})\nonumber\\
&  =\frac{4\tau_{2,T}^{2}}{\lambda^{2}}+\frac{8K_{T}}{\lambda^{3}}%
+O(n^{-1})=\sigma_{n}^{2}+O(n^{-1}), \label{eq:hoeffding-variance}%
\end{align}
since the projection remainders contribute $O(n^{3}T^{-4}+n^{2}T^{-3}%
)=O(n^{-1})$, and replacing $n(n-1)^{2}$ and $n(n-1)$ by $n^{3}$ and $n^{2}$
costs $O(n^{-1})$. (The first line equals $2n(n-1)a_{T}+4n(n-1)(n-2)b_{T}$, in
agreement with \eqref{eq:var-exact-count}.) Equations
\eqref{eq:pqv-concentration} and \eqref{eq:hoeffding-variance} imply
\[
\frac{1}{\sigma_{n}^{2}}\sum_{k=1}^{n}\mathbb{E}(Y_{k}^{2}\mid\mathcal{F}%
_{k-1})\overset{p}{\longrightarrow}1.
\]
For every $\varepsilon>0$,
\[
\sum_{k=1}^{n}\mathbb{E}\left[  \frac{Y_{k}^{2}}{\sigma_{n}^{2}}%
\mathbf{1}_{\{|Y_{k}|>\varepsilon\sigma_{n}\}}|\mathcal{F}_{k-1}\right]
\leq\frac{1}{\varepsilon^{2}\sigma_{n}^{4}}\sum_{k=1}^{n}\mathbb{E}(Y_{k}%
^{4}\mid\mathcal{F}_{k-1})\overset{p}{\longrightarrow}0,
\]
because the expectation of the last sum tends to zero by
\eqref{eq:fourth-sum}. Thus the conditional Lindeberg condition holds in
probability. Together with predictable quadratic-variation convergence to the
constant $1$, the Hall and Heyde (1980, Cor.~3.1) martingale CLT yields
$\sum_{k}Y_{k}/\sigma_{n}\overset{d}{\rightarrow}N(0,1)$, that is,
$(Tr(R^{2})-\mathbb{E}[Tr(R^{2})])/\sigma_{n}\overset{d}{\rightarrow}N(0,1)$.
Because the limit of the conditional variance is constant, the nesting
condition on the $\sigma$-fields in that corollary is not needed, so the array
filtrations $\mathcal{F}_{k}=\mathcal{F}_{n,k}$ are admissible. Finally,
$\sigma_{n}\geq2/\lambda$ and Lemma~2 give $(\mathbb{E}[Tr(R^{2})]-\mu
_{n})/\sigma_{n}=O(n^{-1})$, and Slutsky's theorem yields
$LM\overset{d}{\rightarrow}N(0,1)$.

\clearpage
\setcounter{section}{0}
\setcounter{equation}{0}
\setcounter{lemma}{0}
\setcounter{theorem}{0}
\setcounter{corollary}{0}
\renewcommand{\thesection}{S.\Alph{section}}
\renewcommand{\theequation}{\thesection.\arabic{equation}}
\renewcommand{\thelemma}{S.\arabic{lemma}}
\renewcommand{\thecorollary}{S.\arabic{corollary}}
\renewcommand{\thetheorem}{S.\arabic{theorem}}

\begin{center}
{\Large Online Supplement to\\[0.3em]``A High-Dimensional LM Test for
Cross-Sectional Independence with Serially Correlated Data''}

\bigskip Ulrich Hounyo, Chihwa Kao, and Min Seong Kim
\end{center}

\renewcommand{\thesection}{S.\Alph{section}} \numberwithin{equation}{section}
\renewcommand{\thelemma}{S.\arabic{lemma}}
\renewcommand{\thecorollary}{S.\arabic{corollary}} \renewcommand{\thetheorem}{S.\arabic{theorem}}

\noindent This Online Supplement contains the proofs not given in the Appendix
of the main paper. Section~\ref{app:prelim} collects the model, notation,
Gaussian moment identities, and the expansion of the diagonal normalization
used throughout, including in the Appendix of the main paper;
Section~\ref{app:diag} also proves Lemmas~1 and~4, which are special cases of
results in Hounyo, Kao, and Kim (2026a). Sections~\ref{proof:mean_alt}%
--\ref{proof:boundary} prove the power results: Lemma~5
(Section~\ref{proof:mean_alt}), Lemma~6 (Section~\ref{proof:var_alt}),
Theorem~2 (Section~\ref{proof:clt_alt}), and Corollary~1
(Section~\ref{proof:boundary}). Section~\ref{app:aux-inequalities} collects
the auxiliary Gaussian and operator inequalities. Unprefixed theorem, lemma,
assumption, section, and equation numbers refer to the main paper.

\section{Setup, Notation, and Preliminary Results}

\label{app:prelim}

\subsection{Model, normalization, and notation}

\noindent Throughout the theoretical analysis, except in
Sections~\ref{proof:mean_alt}--\ref{proof:boundary}, which consider
alternatives and use separate notation, we work under the null hypothesis
$H_{0}:U=I_{n}$. Under $H_{0}$, the separable Gaussian model reduces to
$X=Z\Gamma^{1/2}$, where $Z\in\mathbb{R}^{n\times T}$ has i.i.d.\ $N(0,1)$
entries and $\Gamma$ satisfies Assumption~2, with the normalization
$\tau_{1,T}=T^{-1}Tr(\Gamma)=1$ imposed in Section~2. Let $Z_{i}\in
\mathbb{R}^{T}$ denote the transpose of the $i$th row of $Z$. The sample
covariance matrix and its entries are
\[
S=\frac{1}{T}Z\Gamma Z^{\prime},\qquad S_{ij}=\frac{1}{T}Z_{i}^{\prime}\Gamma
Z_{j},
\]
and the corresponding sample correlations satisfy
\[
r_{ij}=\frac{S_{ij}}{\sqrt{S_{ii}S_{jj}}} =\frac{Z_{i}^{\prime}\Gamma Z_{j}}{
\sqrt{(Z_{i}^{\prime}\Gamma Z_{i})(Z_{j}^{\prime}\Gamma Z_{j})}}.
\]
By the normalization of $\Gamma$, $\mathbb{E}(S_{ii})=T^{-1}Tr(\Gamma)=1$.
Under the null, the cross-sectional rows are independent, while temporal
dependence is retained through $\Gamma$. We further define
\[
\delta_{i}=S_{ii}-1,\qquad\Vert\delta\Vert_{\infty}=\max_{1\leq i\leq
n}|\delta_{i}|,\qquad P_{i}=Z_{i}^{\prime}\Gamma^{2}Z_{i}.
\]
Thus $P_{i}$ is distributed as a weighted sum of independent $\chi_{1}^{2}$
variables, with weights given by the eigenvalues of $\Gamma^{2}$. As in
Section~2, $\tau_{k,T}=O(1)$ for each fixed $k$. A symbol with two indices,
such as $\delta_{ts}$, denotes the Kronecker delta; with a single index,
$\delta_{i}$ always denotes the diagonal deviation above.

For a matrix $H$, $\Vert H\Vert_{op}$ and $\Vert H\Vert_{F}$ denote the
operator and Frobenius norms, and $\langle H,G\rangle_{F}=Tr(H^{\prime}G)$.
For a random variable $Y$ and $p\geq1$, $\Vert Y\Vert_{p}=(\mathbb{E}%
|Y|^{p})^{1/p}$. Throughout, $C$ and $c$ denote positive constants that do not
depend on $(n,T)$ and may change from line to line.

\paragraph{Uniformity convention and Gaussian moment bounds.}

For a fixed data-generating temporal covariance sequence satisfying
Assumption~2, put $K_{\Gamma}=\sup_{T}\|\Gamma_{T}\|_{op}<\infty$. All bounds
below, including the constants implicit in $O(\cdot)$, are uniform in the
indices and in $n,T$ for this sequence, with constants allowed to depend on
$K_{\Gamma}$ and on the fixed moment order. A bound uniform over a
\emph{class} of temporal covariance sequences requires the additional
common-envelope restriction $\sup_{\{\Gamma_{T}\}\in\mathcal{G},\,T}
\|\Gamma_{T}\|_{op}<\infty$; Assumption~2 alone does not supply one common
constant over every admissible sequence.

\begin{lemma}
[Uniform Gaussian moment bounds]\label{lem:uniform-gaussian-moments} For each
fixed integer $m\geq1$, under Assumptions~1--2 and $H_{0}$,
\[
\mathbb{E}|S_{ii}-1|^{2m}\leq C_{m} T^{-m},\qquad\mathbb{E}|S_{ij}|^{2m}\leq
C_{m} T^{-m}\quad(i\ne j),
\]
where $C_{m}$ depends only on $m$ and $K_{\Gamma}$. In particular,
$\|\delta_{i}\|_{2m}=O(T^{-1/2})$ and $\|S_{ij}\|_{2m}=O(T^{-1/2})$, uniformly
in the indices.
\end{lemma}

\begin{proof}
Write $\delta_{i}=T^{-1}(Z_{i}^{\prime}\Gamma Z_{i}-Tr\Gamma)$. The centered
Gaussian quadratic-form moment inequality gives $\mathbb{E}|Z_{i}^{\prime
}\Gamma Z_{i}-Tr\Gamma|^{2m} \le C_{m}[Tr(\Gamma^{2})]^{m}$. Since
$Tr(\Gamma^{2})\le K_{\Gamma}T$, the first bound follows. For $i\ne j$,
conditional on $Z_{i}$, $S_{ij}$ is centered Gaussian with variance
$T^{-2}Z_{i}^{\prime}\Gamma^{2}Z_{i}$. Hence
\[
\mathbb{E}|S_{ij}|^{2m}=(2m-1)!!\,T^{-2m} \mathbb{E}(Z_{i}^{\prime}\Gamma
^{2}Z_{i})^{m}.
\]
Minkowski's inequality applied to the nonnegative eigenvalue-weighted
chi-square sum gives $\|Z_{i}^{\prime}\Gamma^{2}Z_{i}\|_{m}\le\|\chi_{1}%
^{2}\|_{m} Tr(\Gamma^{2})\le C_{m}T$, proving the second bound.
\end{proof}

\subsection{Isserlis' (Wick's) Theorem}

\noindent The moment calculations used throughout the proofs rely on the
following standard Gaussian identities. Let $Z\sim N(0, I_{T})$ and let $A$,
$B$, and $C$ be deterministic symmetric $T\times T$ matrices.

\begin{lemma}
[Isserlis identities; Isserlis (1918)]\label{lem:isserlis}
\begin{align}
&  \mathbb{E}[Z^{\prime}AZ]=Tr(A),\label{I1}\\
&  \mathbb{E}[(Z^{\prime}AZ)^{2}]=[Tr(A)]^{2}+2\,Tr(A^{2}),\label{I2a}\\
&  \mathbb{E}[(Z^{\prime}AZ)(Z^{\prime}BZ)]=Tr(A)Tr(B)+2Tr(AB),\label{I2b}\\
&  \mathbb{E}[(Z^{\prime}AZ)(Z^{\prime}BZ)(Z^{\prime}CZ)]\nonumber\\
&  \quad=Tr(A)Tr(B)Tr(C)+2\,Tr(A)Tr(BC)+2\,Tr(B)Tr(AC)\nonumber\\
&  \quad\quad+2\,Tr(C)Tr(AB)+8\,Tr(ABC). \label{I3}%
\end{align}

\end{lemma}

These identities follow from Wick's theorem (Isserlis, 1918) by summing over
the pairings of the Gaussian coordinates; \eqref{I1} is immediate, and
\eqref{I2a} is the case $B=A$ of \eqref{I2b}.

\noindent The following corollary collects the special cases used repeatedly
in the proofs.

\begin{corollary}
\label{cor:isserlis} Under $Z_{i}\sim N(0,I_{T})$ independent across $i$:
\begin{align}
\mathbb{E}[P_{i}]=\mathbb{E}[Z_{i}^{\prime}\Gamma^{2}Z_{i}]  &  =T\tau
_{2,T},\label{EPi}\\
\mathbb{E}[P_{i}\cdot Z_{i}^{\prime}\Gamma Z_{i}]  &  =(T\tau_{2,T}%
)(T)+2T\tau_{3,T}=T^{2}\tau_{2,T}+2T\tau_{3,T},\label{EPiGam}\\
\mathbb{E}[P_{i}^{2}]=\mathbb{E}[(Z_{i}^{\prime}\Gamma^{2}Z_{i})^{2}]  &
=(T\tau_{2,T})^{2}+2T\tau_{4,T}=T^{2}\tau_{2,T}^{2}+2T\tau_{4,T}%
,\label{EPi2}\\
\mathbb{E}[P_{i}(Z_{i}^{\prime}\Gamma Z_{i})^{2}]  &  =T^{3}\tau_{2,T}%
+2T^{2}\tau_{2,T}^{2}+4T^{2}\tau_{3,T}+8T\tau_{4,T}. \label{EPiGam2}%
\end{align}
Identities \eqref{EPiGam}, \eqref{EPi2} and \eqref{EPiGam2} follow from
\eqref{I2b}, \eqref{I2a} and \eqref{I3} with $(A,B)=(\Gamma^{2},\Gamma)$,
$A=\Gamma^{2}$ and $(A,B,C)=(\Gamma^{2},\Gamma,\Gamma)$, respectively.
\end{corollary}

\begin{lemma}
[Rate for the spectral moments]\label{lem:tau-rate}Let $\Gamma_{T}$ satisfy
Assumption~2, write $\gamma_{h}=\gamma_{|h|}$ for $h\in\mathbb{Z}$, and put
$A=\sum_{h\in\mathbb{Z}}|\gamma_{h}|$ and $B=\sum_{h\in\mathbb{Z}}%
|h||\gamma_{h}|$. If $B<\infty$, then for each fixed $k\geq2$ and every $T$,
\[
|\tau_{k,T}-\tau_{k}|\leq\frac{(k-1)BA^{k-2}}{T}.
\]
Consequently $K_{T}-K=O(T^{-1})$.
\end{lemma}

\begin{proof}
Since $f(\omega)=\sum_{h}\gamma_{h}e^{i\,h\omega}$ with $A<\infty$,
\[
\tau_{k}=(2\pi)^{-1}\int_{-\pi}^{\pi}f^{k}=\sum_{h_{1}+\cdots+h_{k}=0}%
\prod_{j=1}^{k}\gamma_{h_{j}},
\]
an absolutely convergent sum. With $t_{k+1}=t_{1}$, $Tr(\Gamma_{T}^{k}%
)=\sum_{t_{1},\ldots,t_{k}=1}^{T}\prod_{j=1}^{k}\gamma_{t_{j}-t_{j+1}}$.
Parametrize by $t_{1}$ and the lags $h_{j}=t_{j}-t_{j+1}$, $j<k$, so that
$h_{k}=-(h_{1}+\cdots+h_{k-1})$. For fixed lags, the points $t_{j}%
=t_{1}-(h_{1}+\cdots+h_{j-1})$ all lie in $\{1,\ldots,T\}$ for exactly
$(T-w)_{+}$ values of $t_{1}$, where $w$ is the range of the partial sums and
satisfies $w\leq\sum_{j<k}|h_{j}|$. Hence, using $|\gamma_{h_{k}}|\leq
\gamma_{0}=1$,
\begin{align*}
|T\tau_{k}-Tr(\Gamma_{T}^{k})|  &  \leq\sum_{h_{1},\ldots,h_{k-1}}\prod
_{j<k}|\gamma_{h_{j}}|\min(T,w)\\
&  \leq\sum_{h_{1},\ldots,h_{k-1}}\prod_{j<k}|\gamma_{h_{j}}|\sum_{j<k}%
|h_{j}|=(k-1)BA^{k-2}.
\end{align*}
Dividing by $T$ gives the bound. Since $K_{T}$ and $K$ are the same polynomial
in $(\tau_{2,T},\tau_{3,T},\tau_{4,T})$ and $(\tau_{2},\tau_{3},\tau_{4})$,
with bounded arguments, $K_{T}-K=O(T^{-1})$.
\end{proof}

\noindent For $i\ne j$, the following conditional moments hold.

\begin{lemma}
[Conditional distribution of $S_{ij}$]\label{lem:condSij} Given $Z_{i}$, the
random variable $S_{ij}=\frac{1}{T}Z_{j}^{\prime}\Gamma Z_{i}$ is Gaussian
with
\[
\mathbb{E}[S_{ij}\,|\,Z_{i}]=0,\qquad Var[S_{ij}\,|\,Z_{i}]=\frac
{Z_{i}^{\prime}\Gamma^{2}Z_{i}}{T^{2}}=\frac{P_{i}}{T^{2}}.
\]
In particular, $\mathbb{E}[S_{ij}^{2}]=\tau_{2,T}/T$ by \eqref{EPi}.
\end{lemma}

This holds because, conditional on $Z_{i}$, $S_{ij}=Z_{j}^{\prime}(\Gamma
Z_{i}/T)$ is a linear form in $Z_{j}\sim N(0,I_{T})$, which is independent of
$Z_{i}$.

\subsection{Diagonal concentration: proof of Lemmas~1 and~4}

\label{app:diag}

\noindent Lemmas~1 and~4 are special cases of Hounyo et al.\ (2026a, Thms.~1
and~5). We give the proof for completeness; this subsection allows a general
$U$ satisfying Assumption~4. Since Lemma~1 is the case $U=I_{n}$ of Lemma~4,
it suffices to prove Lemma~4.

\begin{proof}
[Proof of Lemma~4]Let $X=U^{1/2}Z\Gamma^{1/2}$ and $a_{i}=U^{1/2}e_{i}$, so
that $\Vert a_{i}\Vert^{2}=U_{ii}$ and
\[
S_{ii}=\frac{1}{T}a_{i}^{\prime}Z\Gamma Z^{\prime}a_{i}=\operatorname{vec}%
(Z)^{\prime}B_{i}\operatorname{vec}(Z),\qquad B_{i}=\frac{\Gamma}{T}%
\otimes(a_{i}a_{i}^{\prime}).
\]
Since $\operatorname{diag}(S)-\operatorname{diag}(U)$ is diagonal, its
operator norm is $\max_{i}|S_{ii}-U_{ii}|$. Because $Tr(\Gamma)=T$, we have
$\mathbb{E}(S_{ii})=Tr(B_{i})=U_{ii}$. Moreover,
\[
\Vert B_{i}\Vert_{F}^{2}=\frac{Tr(\Gamma^{2})U_{ii}^{2}}{T^{2}}=\frac
{\tau_{2,T}U_{ii}^{2}}{T}\leq\frac{\bar{\tau}_{2}c_{2}^{2}}{T},\qquad\Vert
B_{i}\Vert_{op}=\frac{\Vert\Gamma\Vert_{op}U_{ii}}{T}\leq\frac{K_{\Gamma}%
c_{2}}{T},
\]
where $\bar{\tau}_{2}=\sup_{T}\tau_{2,T}$ and $K_{\Gamma}=\sup_{T}\Vert
\Gamma\Vert_{op}$ are finite by Assumption~2, and $U_{ii}\leq c_{2}$ by
Assumption~4. Since $\operatorname{vec}(Z)$ has i.i.d.\ $N(0,1)$ entries, the
Hanson--Wright inequality (Rudelson \& Vershynin, 2013, Thm.~1.1) gives, for
every $\varepsilon>0$,
\[
P(|S_{ii}-U_{ii}|>\varepsilon)\leq2\exp\left\{  -c\min\left(  \frac
{\varepsilon^{2}}{\Vert B_{i}\Vert_{F}^{2}},\frac{\varepsilon}{\Vert
B_{i}\Vert_{op}}\right)  \right\}  ,
\]
where $c>0$ is a universal constant independent of $i$, $n$, and $T$. Set
$\varepsilon_{n}=M\sqrt{\log n/T}$. Then
\[
\frac{\varepsilon_{n}^{2}}{\Vert B_{i}\Vert_{F}^{2}}\geq\frac{M^{2}\log
n}{\bar{\tau}_{2}c_{2}^{2}},\qquad\frac{\varepsilon_{n}}{\Vert B_{i}\Vert
_{op}}\geq\frac{M\sqrt{T\log n}}{K_{\Gamma}c_{2}}.
\]
Under Assumption~3, $n\rightarrow\infty$, $T\rightarrow\infty$, and $T/\log
n\rightarrow\infty$, so the first term is the smaller one for all large $n$,
and
\[
P\left(  |S_{ii}-U_{ii}|>M\sqrt{\frac{\log n}{T}}\right)  \leq2n^{-cM^{2}%
/(\bar{\tau}_{2}c_{2}^{2})}%
\]
uniformly in $i$. A union bound over $i=1,\ldots,n$ gives
\[
P\left(  \max_{i}|S_{ii}-U_{ii}|>M\sqrt{\frac{\log n}{T}}\right)
\leq2n^{1-cM^{2}/(\bar{\tau}_{2}c_{2}^{2})}.
\]
Choosing $M>(2\bar{\tau}_{2}c_{2}^{2}/c)^{1/2}$ makes the exponent smaller
than $-1$, so the right-hand side is summable in $n$. By the Borel--Cantelli
lemma, $\max_{i}|S_{ii}-U_{ii}|\leq M\sqrt{\log n/T}$ for all sufficiently
large $n$, almost surely. This proves Lemma~4, and Lemma~1 follows with
$U=I_{n}$ and $c_{2}=1$.
\end{proof}

For every fixed $\varepsilon>0$, the same exponential bound yields
$P(|S_{ii}-U_{ii}|>\varepsilon)\leq C_{\varepsilon}e^{-c_{\varepsilon}T}$
uniformly in $i$, where $C_{\varepsilon},c_{\varepsilon}>0$ may depend on
$\varepsilon$ but not on $i$, $n$, or $T$. This per-unit tail bound is the
form used in the proofs of Lemmas~2, 3, and~5.

\subsection{Diagonal normalization}

\noindent Write $R=D_{S}^{-1/2}SD_{S}^{-1/2}$, where $D_{S}%
=\operatorname{diag}(S_{11},\ldots,S_{nn})$, and let $\mathcal{E}_{n}%
=\{\Vert\delta\Vert_{\infty}\leq1/2\}$.

\begin{lemma}
[Expansion of the diagonal normalization]\label{lem:diag-expansion}Under
Assumptions 1--3 and $H_{0}$, $P(\mathcal{E}_{n})\rightarrow1$, and on
$\mathcal{E}_{n}$ the following hold uniformly in $i\neq j$:

\noindent(a) $S_{jj}^{-1}=1-\delta_{j}+O(\delta_{j}^{2})$;

\noindent(b) $(S_{ii}S_{jj})^{-1}=1-(\delta_{i}+\delta_{j})+O(\Vert\delta
\Vert_{\infty}^{2})$;

\noindent(c) $r_{ij}^{2}=S_{ij}^{2}-S_{ij}^{2}(\delta_{i}+\delta_{j}%
)+O(S_{ij}^{2}\Vert\delta\Vert_{\infty}^{2})$.
\end{lemma}

\begin{proof}
By Lemma~1, $\Vert\delta\Vert_{\infty}=\Vert\operatorname{diag}(S)-I_{n}%
\Vert_{op}=O(\sqrt{\log n/T})$ almost surely, so $\Vert\delta\Vert_{\infty
}\rightarrow0$ almost surely and $P(\mathcal{E}_{n})\rightarrow1$. On
$\mathcal{E}_{n}$, every $S_{jj}\in\lbrack1/2,3/2]$, and Taylor's formula
gives $S_{jj}^{-1}=(1+\delta_{j})^{-1}=1-\delta_{j}+\varepsilon_{j}$ with
$|\varepsilon_{j}|\leq C\delta_{j}^{2}$ uniformly in $j$. This is (a); (b)
follows by multiplying the expansions for $i$ and $j$, and (c) follows from
(b) since $r_{ij}^{2}=S_{ij}^{2}/(S_{ii}S_{jj})$.
\end{proof}

The pairwise bounds in Lemma~\ref{lem:diag-expansion} do not by themselves
make the aggregate remainder negligible. The displayed big-$O$ statements are
pointwise deterministic inequalities on $\mathcal{E}_{n}$, with constants
independent of $i,j,n,T$ once $K_{\Gamma}$ is fixed. Summed over all pairs,
$S_{ij}^{2}\Vert\delta\Vert_{\infty}^{2}$ is only $O(\log n)$, and the
first-order terms $S_{ij}^{2}(\delta_{i}+\delta_{j})$ contribute at the CLT
scale to both the mean and the variance. The proofs of Lemmas~2 and~3
therefore retain these terms explicitly and control the remainders in
expectation after summation over pairs.

\section{Proof of Lemma~5}

\label{proof:mean_alt}

We establish a uniform expansion for a single squared sample correlation and
then sum it. Fix $i\neq j$ and write $\rho=\rho_{ij}^{\ast}$. Put
$v_{a}=Z^{\prime}U^{1/2}e_{a}/U_{aa}^{1/2}$ for $a=i,j$, so that $v_{a}\sim
N(0,I_{T})$ and $Cov(v_{i},v_{j})=\rho I_{T}$. Then $r_{ij}=s/(q_{i}%
q_{j})^{1/2}$ with
\[
s=T^{-1}v_{i}^{\prime}\Gamma v_{j}=\frac{S_{ij}}{\sqrt{U_{ii}U_{jj}}},\qquad
q_{a}=T^{-1}v_{a}^{\prime}\Gamma v_{a}=\frac{S_{aa}}{U_{aa}}.
\]

\begin{lemma}
[Uniform pairwise mean expansion]\label{lem:pair-mean} Under Assumptions~1--4,
uniformly in $i\ne j$ and all sufficiently large $T$,
\begin{equation}
\label{eq:pair-mean-uniform}\mathbb{E}_{A}(r_{ij}^{2})=\rho^{2}+\frac
{\tau_{2,T}}{T}(1-3\rho^{2}+2\rho^{4}) +\frac{4(\tau_{2,T}^{2}-\tau_{3,T}%
)}{T^{2}}+\varepsilon_{ij}, \qquad|\varepsilon_{ij}|\le C\left(  \frac
{\rho^{2}}{T^{2}}+\frac1{T^{3}}\right)  .
\end{equation}
The constant $C$ depends only on $K_{\Gamma}$, not on $i,j,n,T$ or $\rho$.
\end{lemma}

\begin{proof}
The argument does not divide by $1-\rho^{2}$, so the bounds are uniform as
$|\rho|\to1$. In an eigenbasis of $\Gamma$,
\[
(s,q_{i},q_{j})=T^{-1}\sum_{t=1}^{T}\nu_{t,T}(X_{t}Y_{t},X_{t}^{2},Y_{t}%
^{2}),\qquad Y_{t}=\rho X_{t}+\sqrt{1-\rho^{2}}\,Z_{t},
\]
where the $X_{t},Z_{t}$ are independent standard normals. Put $\eta=(\eta
_{1},\eta_{2},\eta_{3})=(s-\rho,q_{i}-1,q_{j}-1)$ and $h(x,y,z)=x^{2}/(yz)$,
so that $r_{ij}^{2}=h((\rho,1,1)+\eta)$ and $0\leq r_{ij}^{2}\leq1$. Let
$\zeta=(XY-\rho,X^{2}-1,Y^{2}-1)$ for a single bivariate normal pair with unit
variances and correlation $\rho$, and let $C_{ab}=\mathbb{E}(\zeta_{a}%
\zeta_{b})$ and $H_{abc}=\mathbb{E}(\zeta_{a}\zeta_{b}\zeta_{c})$;
explicitly,
\[
C=\left(
\begin{array}
[c]{ccc}%
1+\rho^{2} & 2\rho & 2\rho\\
2\rho & 2 & 2\rho^{2}\\
2\rho & 2\rho^{2} & 2
\end{array}
\right)  .
\]
Since $\eta=T^{-1}\sum_{t}\nu_{t,T}\zeta_{t}$ with i.i.d.\ copies $\zeta_{t}$
of $\zeta$, the $k$th joint cumulant of coordinates of $\eta$ equals
$T^{-k}Tr(\Gamma^{k})$ times the corresponding cumulant of $\zeta$, hence is
$O(T^{1-k})$ uniformly in $|\rho|\leq1$, and is a polynomial in $\rho$. In
particular,
\begin{align}
&  \mathbb{E}(\eta_{a}\eta_{b})=\frac{\tau_{2,T}}{T}C_{ab},\qquad
\mathbb{E}(\eta_{a}\eta_{b}\eta_{c})=\frac{\tau_{3,T}}{T^{2}}H_{abc}%
,\nonumber\\
&  \mathbb{E}(\eta_{a}\eta_{b}\eta_{c}\eta_{d})=\frac{\tau_{2,T}^{2}}{T^{2}%
}(C_{ab}C_{cd}+C_{ac}C_{bd}+C_{ad}C_{bc})+O(T^{-3}), \label{eq:eta-moments}%
\end{align}
where the first two identities are exact. By the moment--cumulant formula, a
product of $d$ coordinates of $\eta$ has expectation $O(T^{-\lceil d/2\rceil
})$; in particular $\mathbb{E}\Vert\eta\Vert^{2m}\leq C_{m}T^{-m}$ for each
fixed $m$, and degrees five and six contribute $O(T^{-3})$. Moreover, since
each $q_{a}$ has the null distribution of $S_{aa}$, $P\{\min(q_{i}%
,q_{j})<1/2\}\leq Ce^{-cT}$ by the per-unit tail bound following the proof of Lemma~4.

On $\{\min(q_{i},q_{j})\geq1/2\}$, apply Taylor's theorem to $h$ at
$(\rho,1,1)$ through total degree five. The segment joining $(\rho,1,1)$ and
$(s,q_{i},q_{j})$ lies in $\{(x,y,z):y,z\geq1/2,\ x^{2}\leq yz\}$, because
both endpoints are entries of $2\times2$ positive semidefinite matrices with
diagonal entries at least $1/2$ and that set is convex. Every sixth-order
derivative of $h$ is a finite sum of terms $x^{k}y^{-a}z^{-b}$ with $k\leq2$
and $a,b\geq1$; since $|x|^{k}\leq(yz)^{k/2}$ and $y,z\geq1/2$ there, it is
bounded on this set. The expected sixth-order remainder is therefore
$O(\mathbb{E}\Vert\eta\Vert^{6})=O(T^{-3})$. The complementary event
contributes $O(e^{-cT})$, using $0\leq r_{ij}^{2}\leq1$ for the correlation
and Cauchy--Schwarz for the Taylor polynomial. Writing $P_{d}$ for the
homogeneous degree-$d$ part of $h((\rho,1,1)+\eta)$,
\[
P_{d}(\eta)=\sum_{k=0}^{\min(2,d)}\binom{2}{k}\rho^{2-k}\eta_{1}^{k}%
(-1)^{d-k}\sum_{l=0}^{d-k}\eta_{2}^{l}\eta_{3}^{d-k-l},
\]
we have $\mathbb{E}P_{1}(\eta)=0$, and \eqref{eq:eta-moments} gives
\[
\mathbb{E}_{A}(r_{ij}^{2})=\rho^{2}+T^{-1}\mathcal{A}_{T}(\rho)+T^{-2}%
\mathcal{B}_{T}(\rho)+O(T^{-3}),
\]
uniformly in $|\rho|\leq1$, where $\mathcal{A}_{T}$ collects the
$C$-contraction of $P_{2}$, and $\mathcal{B}_{T}$ collects the $H$-contraction
of $P_{3}$ and the pairwise $C$-contractions of $P_{4}$. For example,
$P_{2}(\eta)=\eta_{1}^{2}-2\rho\eta_{1}(\eta_{2}+\eta_{3})+\rho^{2}(\eta
_{2}^{2}+\eta_{2}\eta_{3}+\eta_{3}^{2})$, whose $C$-contraction is
$(1+\rho^{2})-8\rho^{2}+\rho^{2}(4+2\rho^{2})=1-3\rho^{2}+2\rho^{4}$. Carrying
out the contractions,
\[
\mathcal{A}_{T}(\rho)=\tau_{2,T}(1-3\rho^{2}+2\rho^{4}),\qquad\mathcal{B}%
_{T}(\rho)=2(\rho^{2}-1)\left[  (4\rho^{4}+\rho^{2}-2)\tau_{2,T}^{2}%
+(2-4\rho^{2})\tau_{3,T}\right]  .
\]
At $\rho=0$, $\mathcal{B}_{T}(0)=4(\tau_{2,T}^{2}-\tau_{3,T})$, in agreement
with the null expansion in Appendix~A of the main paper. Moreover,
\[
\mathcal{B}_{T}(\rho)-\mathcal{B}_{T}(0)=2\rho^{2}\left[  (4\rho^{4}-3\rho
^{2}-3)\tau_{2,T}^{2}+(6-4\rho^{2})\tau_{3,T}\right]  ,
\]
so $|\mathcal{B}_{T}(\rho)-\mathcal{B}_{T}(0)|\leq C\rho^{2}$ for $|\rho
|\leq1$, since $\tau_{2,T}$ and $\tau_{3,T}$ are bounded. This proves
\eqref{eq:pair-mean-uniform}. The case $|\rho|=1$ is included: there
$r_{ij}^{2}=1$ and $\mathcal{A}_{T}(\pm1)=\mathcal{B}_{T}(\pm1)=0$.
\end{proof}

We now sum \eqref{eq:pair-mean-uniform} over the $n(n-1)$ ordered pairs. Since
$\sum_{i\ne j}(\rho_{ij}^{\ast})^{2}=n(\phi_{2}-1)=O(n)$ by Assumption~4 and
$T=n\lambda$,
\[
\sum_{i\ne j}|\varepsilon_{ij}| \le C\left\{  \frac{n(\phi_{2}-1)}{T^{2}%
}+\frac{n(n-1)}{T^{3}}\right\}  =O(n^{-1}).
\]
Hence, using $Tr(R^{2})=n+\sum_{i\ne j}r_{ij}^{2}$, $n+n(\phi_{2}%
-1)=Tr(R_{\ast}^{2})$, $\sum_{i\neq j}(1-3\rho_{ij}^{\ast2}+2\rho_{ij}^{\ast
4})=n(n-1)-3n(\phi_{2}-1)+2\sum_{i\neq j}(\rho_{ij}^{\ast})^{4}$, and
$n/T=1/\lambda$,
\begin{align*}
\mathbb{E}_{A}[Tr(R^{2})]  &  =Tr(R_{\ast}^{2})+\frac{(n-1)\tau_{2,T}}%
{\lambda}+\frac{4(n-1)(\tau_{2,T}^{2}-\tau_{3,T})}{n\lambda^{2}}\\
&  \quad-\frac{3\tau_{2,T}(\phi_{2}-1)}{\lambda}+\frac{2\tau_{2,T}}{T}%
\sum_{i\ne j}(\rho_{ij}^{\ast})^{4}+O(n^{-1}),
\end{align*}
which is (3). All constants are uniform under Assumptions~1--4.

%=============================================================

\section{Proof of Lemma 6}

\label{proof:var_alt}

%=============================================================

\noindent We establish the two assertions of Lemma~6 through separate
auxiliary results: a uniform Gaussian Poincar\'e bound and a sharper
covariance comparison under local alternatives. We also prove here the
coupling result used for Theorem~2 (Lemma~\ref{lem:coupling}).

Throughout, $R_{\ast}=I_{n}+E$ with $E$ symmetric and zero-diagonal, so that
$\Vert E\Vert_{F}^{2}=\sum_{i\neq j}(\rho_{ij}^{\ast})^{2}=n(\phi_{2}-1)$ and,
by Assumption~4, $\Vert R_{\ast}\Vert_{op}\leq c_{2}/c_{1}$. Because $R$ is
invariant to rescaling the rows of $X$, its law depends on $U$ only through
$R_{\ast}$; without loss of generality we take
\[
X=WZ\Gamma^{1/2},\qquad W=R_{\ast}^{1/2}=I_{n}+\Theta,
\]
with $W$ the symmetric square root. Since $W^{2}=I_{n}+E$ and $Tr(E)=0$,
\begin{equation}
2Tr(\Theta)+Tr(\Theta^{2})=0\quad\text{exactly},\qquad\Vert\Theta\Vert_{F}%
\leq\Vert E\Vert_{F},\qquad\Vert\Theta\Vert_{op}\leq C, \label{eq:F-identity}%
\end{equation}
where the Frobenius bound follows from $|\sqrt{\mu}-1|\leq|\mu-1|$ applied to
the eigenvalues $\mu$ of $R_{\ast}$.

\medskip\noindent Operator-norm moments. Put $c_{+}=\sup n/T<\infty$ and
recall $K_{\Gamma}=\sup_{T}\|\Gamma_{T}\|_{op}$; for $p\ge2$ let
\[
z_{p}=1+\sqrt{c_{+}}+\sqrt{4p}+2,\qquad m_{p}=1+K_{\Gamma}z_{2p}^{2};
\]
these constants depend only on $c_{+}$ and $K_{\Gamma}$, not on $n,T$. The
Gaussian operator-norm tail inequality $\Pr\{\|Z\|_{op}>\sqrt T+\sqrt n+t\}\le
e^{-t^{2}/2}$ for $t\ge0$ (Davidson \& Szarek, 2001, Thm.~II.13), followed by
integration of the tail, gives, for fixed $p\ge2$ and all sufficiently large
$T$,
\begin{equation}
\|\|Z\Gamma\|_{op}\|_{p}\le K_{\Gamma}z_{p}\sqrt T, \qquad\|\|N\|_{op}%
\|_{p}\le m_{p}, \qquad N=T^{-1}Z\Gamma Z^{\prime}-I_{n}.
\label{eq:explicit-N-moments}%
\end{equation}
To verify these bounds, let $\xi=\|Z\|_{op}-\sqrt T-\sqrt n$. Tail integration
gives $\|\xi_{+}\|_{p}\le\sqrt{4p}+2$, so $\|\|Z\|_{op}\|_{p}\le z_{p}\sqrt
T$; the first bound then follows from $\|Z\Gamma\|_{op}\le K_{\Gamma
}\|Z\|_{op}$, and the second from $\|N\|_{op}\le1+K_{\Gamma}\|Z\|_{op}^{2}/T$
and $\|\|Z\|_{op}^{2}\|_{p}=\|\|Z\|_{op}\|_{2p}^{2}$. No invertibility of
$\Gamma$ is needed.

\begin{lemma}
[Variance bound]\label{lem:poincare}Under Assumptions 1--4, $Var(Tr(R^{2}%
))\leq C$ for a constant $C$ not depending on $(n,T)$.
\end{lemma}

\begin{proof}
Let $x_{i}$ denote the rows of $X$ and $y_{i}=x_{i}/\Vert x_{i}\Vert$, so that
$f=Tr(R^{2})=\sum_{i,j}(y_{i}^{\prime}y_{j})^{2}$. On the full-measure set
where all $x_{i}\neq0$, $f$ is smooth and
\[
\frac{\partial f}{\partial x_{i}}=\frac{4}{\Vert x_{i}\Vert}(I_{T}-y_{i}%
y_{i}^{\prime})\sum_{j}r_{ij}y_{j}.
\]
Since $\Vert\sum_{j}r_{ij}y_{j}\Vert^{2}=(R^{3})_{ii}\leq\Vert R\Vert
_{op}(R^{2})_{ii}$ and $Tr(R^{2})\leq n\Vert R\Vert_{op}$,
\[
\sum_{i=1}^{n}\left\Vert \frac{\partial f}{\partial x_{i}}\right\Vert ^{2}
\leq\frac{16\Vert R\Vert_{op}Tr(R^{2})}{\min_{1\leq i\leq n}\Vert x_{i}%
\Vert^{2}}\leq\frac{16n\Vert R\Vert_{op}^{2}}{T\min_{1\leq i\leq n}S_{ii}}.
\]
The map $Z\mapsto X$ is linear (in Frobenius norm) with operator norm at most
$\Vert W\Vert_{op}\Vert\Gamma^{1/2}\Vert_{op}\leq C$, and $\Vert R\Vert
_{op}\leq\Vert S\Vert_{op}/\min_{i}S_{ii}\leq C\Vert ZZ^{\prime}/T\Vert
_{op}/\min_{i}S_{ii}$. Hence
\[
\Vert\nabla_{Z}f\Vert^{2}\leq C(n/T)\Vert ZZ^{\prime}/T\Vert_{op}^{2}(\min
_{i}S_{ii})^{-3}.
\]
Under Assumption 3, $\mathbb{E}\Vert ZZ^{\prime}/T\Vert_{op}^{4}=O(1)$ by
\eqref{eq:explicit-N-moments} with $\Gamma=I_{T}$. Each $S_{ii}$ is
distributed as $v^{\prime}\Gamma v/T$ with $v\sim N(0,I_{T})$, which has
bounded inverse moments of any fixed order by (B.2) and satisfies
$P(S_{ii}<1/2)\leq Ce^{-cT}$; therefore $\mathbb{E}(\min_{i}S_{ii})^{-6}%
\leq2^{6}+\sum_{i}(\mathbb{E}S_{ii}^{-12})^{1/2}P(S_{ii}<1/2)^{1/2}=O(1)$. By
H\"{o}lder's inequality, $\mathbb{E}\Vert\nabla_{Z}f\Vert^{2}=O(n/T)=O(1)$,
and the Gaussian Poincar\'{e} inequality gives $Var(f)\leq\mathbb{E}%
\Vert\nabla_{Z}f\Vert^{2}\leq C$. (The singular set where some $x_{i}=0$ has
codimension $\operatorname{rank}(\Gamma)\geq3$ for large $T$, so $f$, which is
bounded by $n^{2}$, belongs to the Gaussian Sobolev space and the inequality applies.)
\end{proof}

\begin{lemma}
[Coupling under local alternatives]\label{lem:coupling}Under Assumptions 1--4,
suppose $n(\phi_{2}-1)=O(1)$. Let $R_{0}$ be the sample correlation matrix of
$X_{0}=Z\Gamma^{1/2}$, built from the same $Z$, so that $R_{0}$ has the null
distribution. Then
\[
D_{n}=Tr(R^{2})-Tr(R_{0}^{2})-n(\phi_{2}-1)\quad\text{satisfies}%
\quad\mathbb{E}(D_{n}^{2})=O(n^{-1}).
\]

\end{lemma}

\begin{proof}
Write $S_{0}=X_{0}X_{0}^{\prime}/T=I_{n}+N$. Then $S=WS_{0}W=R_{\ast}+G$ with
\[
G=WNW=N+H,\qquad H=\Theta N+N\Theta+\Theta N\Theta.
\]
The entries of $N$ form a centered Gaussian chaos of order two with the exact
covariance
\begin{equation}
Cov(N_{ab},N_{cd})=(\delta_{ac}\delta_{bd}+\delta_{ad}\delta_{bc})\frac
{\tau_{2,T}}{T}, \label{eq:N-cov}%
\end{equation}
so that, for a deterministic symmetric matrix $A$, $Var(Tr(AN))=2\tau
_{2,T}\Vert A\Vert_{F}^{2}/T$. Joint cumulants of order $m$ of entries of $N$
admit trace-contraction formulas with scale $O(T^{-(m-1)})$ for fixed $m$
under the uniform spectral-norm bound. Moreover, by Gaussian
hypercontractivity $\Vert Q\Vert_{p}\leq C_{p}\Vert Q\Vert_{2}$ for any
polynomial $Q$ of bounded degree in the entries of $N$.

Put $a_{i}=S_{ii}=1+G_{ii}$ and $b_{i}=(S_{0})_{ii}=1+N_{ii}$. On
$\mathcal{E}=\{a_{i},b_{i}\in\lbrack1/2,2]\ \forall i\}$, whose complement has
probability at most $Cne^{-cT}$, the functions $(a_{i},a_{j})\mapsto
1/(a_{i}a_{j})$ have bounded derivatives of all orders; the probability bound
holds because each $a_{i}$ and $b_{i}$ is a unit-mean Gaussian quadratic form
with eigenvalue weights bounded by $K_{\Gamma}$. Since $n\le Tr(R^{2}%
),Tr(R_{0}^{2})\le n^{2}$ and $\|E\|_{F}^{2}\le n(n-1)$, $|D_{n}|\leq3n^{2}$
always, so $\mathcal{E}^{c}$ contributes $O(n^{5}e^{-cT})$ to $\mathbb{E}%
(D_{n}^{2})$, and we work on $\mathcal{E}$; below, $\|T_{j}\|_{2}$ denotes
$\|T_{j}\mathbf{1}_{\mathcal{E}}\|_{2}$. Expanding $S_{ij}^{2}=(E_{ij}%
+N_{ij}+H_{ij})^{2}$ and using $\sum_{i\neq j}E_{ij}^{2}=n(\phi_{2}-1)$,
\[
D_{n}=T_{1}+T_{2}+T_{3}+T_{4},
\]
where
\begin{align*}
T_{1}  &  =\sum_{i\neq j}E_{ij}^{2}\left(  \frac{1}{a_{i}a_{j}}-1\right)  , &
T_{2}  &  =2\sum_{i\neq j}\frac{E_{ij}G_{ij}}{a_{i}a_{j}},\\
T_{3}  &  =\sum_{i\neq j}\frac{2N_{ij}H_{ij}+H_{ij}^{2}}{a_{i}a_{j}}, & T_{4}
&  =\sum_{i\neq j}N_{ij}^{2}\left(  \frac{1}{a_{i}a_{j}}-\frac{1}{b_{i}b_{j}%
}\right)  .
\end{align*}
Let $w_{i}=We_{i}$, so that $\|w_{i}\|^{2}=(R_{\ast})_{ii}=1$, and let
$f_{i}=\|\Theta_{i\cdot}\|_{2}$, where $\Theta_{i\cdot}$ is the $i$th row of
$\Theta$; thus $\sum_{i}f_{i}\leq\sqrt{n}\|\Theta\|_{F}=O(\sqrt{n})$. Since
$G_{ab}=w_{a}^{\prime}Nw_{b}$ and $H_{ii}=Tr\{(w_{i}w_{i}^{\prime}-e_{i}%
e_{i}^{\prime})N\}$ with $\|w_{i}w_{i}^{\prime}-e_{i}e_{i}^{\prime}\|_{F}%
\leq\sqrt{2}\|w_{i}-e_{i}\|=\sqrt{2}f_{i}$ (both vectors have unit length),
\eqref{eq:N-cov} and hypercontractivity give, for each fixed $p\geq2$ and
uniformly in the indices,
\[
\Vert G_{ab}\Vert_{p}+\Vert N_{ab}\Vert_{p}\leq C_{p}T^{-1/2},\qquad\Vert
H_{ii}\Vert_{p}\leq C_{p}T^{-1/2}f_{i}.
\]

\noindent Term $T_{1}$. On $\mathcal{E}$, $|1/(a_{i}a_{j})-1|\le
C(|G_{ii}|+|G_{jj}|)$, so Minkowski's inequality gives
\[
\|T_{1}\mathbf{1}_{\mathcal{E}}\|_{2} \le C\sum_{i\ne j}E_{ij}^{2}%
(\|G_{ii}\|_{2}+\|G_{jj}\|_{2}) \le CT^{-1/2}\|E\|_{F}^{2}=O(T^{-1/2}).
\]
Thus $\|T_{1}\|_{2}=O(T^{-1/2})$.

\noindent Term $T_{2}$. Taylor expansion on $\mathcal{E}$ yields
\[
T_{2}=L_{2}+Q_{2}+R_{2},\quad L_{2}=2Tr(EG),\quad Q_{2}=-4\sum_{i\neq j}%
E_{ij}G_{ij}G_{ii},
\]
where $|R_{2}|\leq C\sum_{i\neq j}|E_{ij}G_{ij}|(G_{ii}^{2}+G_{jj}^{2})$. The
linear term is centered and, by the exact covariance identity,
\[
\mathbb{E}L_{2}^{2}=\frac{8\tau_{2,T}}{T}\Vert WEW\Vert_{F}^{2}=O(T^{-1}).
\]
For the quadratic term set $D=\operatorname{diag}(G)$ and $J_{2}=Tr(DEG)$, so
$Q_{2}=-4J_{2}$ (because $E_{ii}=0$). Wick's second-moment identity gives
\[
\mathbb{E}J_{2}=\frac{2\tau_{2,T}}{T}\sum_{i\neq j}E_{ij}(R_{\ast})_{ij}%
=\frac{2\tau_{2,T}}{T}\Vert E\Vert_{F}^{2}=O(T^{-1}).
\]
Differentiation with respect to the symmetric matrix $G$ gives
\[
dJ_{2}=Tr(B_{G}\,dG),\qquad B_{G}=\operatorname{diag}(\operatorname{diag}%
(EG))+\frac{1}{2}(ED+DE).
\]
Here $G=WNW$, so $B_{N}=WB_{G}W$ and $\nabla_{Z}J_{2}=2T^{-1}B_{N}Z\Gamma$.
For every $i$,
\[
\mathbb{E}|(EG)_{ii}|^{2}=\frac{\tau_{2,T}}{T}\left[  (ER_{\ast}%
E)_{ii}(R_{\ast})_{ii}+(ER_{\ast})_{ii}^{2}\right]  \leq\frac{C}{T}\Vert
E_{i\cdot}\Vert_{2}^{2}.
\]
Also $\mathbb{E}G_{ii}^{2}=2\tau_{2,T}/T$, whence $\mathbb{E}\Vert B_{G}%
\Vert_{F}^{2}\leq CT^{-1}\Vert E\Vert_{F}^{2}=O(T^{-1})$. Since $B_{G}$ is
linear in the Gaussian quadratic chaos $N$, hypercontractivity yields
$\mathbb{E}\Vert B_{N}\Vert_{F}^{4}=O(T^{-2})$; for example, apply the
fourth-moment inequality entrywise and then Minkowski's inequality to the
Frobenius sum. The standard Gaussian operator-norm moment bound $\mathbb{E}%
\Vert Z\Gamma\Vert_{op}^{4}=O(T^{2})$ now gives, by Gaussian Poincar\'{e},
\[
Var(J_{2})\leq\frac{4}{T^{2}}\mathbb{E}\Vert B_{N}Z\Gamma\Vert_{F}^{2}%
\leq\frac{4}{T^{2}}(\mathbb{E}\Vert B_{N}\Vert_{F}^{4})^{1/2}(\mathbb{E}\Vert
Z\Gamma\Vert_{op}^{4})^{1/2}=O(T^{-2}).
\]
Thus $\Vert Q_{2}\Vert_{2}=O(T^{-1})$. Finally, by H\"{o}lder and the uniform
fixed-order quadratic-chaos moments,
\[
\Vert R_{2}\mathbf{1}_{\mathcal{E}}\Vert_{2}\leq CT^{-3/2}\sum_{i\neq
j}|E_{ij}|\leq CnT^{-3/2}\Vert E\Vert_{F}=O(T^{-1/2}).
\]
Therefore $\Vert T_{2}\Vert_{2}=O(T^{-1/2})$.

\noindent Term $T_{3}$. Set $A=R_{\ast}=I_{n}+E=W^{2}$ and write $T_{3}%
=Q_{3}+C_{3}$, where
\[
Q_{3}=\sum_{i\ne j}\{2N_{ij}H_{ij}+H_{ij}^{2}\} =\sum_{i\ne j}\{(WNW)_{ij}%
^{2}-N_{ij}^{2}\}
\]
is the unweighted sum and $C_{3}=\sum_{i\ne j}(2N_{ij}H_{ij}+H_{ij}%
^{2})\{(a_{i}a_{j})^{-1}-1\}$. We first show that
\[
\mathbb{E} Q_{3}=\frac{\tau_{2,T}}{T}\|E\|_{F}^{2}=O(T^{-1}),\qquad
Var(Q_{3})=O(T^{-1}).
\]
To verify the mean, use $Cov(G_{ij},G_{ij})= (\tau_{2,T}/T)(A_{ii}%
A_{jj}+A_{ij}^{2})$, $A_{ii}=1$ and $A_{ij}=E_{ij}$ for $i\ne j$.

For the variance, first consider the full trace difference
\[
V=Tr(G^{2})-Tr(N^{2}) =Tr(ANAN)-Tr(N^{2}).
\]
Its gradient is
\[
\nabla_{Z} V=\frac4T(ANA-N)Z\Gamma, \qquad ANA-N=EN+NE+ENE.
\]
Gaussian Poincar\'e and Cauchy--Schwarz give
\[
Var(V)\leq\mathbb{E}\|\nabla_{Z}V\|_{F}^{2} \leq\frac{16}{T^{2}}
(\mathbb{E}\|ANA-N\|_{F}^{4})^{1/2} (\mathbb{E}\|Z\Gamma\|_{op}^{4})^{1/2}.
\]
Put $K_{E}=\sup_{n}\|E\|_{F}$, which is finite in the local regime, and let
$b=1+K_{E}$ and $a=K_{E}(2+K_{E})$. Since $\|W\|_{op}\le b$ and $\|\Theta
\|_{F}\le K_{E}$, the deterministic Schatten-norm inequalities give,
pathwise,
\[
\|ANA-N\|_{F}\le a\|N\|_{op},\quad\|H\|_{F}\le a\|N\|_{op},\quad\|G\|_{op}\le
b^{2}\|N\|_{op}.
\]
Consequently, by \eqref{eq:explicit-N-moments},
\begin{equation}
\mathbb{E}\|ANA-N\|_{F}^{4}\le a^{4}m_{4}^{4}, \qquad\mathbb{E}\|Z\Gamma
\|_{op}^{4}\le K_{\Gamma}^{4}z_{4}^{4}T^{2}, \label{eq:explicit-V-fourth}%
\end{equation}
and $Var(V)=O(T^{-1})$; the explicit constant is recorded in
\eqref{eq:explicit-VJ-variance} below. The difference between $V$ and $Q_{3}$
is the diagonal sum $\Delta_{3}=\sum_{i}(G_{ii}^{2}-N_{ii}^{2})=\sum
_{i}(2N_{ii}H_{ii}+H_{ii}^{2})$. By the moment bounds stated at the outset,
\[
\|\Delta_{3}\|_{2}\leq\frac C T\sum_{i} (f_{i}+f_{i}^{2}) \leq\frac C T(\sqrt
n\|\Theta\|_{F}+\|\Theta\|_{F}^{2}) =O(T^{-1/2}).
\]
The estimates for $Q_{3}=V-\Delta_{3}$ follow.

With $D=\operatorname{diag}(G)$ as before, put $P=G\circ G-N\circ N$ (Hadamard
product), so that $P_{ij}=H_{ij}(G_{ij}+N_{ij})$, and let $\bar{P}_{i}%
=\sum_{j\ne i}P_{ij}$. On $\mathcal{E}$, the first-order Taylor formula for
$(a_{i}a_{j})^{-1}$ with quadratic remainder and the symmetry of $P$ give
\[
C_{3}=-2\sum_{i}G_{ii}\bar{P}_{i}+\mathcal{R}_{3},\qquad|\mathcal{R}_{3}|\le
C\sum_{i}G_{ii}^{2}\sum_{j\ne i}|P_{ij}|.
\]
The leading term has the exact trace representation
\[
C_{3,1}=-2Tr\{D(G^{2}-N^{2})\} +2\sum_{i}G_{ii}(G_{ii}^{2}-N_{ii}^{2}).
\]

For the trace polynomial $J=Tr\{D(G^{2}-N^{2})\}$, differentiation with
respect to a symmetric perturbation $dN$ yields
\[
dJ=Tr\left(  B_{J}\,dN\right)  ,\quad B_{J}=W\operatorname{diag}\left(
\operatorname{diag}(G^{2}-N^{2})\right)  W+W(GD+DG)W-(ND+DN).
\]
Consequently, since $N=T^{-1}Z\Gamma Z^{\prime}-I_{n}$,
\[
\nabla_{Z}J=\frac{2}{T}B_{J}Z\Gamma,\qquad Var(J)\leq\frac{4}{T^{2}%
}\{\mathbb{E}\Vert B_{J}\Vert_{F}^{4}\}^{1/2}\{\mathbb{E}\Vert Z\Gamma
\Vert_{op}^{4}\}^{1/2}.
\]
To bound $B_{J}$, note that each $G_{ii}=T^{-1}w_{i}^{\prime}Z\Gamma
Z^{\prime}w_{i}-1$ is a centered quadratic form with $\Vert w_{i}\Vert=1$ and
variance $2\tau_{2,T}/T\leq2K_{\Gamma}/T$; degree-two Gaussian
hypercontractivity gives $\Vert G_{ii}\Vert_{8}\leq7\sqrt{2K_{\Gamma}/T}$, and
Minkowski's inequality in $L^{4}$ applied to $\Vert D\Vert_{F}^{2}$ yields
\[
\Vert\Vert D\Vert_{F}\Vert_{8}\leq7\sqrt{2K_{\Gamma}n/T}\leq d_{8},\qquad
d_{8}=7\sqrt{2K_{\Gamma}c_{+}}.
\]
The three terms of $B_{J}$ satisfy the deterministic inequality
\[
\Vert B_{J}\Vert_{F}\leq b^{2}(\Vert G\Vert_{op}+\Vert N\Vert_{op})\Vert
H\Vert_{F}+2(b^{2}\Vert G\Vert_{op}+\Vert N\Vert_{op})\Vert D\Vert_{F},
\]
since $\Vert\operatorname{diag}(\operatorname{diag}M)\Vert_{F}\leq\Vert
M\Vert_{F}$ and $G^{2}-N^{2}=GH+HN$. H\"{o}lder's inequality, with the eighth
moments above, gives the dimension-free bound
\begin{equation}
\left(  \mathbb{E}\Vert B_{J}\Vert_{F}^{4}\right)  ^{1/4}\leq b^{2}%
(b^{2}+1)a\,m_{8}^{2}+2(b^{4}+1)m_{8}d_{8}=:C_{J}<\infty.
\label{eq:explicit-J-fourth}%
\end{equation}
The two gradient bounds, together with \eqref{eq:explicit-V-fourth} and
\eqref{eq:explicit-J-fourth} therefore give
\begin{equation}
Var(V)\leq16a^{2}m_{4}^{2}K_{\Gamma}^{2}z_{4}^{2}/T,\qquad Var(J)\leq
4C_{J}^{2}K_{\Gamma}^{2}z_{4}^{2}/T. \label{eq:explicit-VJ-variance}%
\end{equation}
In particular $Var(J)=O(T^{-1})$. For the mean, recall $\tau_{3,T}%
=T^{-1}Tr(\Gamma^{3})$. For deterministic symmetric $B_{1},B_{2},B_{3}$,
diagonalize $\Gamma$ and write $N=T^{-1}\sum_{t=1}^{T}\nu_{t,T}(z_{t}%
z_{t}^{\prime}-I_{n})$, where the $z_{t}$ are independent $N(0,I_{n})$
vectors. Independence and centering imply that only equal temporal indices
survive in the third product. Wick's formula (equivalently the third joint
cumulant of Gaussian quadratic forms) therefore gives the \emph{exact}
identity
\begin{equation}
\mathbb{E}\prod_{r=1}^{3}Tr(B_{r}N)=\frac{4\tau_{3,T}}{T^{2}}Tr(B_{1}%
B_{2}B_{3}+B_{1}B_{3}B_{2}). \label{eq:C-third-cumulant}%
\end{equation}
For deterministic $x,y$ and symmetric $A$, let $B=xx^{\prime}$ and
$C_{k}=(ye_{k}^{\prime}+e_{k}y^{\prime})/2$. Since $y^{\prime}NANy=\sum
_{k,l}A_{kl}Tr(C_{k}N)Tr(C_{l}N)$, equation~\eqref{eq:C-third-cumulant}
yields
\begin{align}
\mathbb{E}\{(x^{\prime}Nx)(y^{\prime}NANy)\}  &  =\frac{8\tau_{3,T}}{T^{2}%
}x^{\prime}\left(  \sum_{k,l}A_{kl}C_{k}C_{l}\right)  x\nonumber\\
&  =\frac{2\tau_{3,T}}{T^{2}}\{Tr(A)(x^{\prime}y)^{2}+\Vert y\Vert
^{2}x^{\prime}Ax+2(x^{\prime}y)(x^{\prime}Ay)\}. \label{eq:C-contraction}%
\end{align}
Here the matrix product in the first line is exactly
\[
\sum_{k,l}A_{kl}C_{k}C_{l}=\frac{1}{4}\{y(Ay)^{\prime}+(Ay)y^{\prime
}+Tr(A)yy^{\prime}+\Vert y\Vert^{2}A\}.
\]
Let $w_{i}=We_{i}$, $A=W^{2}=I_{n}+E$, and recall $\Vert w_{i}\Vert^{2}%
=A_{ii}=1$. The exact trace identity is
\[
J=\sum_{i}\{(w_{i}^{\prime}Nw_{i})(w_{i}^{\prime}NANw_{i})-(w_{i}^{\prime
}Nw_{i})(e_{i}^{\prime}N^{2}e_{i})\}.
\]
Apply~\eqref{eq:C-contraction} to the first summand with $x=y=w_{i}$ and
matrix $A$, and to the second with $x=w_{i}$, $y=e_{i}$, and matrix $I_{n}$.
The respective expectations are
\[
\frac{2\tau_{3,T}}{T^{2}}\{n+3(A^{2})_{ii}\},\qquad\frac{2\tau_{3,T}}{T^{2}%
}\{(n+2)W_{ii}^{2}+1\}.
\]
Thus
\begin{equation}
\mathbb{E}J=\frac{2\tau_{3,T}}{T^{2}}\{n^{2}+3Tr(A^{2})-(n+2)\sum_{i}%
W_{ii}^{2}-n\}. \label{eq:C-EJ-raw}%
\end{equation}
Since $Tr(E)=0$, $Tr(A^{2})=n+\Vert E\Vert_{F}^{2}$. Also $W=I_{n}+\Theta$ and
$\operatorname{diag}(W^{2})=\mathbf{1}$ imply $2\Theta_{ii}+(\Theta^{2}%
)_{ii}=0$, hence $W_{ii}^{2}=1-\sum_{j\neq i}\Theta_{ij}^{2}$. Substitution
into~\eqref{eq:C-EJ-raw} gives
\begin{equation}
\mathbb{E}J=\frac{2\tau_{3,T}}{T^{2}}\left\{  3\Vert E\Vert_{F}^{2}%
+(n+2)\sum_{i\neq j}\Theta_{ij}^{2}\right\}  . \label{eq:C-EJ}%
\end{equation}
Because $\Vert E\Vert_{F}=O(1)$, $\Vert\Theta\Vert_{F}\leq\Vert E\Vert
_{F}=O(1)$, $n\asymp T$, and $\tau_{3,T}\leq\Vert\Gamma\Vert_{op}^{2}=O(1)$,
we obtain the explicit bound
\[
|\mathbb{E}J|\leq CT^{-2}\{\Vert E\Vert_{F}^{2}+n\Vert\Theta\Vert_{F}%
^{2}\}=O(T^{-1}).
\]
Combining this with the preceding Gaussian Poincar\'{e} estimate
$Var(J)=O(T^{-1})$ proves $\Vert J\Vert_{2}=O(T^{-1/2})$.

The diagonal correction $\sum_{i}G_{ii}(G_{ii}^{2}-N_{ii}^{2})$ has $L^{2}$
norm $O(T^{-1/2})$ by $\|G_{ii}\|_{6}=O(T^{-1/2})$, $\|N_{ii}\|_{6}%
=O(T^{-1/2})$, and Minkowski's inequality. Therefore $\|C_{3,1}\|_{2}%
=O(T^{-1/2})$.

For the Taylor remainder, write $f=\|\Theta\|_{F}=O(1)$ and use the row
estimates
\begin{equation}
\|\|N_{i\cdot}\|_{2}\|_{8}+ \|\|G_{i\cdot}\|_{2}\|_{8}\le C, \qquad
\|\|H_{i\cdot}\|_{2}\|_{8} \le C\left(  f_{i}+\frac{f}{\sqrt T}\right)  .
\label{eq:R-row}%
\end{equation}
To prove \eqref{eq:R-row}, use the exact covariance formula
\[
\mathbb{E}(N_{ab}N_{cd})=\frac{\tau_{2,T}}{T} (\delta_{ac}\delta_{bd}%
+\delta_{ad}\delta_{bc}).
\]
For a deterministic symmetric matrix $\Theta$ this gives
\[
\mathbb{E}\|e_{i}^{\prime}N\|_{2}^{2}=\frac{\tau_{2,T}}{T}(n+1),
\quad\mathbb{E}\|e_{i}^{\prime}\Theta N\|_{2}^{2}=\frac{\tau_{2,T}}%
{T}(n+1)f_{i}^{2}, \quad\mathbb{E}\|e_{i}^{\prime}N\Theta\|_{2}^{2}=\frac
{\tau_{2,T}}{T}(f^{2}+\|\Theta e_{i}\|_{2}^{2}).
\]
Since $\Theta$ is symmetric, $\|\Theta e_{i}\|_{2}=f_{i}$. Moreover,
$e_{i}^{\prime}\Theta N\Theta=(e_{i}^{\prime}\Theta N)\Theta$, so its
Euclidean norm is bounded by $\|\Theta\|_{op}\|e_{i}^{\prime}\Theta N\|_{2}$.
Thus the $L^{2}$ row estimates follow from $H=\Theta N+N\Theta+\Theta N\Theta$
and $\|\Theta\|_{op}=O(1)$. Since the entries of $N$, $H$, and $G$ are
centered Gaussian quadratic forms, the dimension-free chaos bound in
Section~\ref{app:aux-inequalities}, with $p=8$, upgrades these $L^{2}$ row
estimates to $L^{8}$ bounds with constants independent of $n$ and $T$. This
proves the $N$ and $H$ bounds; the $G$ bound follows from $G=N+H$ and
$\sup_{i} f_{i}\le f=O(1)$. Also, by the moment bounds stated at the outset,
\begin{equation}
\sup_{i}\|G_{ii}\|_{8}\le CT^{-1/2}. \label{eq:R-Gii}%
\end{equation}

By Cauchy--Schwarz in the row index,
\[
\sum_{j\neq i}|P_{ij}|\leq\Vert H_{i\cdot}\Vert_{2}(\Vert G_{i\cdot}\Vert
_{2}+\Vert N_{i\cdot}\Vert_{2}).
\]
H\"{o}lder's inequality with exponents $4,8,8$, followed by
\eqref{eq:R-row}--\eqref{eq:R-Gii}, yields
\begin{align}
\Vert\mathcal{R}_{3}\mathbf{1}_{\mathcal{E}}\Vert_{2}  &  \leq C\sum_{i}\Vert
G_{ii}\Vert_{8}^{2}\Vert\Vert H_{i\cdot}\Vert_{2}\Vert_{8}\left(  \Vert\Vert
G_{i\cdot}\Vert_{2}\Vert_{8}+\Vert\Vert N_{i\cdot}\Vert_{2}\Vert_{8}\right)
\nonumber\\
&  \leq\frac{C}{T}\sum_{i}\left(  f_{i}+\frac{f}{\sqrt{T}}\right)  \leq
\frac{C}{T}\left(  \sqrt{n}f+\frac{nf}{\sqrt{T}}\right)  =O(T^{-1/2}).
\label{eq:R-remainder}%
\end{align}
Combining \eqref{eq:R-remainder} with $\|C_{3,1}\|_{2}=O(T^{-1/2})$ gives
\[
\|C_{3}\|_{2}=O(T^{-1/2}),
\]
and hence $\|T_{3}\|_{2}\leq\|Q_{3}\|_{2}+\|C_{3}\|_{2}=O(T^{-1/2})$.

\noindent Term $T_{4}$. On $\mathcal{E}$, Taylor's formula gives
\[
\frac1{a_{i} a_{j}}-\frac1{b_{i} b_{j}} =-(H_{ii}+H_{jj})+\omega^{\prime}%
_{ij},\qquad|\omega^{\prime}_{ij}|\le C(|H_{ii}|+|H_{jj}|) (|G_{ii}%
|+|G_{jj}|+|N_{ii}|+|N_{jj}|).
\]
Put $m_{i}=\sum_{j\ne i}N_{ij}^{2}$ and $\bar m=(n-1)\tau_{2,T}/T$. The
first-order part is $-2\sum_{i}H_{ii}m_{i}$. Conditional on row $Z_{i}=z$, the
variables $N_{ij}=T^{-1}z^{\prime}\Gamma Z_{j}$ for $j\ne i$ are independent
centered normals with common variance $\varsigma_{i}=T^{-2}z^{\prime}%
\Gamma^{2}z$. Hence
\[
m_{i}\mid Z_{i}\ \overset{d}{=}\varsigma_{i}\chi^{2}_{n-1},\qquad
\mathbb{E}\varsigma_{i}=\frac{\tau_{2,T}}{T},\qquad\|\varsigma_{i}%
-\mathbb{E}\varsigma_{i}\|_{4}=O(T^{-3/2}).
\]
The last bound follows from the fourth moment of a centered Gaussian quadratic
form and $Tr(\Gamma^{2k})=O(T)$ for fixed $k$. Using conditional fourth
moments of the chi-square variable and H\"older,
\[
\|m_{i}-\bar m\|_{4} \le C\left(  \frac{\sqrt n}{T}+\frac{n}{T^{3/2}}\right)
=O(T^{-1/2}),\qquad\|m_{i}\|_{4}=O(1).
\]
Because $\sum_{i}H_{ii}=Tr((2\Theta+\Theta^{2})N)=Tr(EN)$,
\[
\left\|  \bar m\sum_{i}H_{ii}\right\|  _{2} \le C\sqrt{\frac{\tau_{2,T}}{T}%
}\|E\|_{F}=O(T^{-1/2}).
\]
For the centered row sums, Minkowski and H\"older yield
\[
\left\|  \sum_{i}H_{ii}(m_{i}-\bar m)\right\|  _{2} \le C T^{-1}\sum_{i}f_{i}
\le C\frac{\sqrt n}{T}\|\Theta\|_{F}=O(T^{-1/2}).
\]
For the Taylor remainder, set $u_{i}=|G_{ii}|+|N_{ii}|$ and $h_{i}=|H_{ii}|$.
On $\mathcal{E}$, the two-variable Taylor formula, applied to
$g(x,y)=[(1+x)(1+y)]^{-1}$, gives
\[
|\omega^{\prime}_{ij}|\le C(h_{i}+h_{j})(u_{i}+u_{j}).
\]
This bound holds because the gradient of $g$ is Lipschitz on $[-1/2,1]^{2}$,
and $G_{ii}-N_{ii}=H_{ii}$; this is a difference of Taylor expansions at
$(G_{ii},G_{jj})$ and $(N_{ii},N_{jj})$. By the moment bounds stated at the
outset, $\|u_{i}\|_{p}\le C_{p}T^{-1/2}$ and $\|h_{i}\|_{p}\le C_{p}%
T^{-1/2}f_{i}$ for each fixed $p$, so H\"older with exponents $(4,8,8)$ gives
the per-pair bound
\[
\|N_{ij}^{2}h_{i} u_{j}\|_{2} \le\|N_{ij}^{2}\|_{4}\,\|h_{i}\|_{8}%
\,\|u_{j}\|_{8} \le C T^{-2}f_{i}.
\]
The same estimate holds for $N_{ij}^{2}h_{i}u_{i}$, and, upon exchanging
$i,j$, for $N_{ij}^{2}h_{j}u_{i}$ and $N_{ij}^{2}h_{j}u_{j}$. By Minkowski's
inequality, and since $n/T$ is bounded,
\begin{align*}
\left\|  \mathbf{1}_{\mathcal{E}}\sum_{i\ne j}N_{ij}^{2}\omega_{ij}^{\prime
}\right\|  _{2}  &  \le CT^{-2}\sum_{i\ne j}(f_{i}+f_{j})\\
&  \le C\frac{n}{T^{2}}\sum_{i=1}^{n}f_{i}\le C\frac{n^{3/2}}{T^{2}}%
\|\Theta\|_{F}=O(T^{-1/2}).
\end{align*}
Consequently $\|T_{4}\|_{2}=O(T^{-1/2})$.

The four contributions have the following $L^{2}$ bounds, all on $\mathcal{E}%
$:
\begin{align*}
\|T_{1}\|_{2}  &  \le CT^{-1/2}\|E\|_{F}^{2},\\
\|T_{2}\|_{2}  &  \le C\{T^{-1/2}+nT^{-3/2}\|E\|_{F}\},\\
\|T_{3}\|_{2}  &  \le C\{T^{-1/2}+T^{-1}\textstyle\sum_{i}f_{i}+T^{-1}%
nf/\sqrt{T}\},\\
\|T_{4}\|_{2}  &  \le C\{T^{-1/2}+T^{-1}\textstyle\sum_{i}f_{i}+nT^{-2}%
\sum_{i}f_{i}\}.
\end{align*}
Here $f_{i}=\|\Theta_{i\cdot}\|_{2}$, $f=\|\Theta\|_{F}$, and $\sum_{i}
f_{i}\le\sqrt n f$. In the local regime $\|E\|_{F}=O(1)$ and $n/T$ bounded
above and away from zero, every displayed term is $O(T^{-1/2})$, so
$\|D_{n}\mathbf{1}_{\mathcal{E}}\|_{2}\le\sum_{j=1}^{4}\|T_{j}\mathbf{1}%
_{\mathcal{E}}\|_{2}=O(T^{-1/2})$. Together with the exceptional-event bound,
$\mathbb{E}(D_{n}^{2})=O(T^{-1})=O(n^{-1})$. All constants depend only on the
bounds for $n/T$, $K_{\Gamma}$, $K_{E}$, and fixed moment orders.
\end{proof}

\begin{remark}
The unit-diagonal restriction on $W^{2}=R_{\ast}$, namely $2\Theta
_{ii}+(\Theta^{2})_{ii}=0$ for each $i$ (whose trace is the identity in
\eqref{eq:F-identity}), is used in the exact evaluation of $\mathbb{E}J$ in
\eqref{eq:C-EJ}; no assumption $\|R_{\ast}-I\|_{op}=o(1)$ is needed. Section~5
of the main paper reports a Monte Carlo check of the size of $D_{n}$.
\end{remark}

\medskip\noindent Auxiliary lemma (coefficientwise Gaussian cumulant bound).
\emph{Let $V$ be a fixed index set with $|V|\leq4$, let $R_{V}$ be a
correlation matrix, and let $g\sim N(0,R_{V}\otimes I_{T})$. For
$j=1,\ldots,r$, let
\[
Q_{j}=T^{-1}\{g^{\prime}(C_{j}\otimes\Gamma)g -\operatorname{Tr}[(C_{j}%
R_{V})\otimes\Gamma]\},
\]
where the $C_{j}$ are fixed symmetric matrices of order $|V|$ with uniformly
bounded entries. Assume $\operatorname{Tr}\Gamma=T$ and $\|\Gamma\|_{op}\leq
K$. For every fixed $r\geq2$, the joint cumulant $\kappa_{r}(Q_{1}%
,\ldots,Q_{r})$ is a polynomial in the off-diagonal entries of $R_{V}$ of
degree at most $r$, and every coefficient of this polynomial is $O(T^{1-r})$,
uniformly in $R_{V}$, $T$, and the choice of indices in $V$. More generally,
each coefficient of a joint moment of $d$ such centered forms is $O(T^{-\lceil
d/2\rceil})$; the same bound holds for a covariance of two nonempty products,
with only partitions connecting the products retained.}

\noindent Proof. For a real Gaussian vector, polarization of the
quadratic-form cumulant identity yields
\[
\kappa_{r}(Q_{1},\ldots,Q_{r}) =\frac{2^{r-1}}{r}\,T^{-r} \sum_{\pi
\in\mathfrak{S}_{r}} \operatorname{Tr}\!\left(  (C_{\pi(1)}R_{V})\cdots
(C_{\pi(r)}R_{V})\right)  \operatorname{Tr}(\Gamma^{r}).
\]
(The factor $1/r$ accounts for the $r$ cyclic rotations; in particular, the
formula gives $2T^{-2}\operatorname{Tr}(C_{1}R_{V}C_{2}R_{V})
\operatorname{Tr}(\Gamma^{2})$ when $r=2$.) Each trace in dimension $|V|\leq4$
is a sum of at most a constant (depending on $r$) number of products of $r$
entries of $R_{V}$ and entries of the $C_{j}$. On substituting $R_{V}%
=I_{V}+E_{V}$, this is a polynomial of degree at most $r$ in the off-diagonal
correlations, with each coefficient bounded independently of $T$ and $V$.
Since $|\operatorname{Tr}(\Gamma^{r})|\leq TK^{r-1}$, every coefficient has
absolute value at most $C_{r}T^{1-r}$.

For the moment assertion, expand by the moment--cumulant formula. Singleton
blocks vanish because the $Q_{j}$ are centered. A partition of $d$ factors
into $b$ blocks, each of size at least two, has coefficientwise order
$\prod_{B}T^{1-|B|}=T^{b-d}$, and $b\leq\lfloor d/2\rfloor$. The number of
partitions and the polynomial degrees depend only on $d$, so products and sums
retain the stated coefficientwise order. For a covariance of two products,
partitions with no block meeting both products cancel against the product of
expectations. This leaves the same bound for the remaining connecting
partitions. \hfill$\square$ \medskip

\noindent Auxiliary lemma (Eulerian correlation-monomial summation). \emph{Let
$E=E^{\prime}$ have zero diagonal and $\|E\|_{F}\le M$. Let $\mu$ be a
monomial in off-diagonal entries $E_{ab}$ of fixed positive degree $m$, whose
associated multigraph has even degree at every incident vertex. If $\ell$
additional labels do not occur in $\mu$, then the sum of $|\mu|$ over all
assignments of its incident vertices and those $\ell$ labels is bounded by
\[
n^{\ell}\prod_{j=1}^{q}\|E\|_{F}^{m_{j}} = n^{\ell}\|E\|_{F}^{m}, \qquad
m_{1}+\cdots+m_{q}=m,
\]
where $m_{j}$ are the edge counts of its connected components. In particular,
a coefficient bounded by $CT^{-\ell-1}$ contributes at most $CM^{m} n^{\ell
}T^{-\ell-1}=O(n^{-1})$ when $n/T$ is bounded above and below.}

\noindent Proof. Every connected component with nonzero degree has an Eulerian
closed walk, since all degrees are even. Summing the absolute value of the
product of edge weights over its vertex assignments is no larger than the sum
over all closed walks of length $m_{j}$ in the nonnegative symmetric matrix
$|E|=(|E_{ab}|)$, namely $\operatorname{Tr}(|E|^{m_{j}})$. As $|E|$ is
symmetric and $m_{j}\ge2$,
\[
0\le\operatorname{Tr}(|E|^{m_{j}}) \le\sum_{k}|\lambda_{k}(|E|)|^{m_{j}}
\le\|E\|_{F}^{m_{j}}.
\]
The label assignments of distinct components can be summed independently to
obtain an upper bound, even if the original configuration requires distinct
labels. Each absent vertex adds at most $n$ choices. Multiplication gives the
asserted inequality. The zero diagonal makes assignments collapsing an edge
vanish. \hfill$\square$ \medskip

\begin{lemma}
[Local covariance expansion]\label{lem:var-local}Under Assumptions 1--4,
suppose $n(\phi_{2}-1)=O(1)$. Then $Var_{A}(Tr(R^{2}))-Var_{0}(Tr(R^{2}%
))=O(n^{-1})$, where $Var_{0}$ denotes the variance under $H_{0}$.
\end{lemma}

\begin{proof}
With $w_{a}=We_{a}$, write $v_{a}=Z^{\prime}w_{a}$, so that $v_{a}\sim
N(0,I_{T})$, $Cov(v_{a},v_{b})=(R_{\ast})_{ab}I_{T}$, and $r_{ab}%
=v_{a}^{\prime}\Gamma v_{b}/(v_{a}^{\prime}\Gamma v_{a}\,v_{b}^{\prime}\Gamma
v_{b})^{1/2}$. Put $\rho_{ab}=E_{ab}$ for $a\neq b$, so that $\sum_{a\neq
b}\rho_{ab}^{2}=\Vert E\Vert_{F}^{2}=n(\phi_{2}-1)=O(1)$. All constants below
depend only on fixed moment orders, the uniform bound on $\|\Gamma\|_{op}$,
and upper and lower bounds for $T/n$; they do not depend on the individual
correlations or on $n,T$.

\noindent Step 1: Locality. Let $\mathcal{I}$ be the set of ordered pairs
$e=(a,b)$ with $a\neq b$. Then $Tr(R^{2})-n=\sum_{e\in\mathcal{I}}r_{e}^{2}$
and
\[
Var_{A}(Tr(R^{2}))=\sum_{e,e^{\prime}\in\mathcal{I}}c_{ee^{\prime}}(R_{\ast
}),\qquad c_{ee^{\prime}}(R)=Cov_{R}(r_{e}^{2},r_{e^{\prime}}^{2}).
\]
Since $r_{e}^{2}$ depends only on the two vectors $v_{a}$, $a\in e$,
$c_{ee^{\prime}}(R)$ depends on $R$ only through the principal submatrix
$R_{V}$, $V=e\cup e^{\prime}$, and $|V|\in\{2,3,4\}$. The same identity with
$R=I_{n}$ gives $Var_{0}(Tr(R^{2}))$. Ordered pairs are retained throughout.

\noindent Step 2: Polynomial truncation. For $e=(a,b)$ put $s_{e}=T^{-1}%
v_{a}^{\prime}\Gamma v_{b}$, $\delta_{a}=T^{-1}v_{a}^{\prime}\Gamma v_{a}-1$,
define the degree-five Taylor polynomial by
\[
\Pi_{5}(x,y)=\sum_{k=0}^{5}(-1)^{k}\sum_{j=0}^{k}x^{j} y^{k-j}, \qquad\frac
{1}{(1+x)(1+y)}=\Pi_{5}(x,y)+Q_{6}(x,y),
\]
where, for $|x|\vee|y|\le1/2$,
\[
|Q_{6}(x,y)|\le C(|x|+|y|)^{6}.
\]
Put $\Pi_{e}=\Pi_{5}(\delta_{a},\delta_{b})$ and set $\mathcal{P}_{e}%
=s_{e}^{2}\Pi_{e}$ and $\mathcal{R}_{e}=r_{e}^{2}-\mathcal{P}_{e}$. On the
event $|\delta_{a}|\vee|\delta_{b}|\leq1/2$ the multivariate Taylor formula
gives $|\mathcal{R}_{e}|\leq C|s_{e}|^{2}(|\delta_{a}|+|\delta_{b}|)^{6}$. For
every fixed $q$, Gaussian quadratic-form moment bounds imply $\Vert\delta
_{a}\Vert_{q}\leq C_{q}T^{-1/2}$ and $\Vert s_{e}\Vert_{q}\leq C_{q}(|\rho
_{e}|+T^{-1/2})$. H\"{o}lder's inequality therefore yields the first $L^{2}$
bound displayed below. The complement has probability at most $Ce^{-cT}$
uniformly in the alternative, since each $v_{a}$ is marginally standard
Gaussian; as $0\leq r_{e}^{2}\leq1$ and $\mathcal{P}_{e}$ has uniformly
bounded moments, its contribution is exponentially small. Thus
\[
\Vert\mathcal{R}_{e}\Vert_{2}\leq C(\rho_{e}^{2}+T^{-1})T^{-3},\qquad\Vert
r_{e}^{2}-\mathbb{E}r_{e}^{2}\Vert_{2}+\Vert\mathcal{P}_{e}-\mathbb{E}%
\mathcal{P}_{e}\Vert_{2}\leq C(|\rho_{e}|T^{-1/2}+T^{-1}),
\]
where $\rho_{e}=\rho_{ab}$. Writing $r_{e}^{2}=\mathcal{P}_{e}+\mathcal{R}%
_{e}$, Cauchy--Schwarz and symmetry in $(e,e^{\prime})$ give
\[
\sum_{e,e^{\prime}}\left\vert c_{ee^{\prime}}(R_{\ast})-Cov(\mathcal{P}%
_{e},\mathcal{P}_{e^{\prime}})\right\vert \leq CT^{-3}\sum_{e}(\rho_{e}%
^{2}+T^{-1})\sum_{e^{\prime}}(|\rho_{e^{\prime}}|T^{-1/2}+T^{-1}).
\]
For completeness, the preceding covariance approximation is uniform in all
pair configurations, including $e=e^{\prime}$, pairs sharing one index, and
disjoint pairs. Indeed, with $\mathcal{E}_{e}=\{\max(|\delta_{a}|,|\delta
_{b}|)\le1/2\}$, the fixed-order Gaussian quadratic-form moment bounds imply,
uniformly in $a,b$ and the admissible $R_{*}$,
\[
\|s_{e}\|_{8}\le C(|\rho_{e}|+T^{-1/2}),\qquad\|\delta_{a}\|_{24}+\|\delta
_{b}\|_{24}\le CT^{-1/2}.
\]
Consequently, H\"older's inequality on $\mathcal{E}_{e}$ yields
\[
\|\mathcal{R}_{e}\mathbf{1}_{\mathcal{E}_{e}}\|_{2} \le C\|s_{e}\|_{8}^{2}
(\|\delta_{a}\|_{24}+\|\delta_{b}\|_{24})^{6} \le C(\rho_{e}^{2}+T^{-1}%
)T^{-3}.
\]
On $\mathcal{E}_{e}^{c}$, Gaussian concentration for the normalized quadratic
form $T^{-1}v_{a}^{\prime}\Gamma v_{a}$ gives $\mathbb{P}(\mathcal{E}_{e}%
^{c})\le Ce^{-cT}$; the marginal law of each $v_{a}$ is standard Gaussian, so
this estimate is independent of the cross-sectional correlation. The
inequality $0\le r_{e}^{2}\le1$ and uniform fixed-order moments of the
polynomial $\mathcal{P}_{e}$ give an $O(e^{-c^{\prime}T})$ contribution to the
$L^{2}$ norm by H\"older's inequality (with possibly smaller $c^{\prime}>0$).

To check the centered-factor estimate, write $s_{e}=\rho_{e}+\widetilde{s}%
_{e}$, where $\|\widetilde{s}_{e}\|_{q}\le C_{q}T^{-1/2}$, and note that
$\|\Pi_{e}-1\|_{q}\le C_{q}T^{-1/2}$ for fixed $q$. Then
\[
\mathcal{P}_{e}-\rho_{e}^{2} =2\rho_{e}\widetilde{s}_{e}+\widetilde{s}_{e}^{2}
+s_{e}^{2}(\Pi_{e}-1),
\]
so centering and H\"older's inequality yield $\|\mathcal{P}_{e}-\mathbb{E}%
\mathcal{P}_{e}\|_{2} \le C(|\rho_{e}|T^{-1/2}+T^{-1})$. The same estimate for
$r_{e}^{2}$ follows by adding $\mathcal{R}_{e}$; $|\rho_{e}|\le1$ ensures that
its remainder is no larger in order. Finally, Cauchy--Schwarz applied to
$\operatorname{Cov}(\mathcal{R}_{e},\mathcal{P}_{e^{\prime}})$,
$\operatorname{Cov}(\mathcal{P}_{e},\mathcal{R}_{e^{\prime}})$ and
$\operatorname{Cov}(\mathcal{R}_{e},\mathcal{R}_{e^{\prime}})$ produces the
displayed double-sum bound. In particular, this argument does not assume
independence between $e$ and $e^{\prime}$.

The first sum is $\Vert E\Vert_{F}^{2}+n(n-1)/T=O(n)$, and the second is
$O(n)$ because $\sum_{e}|\rho_{e}|\leq n\Vert E\Vert_{F}$. The right-hand side
is therefore $O(n^{2}T^{-3})=O(n^{-1})$. The same bound holds under $H_{0}$.
Hence
\begin{equation}
Var_{A}(Tr(R^{2}))-Var_{0}(Tr(R^{2}))=\sum_{e,e^{\prime}}\{\vartheta
_{ee^{\prime}}(R_{\ast})-\vartheta_{ee^{\prime}}(I_{n})\}+O(n^{-1}),
\label{eq:local-reduction}%
\end{equation}
where $\vartheta_{ee^{\prime}}(R)=Cov_{R}(\mathcal{P}_{e},\mathcal{P}%
_{e^{\prime}})$.

\noindent Step 3: Polynomial coefficients and the free-vertex bound. For
$V=e\cup e^{\prime}$ write $R_{V}=I_{V}+(\rho_{ab})_{a,b\in V}$, with zero
diagonal in the perturbation. Each $\delta_{a}$ and $\widetilde{s}_{ab}%
=s_{ab}-\rho_{ab}$ is a centered Gaussian quadratic form. For fixed $r\ge2$,
the auxiliary lemma above applies to their joint cumulants and to the moments
and covariances of products of these centered forms. In particular, each
coefficient in the finite-dimensional correlation polynomial of an $r$th
cumulant is $O(T^{1-r})$, and a connected partition of $d$ centered factors
into $b$ blocks contributes coefficientwise $O(T^{b-d})$.

To account explicitly for powers of $\rho_{e}$, write
\[
\mathcal{P}_{e}=\sum_{\alpha+\beta\le5}c_{\alpha\beta} (\rho_{e}^{2}+2\rho
_{e}\widetilde{s}_{e}+\widetilde{s}_{e}^{2}) \delta_{a}^{\alpha}\delta
_{b}^{\beta}.
\]
Every summand is a deterministic correlation monomial multiplied by a product
of centered quadratic forms. Terms with no centered factors are deterministic
and disappear from covariances. More generally, a nonzero covariance
contribution requires at least one centered factor from each of the two
pairwise polynomials. The explicit factors $\rho_{e}$ and $\rho_{e^{\prime}}$
are counted as edges of the correlation monomial, not as centered factors.
Consequently, their presence never improves the $T$-order of a coefficient; it
can only reduce the number of absent vertices. The preceding
connected-partition bound applies term by term to all three choices $\rho
_{e}^{2}$, $2\rho_{e}\widetilde{s}_{e}$, and $\widetilde{s}_{e}^{2}$ and
likewise for $e^{\prime}$.

It follows that $\vartheta_{ee^{\prime}}(R_{V})-\vartheta_{ee^{\prime}}%
(I_{V})$ is a polynomial $\sum_{\mu\ne1}a_{\mu}\mu$ in the off-diagonal
correlations, of degree bounded independently of $n,T$. Let $\ell(\mu)$ count
the vertices of $V$ absent from $\mu$. We claim
\begin{equation}
\label{eq:local-coeff-bound}\deg_{a}(\mu)\ \text{is even for every vertex
}a,\qquad|a_{\mu}|\le CT^{-\ell(\mu)-1}.
\end{equation}
The parity follows from simultaneous sign changes of $v_{a}$ and the
correlations incident to $a$: each $\mathcal{P}_{e}$ is invariant. For
$\ell=0$, a nonzero covariance term contains at least two centered factors,
and the connected-partition bound is $O(T^{-1})$. For $\ell\ge1$, choose an
absent vertex $a$ and an edge $e$ containing it. Neither an explicit factor
$\rho_{e}$ nor a correlation factor incident to $a$ can occur in the
coefficient of $\mu$. Therefore the contribution from $\mathcal{P}_{e}$
necessarily uses $\widetilde{s}_{e}^{2}$ and at least two centered factors.
The other pairwise polynomial must contribute at least one centered factor for
a nonzero covariance. Thus $d\ge3$ and the connected-partition bound is
$O(T^{-2})$, proving the claim for $\ell=1$.

For $\ell=2$, a nonconstant monomial with even degree at each vertex has at
least two vertices and hence $|V|=4$. The pairs are disjoint. If the two
absent vertices form one of the pairs, say $e$, then $\mu$ involves only
$\rho_{e^{\prime}}$; setting every other correlation to zero makes
$(v_{a})_{a\in e}$ independent of $(v_{a})_{a\in e^{\prime}}$, so
$\vartheta_{ee^{\prime}}$ vanishes identically and all such coefficients are
zero. Otherwise, after relabeling, the only possible monomial uses
correlations between one endpoint $i$ of $e=(i,j)$ and one endpoint $k$ of
$e^{\prime}=(k,l)$, with $j,l$ absent. To bound all its coefficients
simultaneously, set $\rho_{ik}=t$ and every other off-diagonal correlation to
zero. This is a valid correlation matrix for $|t|\le1/2$. Conditional on
$(v_{i},v_{k})$, $v_{j},v_{l}$ are independent, so
\[
\vartheta_{ee^{\prime}}(t)=Cov\{h(v_{i}),h^{\prime}(v_{k})\},\qquad
h(v_{i})=\mathbb{E}(\mathcal{P}_{ij}\mid v_{i}).
\]
Since $\rho_{ij}=0$, expand $\Pi_{ij}=\sum_{\alpha,\beta}c_{\alpha\beta}
\delta_{i}^{\alpha}\delta_{j}^{\beta}$ and write
\[
h(v_{i})=\sum_{\alpha,\beta}c_{\alpha\beta}\delta_{i}^{\alpha}\frac
{v_{i}^{\prime}A_{\beta}v_{i}}{T^{2}},\qquad A_{\beta}=\Gamma\mathbb{E}%
(vv^{\prime}\delta^{\beta})\Gamma,
\]
where $v\sim N(0,I_{T})$ and $\delta=T^{-1}v^{\prime}\Gamma v-1$. For every
fixed $\beta$, $\|A_{\beta}\|_{op}\le C_{\beta}$ uniformly: for a unit vector
$u$, Cauchy--Schwarz bounds $\mathbb{E}\{(u^{\prime}v)^{2}|\delta|^{\beta}\}$
by a constant depending only on $\beta$ and $K_{\Gamma}$, so $\|\mathbb{E}%
(vv^{\prime}\delta^{\beta})\|_{op}\le C_{\beta}$, and $\|A_{\beta}\|_{op}%
\le\|\Gamma\|_{op}^{2}\|\mathbb{E}(vv^{\prime}\delta^{\beta})\|_{op}$. Set
$\varpi_{\beta}=T^{-1}Tr(A_{\beta})$ and $B_{\beta}=A_{\beta}-\varpi_{\beta
}\Gamma$. Then $Tr(B_{\beta})=0$, $\|B_{\beta}\|_{F}\le C\sqrt T$, and
\[
\frac{v_{i}^{\prime}A_{\beta}v_{i}}{T^{2}} =\frac{\varpi_{\beta}}{T}%
(1+\delta_{i}) +\frac{v_{i}^{\prime}B_{\beta}v_{i}}{T^{2}},\qquad\left\|
\frac{v_{i}^{\prime}B_{\beta}v_{i}}{T^{2}}\right\|  _{q} \le C_{q}T^{-3/2}.
\]
Here the last estimate follows from Gaussian quadratic-form moments, and
$\|\delta_{i}\|_{q}\le C_{q}T^{-1/2}$ for each fixed $q$. To make the
centering argument explicit, set $U_{\beta,i}=T^{-2}v_{i}^{\prime}B_{\beta
}v_{i}$ and $H_{\alpha\beta,i}=\delta_{i}^{\alpha}\{\varpi_{\beta}%
T^{-1}(1+\delta_{i})+U_{\beta,i}\}$. The deterministic coefficients satisfy
$|\varpi_{\beta}|\le C_{\beta}$. For $\alpha=0$, subtracting the expectation
gives
\[
\|H_{0\beta,i}-\mathbb{E}H_{0\beta,i}\|_{2} \le CT^{-1}\|\delta_{i}%
\|_{2}+2\|U_{\beta,i}\|_{2} \le CT^{-3/2}.
\]
For every fixed $\alpha\ge1$, H\"older's inequality and the fixed-order
Gaussian moment estimates give
\begin{align*}
\|H_{\alpha\beta,i}-\mathbb{E}H_{\alpha\beta,i}\|_{2}  &  \le2\|H_{\alpha
\beta,i}\|_{2}\\
&  \le CT^{-1}\|\delta_{i}\|_{4\alpha}^{\alpha} (1+\|\delta_{i}\|_{4})
+C\|\delta_{i}\|_{4\alpha}^{\alpha}\|U_{\beta,i}\|_{4}\\
&  \le CT^{-1-\alpha/2}+CT^{-(3+\alpha)/2} \le CT^{-3/2}.
\end{align*}
Here the case $\alpha=0$ is treated separately so that the nonvanishing
deterministic part of $T^{-1}(1+\delta_{i})$ is not mistaken for a centered
fluctuation. Since the Taylor polynomial contains only finitely many
$(\alpha,\beta)$, summing these bounds shows $\|h-\mathbb{E}h\|_{2}\le
CT^{-3/2}$; the same argument applies to $h^{\prime}$. The marginal bounds are
uniform in $t$, because each $v_{i}$ has the same $N(0,I_{T})$ law for
$|t|\le1/2$. Cauchy--Schwarz therefore gives $|\vartheta_{ee^{\prime}}(t)|\le
CT^{-3}$ uniformly on $[-1/2,1/2]$. Because this is a polynomial in $t$ of
fixed degree, its individual coefficients are $O(T^{-3})$ by
finite-dimensional polynomial norm equivalence. To make the coefficient
identification explicit, fix the ordered-pair configuration $e=(i,j)$,
$e^{\prime}=(k,l)$ with four distinct indices. A monomial in which $j$ and $l$
are absent can involve only the single off-diagonal variable $\rho_{ik}$: all
other correlations have an endpoint in $\{j,l\}$. Thus the specialization
$\rho_{ik}=t$, with all other off-diagonal entries zero, preserves every
coefficient in this class and eliminates every monomial outside it. In
particular, no cancellation with monomials involving an absent vertex is
possible. For the fixed Taylor degree, the univariate polynomial has uniformly
bounded degree; evaluating it at a fixed set of distinct points in
$[-1/2,1/2]$ and inverting the corresponding Vandermonde matrix bounds each
coefficient by $CT^{-3}$. This establishes the claimed bound for the
two-absent-vertex class. No nonconstant even-degree monomial has $\ell>2$ when
$|V|\le4$. This completes the coefficient proof.

\noindent Step 4: Summation. We apply the Eulerian correlation-monomial
summation lemma above. Fix one of the finitely many types of pairs
$(e,e^{\prime})$, a nonconstant monomial $\mu$, and its factor multigraph
$\mathcal{G}_{\mu}$ on the vertices of $\mu$. By the parity property in
\eqref{eq:local-coeff-bound}, every connected component of $\mathcal{G}_{\mu}$
is Eulerian. The auxiliary summation lemma bounds the sum of $|\mu|$ over all
assignments of its incident vertices and $\ell(\mu)$ absent vertices by
$n^{\ell(\mu)}\|E\|_{F}^{\deg(\mu)}$. Combining this with
\eqref{eq:local-coeff-bound}, each fixed monomial type contributes at most
\[
C T^{-\ell(\mu)-1} n^{\ell(\mu)} \|E\|_{F}^{\deg(\mu)}=O(n^{-1}).
\]
This estimate allows disconnected factor graphs and repeated labels; no
distinctness assumption is needed for the upper bound. The number of abstract
graph types and monomials is bounded by an absolute constant because $|V|\le4$
and the Taylor degree is fixed; thus this summation does not introduce a
dimension-dependent factor. Summing over these types and monomials proves the lemma.
\end{proof}

\noindent Completion of the proof of Lemma~6. Lemma~\ref{lem:poincare}
establishes the asserted uniform $O(1)$ variance bound. If $n(\phi
_{2}-1)=O(1)$, the covariance comparison in Lemma~\ref{lem:var-local} gives
\[
Var_{A}(Tr(R^{2})) =Var_{0}(Tr(R^{2}))+O(n^{-1}).
\]
Applying the null expansion of Lemma~3 proves (4).

\section{Proof of Theorem~2}

\label{proof:clt_alt}

\noindent Put $F_{n}=Tr(R^{2})$, $F_{0n}=Tr(R_{0}^{2})$, and
$s_{n}=n(\phi_{2}-1)$. Under the local-alternative condition $s_{n}=O(1)$,
Lemma~\ref{lem:coupling} constructs $F_{n}$ and $F_{0n}$ from the same
Gaussian innovations and gives
\[
D_{n}=F_{n}-F_{0n}-s_{n},\qquad\mathbb{E}_{A}(D_{n}^{2})=O(n^{-1}).
\]
Thus $D_{n}=o_{P_{A}}(1)$, while $F_{0n}$ has exactly the null distribution.
Since $\tau_{2,T}\geq1$, $K_{T}\geq0$ (Remark~2), $\tau_{2,T},K_{T}=O(1)$, and
$\lambda\to\lambda_{0}\in(0,\infty)$, there are constants $0<c\leq C$ with
$c\leq\sigma_{n}^{2}\leq C$. With $\Delta_{n}=s_{n}/\sigma_{n}$, the
definition of $LM$ yields the exact identity
\[
LM-\Delta_{n}=\frac{F_{0n}-\mu_{n}}{\sigma_{n}}+\frac{D_{n}}{\sigma_{n}}.
\]
Since $F_{0n}$ has the null distribution, Theorem~1 applies to the first term
on the right, and the second is $o_{P_{A}}(1)$. Slutsky's theorem gives
$LM-\Delta_{n}\overset{d}{\longrightarrow}N(0,1)$; no convergence of the
bounded sequence $\Delta_{n}$ is required.

%=============================================================

\section{Proof of Corollary~1 (Detection Boundary)}

\label{proof:boundary}

Write $s_{n}=n(\phi_{2}-1)=Tr(R_{*}^{2})-n=\|R_{*}-I_{n}\|_{F}^{2}\geq0$ and
$\Delta_{n}=s_{n}/\sigma_{n}$. As in the proof of Theorem~2, there exist
constants $0<c<C<\infty$ such that $c\leq\sigma_{n}\leq C$ for all
sufficiently large $n$. Thus $s_{n}\to0$ (respectively, $s_{n}\to\infty$) if
and only if $\Delta_{n}\to0$ (respectively, $\Delta_{n}\to\infty$). Set
\[
b_{\alpha}(d)=\Phi(d-z_{\alpha/2})+\Phi(-d-z_{\alpha/2}),\qquad d\geq0.
\]
For $0<\alpha<1$, $b_{\alpha}(0)=\alpha$, $b_{\alpha}(d)<1$ for every finite
$d$, and $b_{\alpha}(d)\to1$ as $d\to\infty$. Moreover,
\[
b_{\alpha}^{\prime}(d)=\varphi(d-z_{\alpha/2})-\varphi(d+z_{\alpha/2})>0
\quad(d>0),
\]
where $\varphi$ denotes the standard normal density. Hence $b_{\alpha}$ is
strictly increasing on $[0,\infty)$.

\noindent Local power. Suppose $s_{n}=O(1)$. Let $\mathcal{W}_{n}=LM-\Delta_{n}$
and $z=z_{\alpha/2}$. By Theorem~2, $\mathcal{W}_{n}\overset{d}{\longrightarrow
}N(0,1)$. Since $\Phi$ is continuous, P\'olya's theorem gives $\sup_{x}%
|P_{A}(\mathcal{W}_{n}\leq x)-\Phi(x)|\rightarrow0$, and likewise for
$P_{A}(\mathcal{W}_{n}<x)$; this covers the moving thresholds $z-\Delta_{n}$
and $-z-\Delta_{n}$. Since the two rejection events are disjoint,
\[
P_{A}(|LM|>z)=P_{A}(\mathcal{W}_{n}>z-\Delta_{n})+P_{A}(\mathcal{W}%
_{n}<-z-\Delta_{n}),
\]
and consequently
\begin{equation}
P_{A}(|LM|>z_{\alpha/2})=b_{\alpha}(\Delta_{n})+o(1). \label{eq:local-power}%
\end{equation}
No convergence of the bounded sequence $\Delta_{n}$ is required.

\noindent(a) Sufficiency. Suppose $s_{n}\to\infty$. The mean expansion of
Lemma~5, together with the boundedness of $\phi_{2}-1$ under Assumption~4,
implies
\[
m_{n}:=\mathbb{E}_{A}[Tr(R^{2})]-\mu_{n}=s_{n}+O(1)+O(n^{-1}) =s_{n}%
\{1+o(1)\}.
\]
Lemma~6 supplies $Var_{A}\{Tr(R^{2})\}=O(1)$. Consequently $m_{n}-z_{\alpha
/2}\sigma_{n}\to\infty$, and Chebyshev's inequality yields
\begin{align*}
P_{A}\{|LM|\leq z_{\alpha/2}\}  &  \leq P_{A}\{Tr(R^{2})-\mathbb{E}%
_{A}Tr(R^{2}) \leq-(m_{n}-z_{\alpha/2}\sigma_{n})\}\\
&  \leq\frac{Var_{A}\{Tr(R^{2})\}} {(m_{n}-z_{\alpha/2}\sigma_{n})^{2}%
}\longrightarrow0.
\end{align*}
Thus the rejection probability tends to one. This covers every diverging
signal, including $s_{n}=o(n)$.

\noindent(a) Necessity. If $s_{n}$ does not tend to infinity, there is a
subsequence on which $s_{n}$ is bounded. Passing to a further subsequence,
$s_{n}\to s\in[0,\infty)$ and $\sigma_{n}\to\sigma\in[c,C]$. This is a
local-alternative subsequence, so \eqref{eq:local-power} gives rejection
probability tending to $b_{\alpha}(s/\sigma)<1$. Hence power cannot tend to one along the
full sequence.

\noindent(b) Sufficiency. If $s_{n}\to0$, then $\Delta_{n}\to0$; \eqref{eq:local-power} and
continuity of $b_{\alpha}$ give power tending to $b_{\alpha
}(0)=\alpha$.

\noindent(b) Necessity. Suppose instead that $s_{n}$ does not tend to zero.
There are $\varepsilon>0$ and a subsequence with $s_{n}\geq\varepsilon$. If
this subsequence has an unbounded further subsequence, extract one with
$s_{n}\to\infty$; part (a) then gives power tending to one, contrary to
convergence to $\alpha$. Otherwise $s_{n}$ is bounded on a further
subsequence; compactness permits $s_{n}\to s\geq\varepsilon$ and $\sigma
_{n}\to\sigma\in[c,C]$. Then \eqref{eq:local-power} gives limiting power $b_{\alpha
}(s/\sigma)>\alpha$, again a contradiction. Thus power tending to $\alpha$
forces $s_{n}\to0$.

\noindent(c) If $\Delta_{n}\to\Delta\in(0,\infty)$, then $s_{n}=O(1)$ because
$\sigma_{n}$ is bounded, and \eqref{eq:local-power} yields power tending to
$b_{\alpha}(\Delta)\in(\alpha,1)$.

\begin{remark}
The equivalences above concern the rejection probability of the stated
two-sided LM test at a fixed nominal level $\alpha\in(0,1)$. They rely on the
boundedness and nondegeneracy of $\sigma_{n}$, the uniform $O(1)$ alternative
variance bound in Lemma~6, and the local Gaussian shift approximation of
Theorem~2. Theorem~2 alone does not justify the $s_{n}\to\infty$ case; this is
why part (a) uses Chebyshev's inequality.
\end{remark}

\section{Auxiliary Gaussian and operator inequalities}

\label{app:aux-inequalities} This section proves the analytic inputs used in
the null martingale argument and in the local-alternative coupling. Constants
below depend only on $K_{\Gamma}$ and on fixed moment orders.

\paragraph{Quadratic-form cumulants and moments.}

Let $z\sim N(0,I_{d})$ and $B=B^{\prime}$ be deterministic. Diagonalization
and the moment generating function of a squared standard normal yield, for
$|t|$ sufficiently small,
\[
\log\mathbb{E}\exp\{t(z^{\prime}Bz-Tr\,B)\} =\sum_{r=2}^{\infty}\frac{2^{r-1}%
}r t^{r}Tr(B^{r}).
\]
Thus $\kappa_{r}=2^{r-1}(r-1)!Tr(B^{r})$ for $r\ge2$. For each fixed integer
$m$, the moment--cumulant formula, together with $|Tr(B^{r})|\le\|B\|_{F}%
^{2}\|B\|_{op}^{r-2}$, implies
\[
\mathbb{E}|z^{\prime}Bz-Tr\,B|^{2m}\le C_{m}\|B\|_{F}^{2m}.
\]
For symmetric $B_{1},B_{2},B_{3}$, polarization of the third cumulant gives
\[
\mathbb{E}\prod_{j=1}^{3}(z^{\prime}B_{j}z-Tr\,B_{j}) =4Tr(B_{1}B_{2}%
B_{3}+B_{1}B_{3}B_{2}).
\]
These are exact identities, including when the matrices are singular.

\paragraph{Inverse quadratic forms without a spectral gap.}

Suppose $\Gamma\succeq0$, $Tr\Gamma=T$, and $\|\Gamma\|_{op}\le K_{\Gamma}$.
If exactly $k$ eigenvalues are at least $1/2$, then $T\le kK_{\Gamma}%
+(T-k)/2$, so $k\ge T/(2K_{\Gamma}-1)$. Consequently $z^{\prime}\Gamma
z\ge\frac12\chi_{k}^{2}$ under an eigenbasis coupling, and for every fixed
positive integer $m$ and all sufficiently large $T$,
\[
\mathbb{E}(z^{\prime}\Gamma z)^{-m}\le\mathbb{E}(\tfrac12\chi_{k}^{2})^{-m}
=\frac{\Gamma(k/2-m)}{\Gamma(k/2)}=O(T^{-m}).
\]
This is the bound (B.2); in particular, no bound on the smallest eigenvalue of
$\Gamma$ is needed.

\paragraph{Hilbert-valued Gaussian-chaos moment bound.}

Let $Y=(Y_{1},\ldots,Y_{d})$ be a vector of polynomials of total Gaussian
degree at most two, including deterministic terms. Scalar Gaussian
hypercontractivity gives $\|Y_{j}\|_{L^{p}}\le(p-1)\|Y_{j}\|_{L^{2}}$ for
$p\ge2$. Minkowski's inequality in $L^{p/2}$ then yields
\[
\big\|\|Y\|_{2}\big\|_{L^{p}}^{2} =\big\|\sum_{j}Y_{j}^{2}\big\|_{L^{p/2}}
\le\sum_{j}\|Y_{j}\|_{L^{p}}^{2} \le(p-1)^{2} \mathbb{E}\|Y\|_{2}^{2}.
\]
Thus the constant is independent of the dimension $d$. This yields the rowwise
$L^{8}$ bounds in the coupling proof once their second moments have been
computed from the exact covariance of $N$.

\section*{References}

%References: uniform 2em hanging indentation; no list environment.
\begingroup
\setlength{\leftskip}{2em}\setlength{\parindent}{-2em}\setlength{\parskip}{0.4\baselineskip}\raggedright
Davidson, K. R., \& Szarek, S. J. (2001). Local operator theory, random
matrices and Banach spaces. In W. B. Johnson \& J. Lindenstrauss (Eds.),
\textit{Handbook of the geometry of Banach spaces} (Vol.~1, pp.~317--366). North-Holland.

Hounyo, U., Kao, C., \& Kim, M. S. (2026a). \textit{Limiting spectral
distribution for large sample correlation matrices under temporal dependence}
[Working paper]. University of Connecticut.

Isserlis, L. (1918). On a formula for the product-moment coefficient of any
order of a normal frequency distribution in any number of variables.
\textit{Biometrika}, \textit{12}(1--2), 134--139.
\href{https://doi.org/10.1093/biomet/12.1-2.134}{\mbox{doi:10.1093/biomet/12.1-2.134}}%

Rudelson, M., \& Vershynin, R. (2013). Hanson--Wright inequality and
sub-Gaussian concentration. \textit{Electronic Communications in Probability},
\textit{18}, Article 82.
\href{https://doi.org/10.1214/ECP.v18-2865}{doi:10.1214/ECP.v18-2865}

\endgroup

\end{document}